\documentclass[moor,nonblindrev]{informs1}              

\usepackage{amsmath,amssymb}
\usepackage{empheq}
\usepackage[mathscr]{euscript}
\usepackage{stmaryrd}
\usepackage{dsfont}
\usepackage{graphicx}
\usepackage{textgreek}
\usepackage{scalerel}

\newcommand{\abs}[1]{\left|#1\right|}                       
\newcommand{\bra}[1]{\left(#1\right)}
\newcommand{\set}[1]{\left\{#1\right\}}

\newcommand{\scrO}{\mathscr{O}}

\newcommand{\scrR}{\mathscr{R}}
\newcommand{\scrS}{\mathscr{S}}

\newcommand{\scrV}{\mathscr{V}}

\newcommand{\bbR}{\mathbb{R}}
\newcommand{\bbZ}{\mathbb{Z}}
\newcommand{\bbone}{\mathds{1}}

\newcommand{\CDC}{\operatorname{CDC}}
\newcommand{\Conv}{\operatorname{Conv}}
\renewcommand{\div}{\operatorname{div}}
\newcommand{\ext}{\operatorname{ext}}
\newcommand{\graph}{\operatorname{gr}}

\renewcommand{\int}{\operatorname{relint}}

\newcommand{\Pre}{\operatorname{Pre}}
\newcommand{\Proj}{\operatorname{Proj}}

\newcommand{\relint}{\operatorname{relint}}
\newcommand{\rem}{\operatorname{mod}}

\newcommand\defeq{\mathrel{\overset{\makebox[0pt]{\mbox{\normalfont\tiny\sffamily def}}}{=}}}
\usepackage[normalem]{ulem}

\newcommand{\sidx}[1]{\left\llbracket     #1 \right\rrbracket}
\usepackage{natbib}
 \NatBibNumeric
 \def\bibfont{\small}%
 \bibpunct[, ]{[}{]}{,}{n}{}{,}%

\usepackage[colorlinks=true,breaklinks=true,bookmarks=true,urlcolor=blue,
     citecolor=blue,linkcolor=blue,bookmarksopen=false,draft=false]{hyperref}

\def\EMAIL#1{\href{mailto:#1}{#1}}

\TheoremsNumberedThrough     

\EquationsNumberedThrough    

\begin{document}


\RUNAUTHOR{Huchette and Vielma}

\RUNTITLE{A combinatorial approach for disjunctive constraints}

\TITLE{A combinatorial approach for small and strong formulations of disjunctive constraints}

\ARTICLEAUTHORS{%
\AUTHOR{Joey Huchette}
\AFF{Operations Research Center, MIT, \EMAIL{huchette@mit.edu}}
\AUTHOR{Juan Pablo Vielma}
\AFF{Sloan School of Management, MIT, \EMAIL{jvielma@mit.edu}}
} 

\ABSTRACT{
    We present a framework for constructing strong mixed-integer programming formulations for logical disjunctive constraints. Our approach is a generalization of the logarithmically-sized formulations of Vielma and Nemhauser for SOS2 constraints \cite{Vielma:2009a}, and we offer a complete characterization of its expressive power. We apply the framework to a variety of disjunctive constraints, producing novel small and strong formulations for outer approximations of multilinear terms, generalizations of special ordered sets, piecewise linear functions over a variety of domains, and obstacle avoidance constraints.
}%


\MSCCLASS{90C11}
\ORMSCLASS{Primary: Integer Programming; secondary: Piecewise Linear, Polyhedra}
\HISTORY{}

\maketitle

%


%
%
%

\section{Introduction}
A central modeling primitive in mathematical optimization is the disjunctive constraint: any feasible solution must satisfy at least one of some fixed, finite collection of alternatives. This type of constraint is general enough to capture structures as diverse as boolean satisfiability, complementarity constraints, special ordered sets, and (bounded) integrality. The special case of polyhedral disjunctive constraints corresponds to the form
\begin{equation}\label{disconst}
x \in \bigcup\nolimits_{i=1}^d P^i,
\end{equation}
where we have that each $P^i\subseteq \mathbb{R}^n$ is a polyhedron. In this work, we will focus on $\scrV$-polyhedra; that is, we have a description of the $P^i$ in terms of their extreme points $\ext(P^i)$.

We are particularly interested in the case where constraint \eqref{disconst} is \emph{primitive}, or a basic building block for a much more complex optimization problem. For this reason, we are interested in modeling it in a generic, composable way. In particular, if \eqref{disconst} is embedded in a larger, more complex optimization problem
\begin{equation} \label{eqn:larger-optimization-problem}
    \min_{(x,y) \in Q : \eqref{disconst}} f(x,y),
\end{equation}
we hope for a mathematical description sufficiently structured such that we may use more advanced algorithmic approaches, beyond na\"ive enumeration, to solve \eqref{eqn:larger-optimization-problem}. In particular, $Q$ could be described by any number of different types of constraints: linear inequalities, conic constraints, integrality conditions, or additional disjunctive constraints. In this context, it is well known that merely constructing the convex hull $\Conv(\bigcup_{i=1}^d P^i)$ is not sufficient for solving \eqref{eqn:larger-optimization-problem}; we will need a formulation for $\bigcup_{i=1}^d P^i$ directly.

Mixed-integer programming (MIP) has emerged as an incredibly expressive modeling methodology, with advanced computational methods capable of solving many problems of practical interest, often at very large scale \cite{Bixby:2007,Junger:2010}. Constraint \eqref{disconst} with polyhedral sets $P^i$ is particularly well-suited for a mixed-integer programming approach. Indeed, standard formulations for \eqref{disconst} were presented in \cite{Jeroslow:1984}, and are \emph{ideal}, or as strong as possible with respect to their continuous linear programming relaxations (see Section~\ref{ss:formulation-definitions} for a formal definition). However, this formulation requires introducing $d$ auxiliary binary variables, which may be impractically large, especially in the context of the larger problem~\eqref{eqn:larger-optimization-problem}.

However, it is sometimes possible to construct ideal formulations with considerably fewer auxiliary variables. In particular, a string of recent work \cite{Adams:2012,Muldoon:2013,Vielma:2016,Vielma:2009a} has presented ideal formulations for certain highly structured constraints such as SOS2 \cite{Beale:1970} with only $\scrO(\log(d))$ auxiliary binary variables and additional constraints (excluding variable bounds). Moreover, these formulations have proven practically useful, and indeed are the most performant by a significant margin for a large swath of instances of the problem classes to which they have been applied \cite{Vielma:2010}. However, building these formulations is complex and ad-hoc, hindering the construction, analysis, and implementation of such formulations for new constraints.

One relatively generic and versatile approach to construct small ideal formulations is the \emph{independent branching} (IB) scheme framework introduced by Vielma and Nemhauser \citep{Vielma:2009a}. The approach is to find some (particularly structured) polyhedra $Q^{1,j}$ and $Q^{2,j}$ such that \eqref{disconst} can be rewritten as
\begin{equation}\label{eqn:disjunction-via-IB}
    \bigcup\nolimits_{i=1}^d P^i = \bigcap\nolimits_{j=1}^t \left(Q^{1,j} \cup Q^{2,j}\right).
\end{equation}
This represents the disjunctive constraint in term of a series of simple choices between two alternatives. Given such a representation, it is often  straightforward to construct a simple, small, and ideal formulation for \eqref{disconst} by formulating each of the $t$ alternatives separately, and then combining them. Furthermore, when the polyhedra $P^i$ are $\scrV$-polyhedra, the construction of the independent branching scheme-based formulation is purely combinatorial, based on the extreme points that are shared between the different polyhedra $P^i$. As we will see, we can therefore approach formulating \eqref{disconst} combinatorially, by studying the shared structure amongst the extreme points.

In this work we generalize and provide a systematic study of the applicability and limitations of the independent branching approach. The contributions of this work can be categorized in the following way.
\begin{enumerate}
    \item We generalize the notion of independent branching schemes to allow for multiple alternatives, and provide an exact characterization of when there exists \emph{any} independent branching representation for \eqref{disconst}, in terms of the graphical representation of the shared extreme points amongst the polyhedra $P^i$. In particular:
    \begin{enumerate}
        \item We demonstrate that the widely-used cardinality constraints cannot be expressed by any independent branching scheme with few alternatives. We argue that this negative result provides theoretical justification for the practical observation that both MIP formulations and simple constraint branching schemes struggle with modeling cardinality constraints effectively.
        \item We show that arbitrary piecewise linear functions in the plane can be modeled with at most three alternatives, and provide a polynomial-time verifiable condition for representability with two alternatives.
        \item We argue that nonconvex polygonal set avoidance constraints are always representable with two alternatives.
    \end{enumerate}
    \item We provide an exact characterization for when there exists a two-alternative independent branching representation for \eqref{disconst} of size $t$, in terms of the classical biclique covering problem. This relation allows the \emph{algorithmic} construction of small independent branching formulations for~\eqref{disconst}. In particular, we study and apply simple properties of biclique covers and their composition to systematically construct explicit descriptions of small independent branching formulations for special structures.

    \item We apply our framework to a variety of constraints of the form \eqref{disconst} to give an indication of the expressive power of the IB scheme approach and the advantage of using biclique cover techniques to construct formulations. In particular:
    \begin{enumerate}
        \item We review and develop several generic properties of biclique covers that lead to systematic construction techniques. Using these techniques we construct explicit, small (i.e. logarithmic in $d$), ideal formulations for generalizations of the special ordered sets of Beale and Tomlin \cite{Beale:1970}, piecewise linear functions over arbitrary 2-dimensional grid triangulations, and outer-approximating discretizations of multilinear terms. This last formulation generalizes the popular logarithmically-sized formulation of Misener et al. \cite{Misener:2012,Misener:2011} for nonconvex quadratic optimization, which we show is not ideal in general.
        \item We provide matching lower bounds for these constructions, showing that they are asymptotically optimal with respect to the size of \emph{any} possible MIP formulation.
    \end{enumerate}
\end{enumerate}

\section{Preliminaries: Definitions, notation, and nomenclature}\label{prelimsec}
A (bounded) $\scrV$-polyhedra (or polyhedra in $\scrV$-form) is a set $P \subset \bbR^n$ that can be expressed as
\[
    P = \Conv(V) \defeq \left\{\sum_{v \in V} \lambda_v v : \lambda \in \Delta^V \right\}
\]
for some finite set of vectors $V \subset \bbR^n$, where $\Delta^V \defeq \{\lambda \in \bbR_{+}^V : \sum_{v \in V} \lambda_v = 1\}$ is the standard simplex. According to the celebrated Minkowski-Weyl Theorem (e.g. \cite[Corollary 3.14]{Conforti:2014}), any polyhedral disjunctive constraint \eqref{disconst} can be expressed as the union of $\scrV$-polyhedra in terms of their extreme points $\ext(P^i)$\footnote{For the moment we are assuming that the $P^i$ are bounded; the unbounded case is more delicate, as we will discuss shortly.}.

When the disjunctive constraint \eqref{disconst} is a union of $\scrV$-polyhedra, it suffices to consider only the combinatorial structure of the extreme points of the polyhedra $P^i$. To see why, consider $J = \bigcup_{i=1}^d \ext(P^i)$ as the \emph{ground set} and take $\scrS = \{\ext(P^i)\}_{i=1}^d \subseteq 2^J$ as the collection of extreme points for each of the polyhedra. We can then define a corresponding disjunctive constraint that is purely combinatorial on the sets $\scrS$.

\begin{definition}\label{CDCDef}
    A \emph{combinatorial disjunctive constraint} (CDC) induced by the sets $\scrS$ is
    \[
        \lambda \in \CDC(\scrS) \defeq \bigcup_{S \in \scrS} Q(S),
    \]
    where $Q(S) \defeq \{\lambda \in \Delta^J : \lambda_{J \backslash S} \leq 0\}$ is the face that $S \subseteq J$ induces on the standard simplex.
\end{definition}

Combinatorial disjunctive constraints may also appear as natural primitive constraints that do not explicitly arise from  unions of $\scrV$-polyhedra, as we will see in Sections~\ref{soskintrosec} and \ref{cardintrosec}. However, when they do arise from unions of $\scrV$-polyhedra, it is straightforward to construct a corresponding formulation for \eqref{disconst} as
\begin{equation} \label{eqn:disaggregated-CDC-to-DC}
    \left\{\sum_{v \in J} \lambda_v v \: : \: \lambda \in \CDC(\scrS) \right\}.
\end{equation}
One advantage of this approach is that formulation \eqref{eqn:disaggregated-CDC-to-DC} allows us to divorce the problem-specific data (i.e. the values $v \in J$) from the underlying combinatorial structure encapsulated in $\CDC(\scrS)$. As such, we can construct a single, strong formulation for a given structure $\CDC(\scrS)$ and this formulation will remain valid for transformations of the data, so long as this transformation sufficiently preserves the  combinatorial structure of $\CDC(\scrS)$. For instance, if $\set{P^i}_{i=1}^d$ are the polyhedra for the original constraint represented by $\CDC(\scrS)$, and $\{\hat{P}^i\}_{i=1}^d$ are those associated with the new data, then a sufficient condition for the formulation of $\CDC(\scrS)$ yielding a valid formulation for $\bigcup_{i=1}^d \hat{P}^i$ is the existence of a bijection $\pi:J \to \hat{J}$ (with $J = \bigcup\nolimits_{i=1}^d \ext(P^i)$ and $\hat{J} = \bigcup\nolimits_{i=1}^d \ext(\hat{P}^i)$) such that
\begin{equation}\label{biyectioncond}
    v \in \ext(P^i) \Longleftrightarrow \pi(v) \in \ext(\hat{P}^i) \quad \forall i \in \llbracket d \rrbracket, \: v \in J,
\end{equation}
where $\llbracket d \rrbracket\defeq \{1,\ldots,d\}$. In this way, we can construct a single small, strong formulation for $\CDC(\scrS)$, and use it repeatedly for many different ``combinatorially equivalent'' instances of the same constraint.

We note that one subtle disadvantage of this data-agnostic approach is that, even if condition \eqref{biyectioncond} is satisfied, the resulting formulation for $\bigcup_{i=1}^d \hat{P}^i$ may be larger than necessary. An extreme manifestation of this would be when the new polyhedra $\{\hat{P}^i\}_{i=1}^d$ are such that $\hat{P}_i\subseteq \hat{P}_1$ for all $i\in \llbracket d \rrbracket$. In this case, $\bigcup\nolimits_{i=1}^d \hat{P}^i = \hat{P}^1$, and so formulating this does not require a MIP formulation at all. Less pathological cases could occur where some subset of the disjunctive sets become redundant after changing the problem data. However, we note that in many of the applications considered in this work, the combinatorial representation leads to redundancy of this form only in rare pathological cases (e.g. Sections~\ref{sos2introsec} and \ref{gridintrosec}). In the remaining cases we will take care to consider, for example, the geometric structure of the data before constructing the disjunctive constraint (e.g. Section~\ref{ss:partitions-of-the-plane}).

Finally, we also note that if we wish to model the case where the polyhedra $P^i$ are unbounded, a result of Jeroslow and Lowe \cite{Jeroslow:1984} \cite[Proposition 11.2]{Vielma:2015} tells us that we may only construct a (binary) MIP formulation for \eqref{disconst} if the recession cones coincide for each $P^i$. In the case this condition is met, we may formulate \eqref{disconst} with
\begin{equation} \label{eqn:disaggregated-CDC-to-DC-with-rays}
    \left\{\sum_{v \in J} \lambda_v v + \sum_{r \in R} \mu_r r \: : \: \lambda \in \CDC(\scrS), \: \mu \in \bbR^R_{+} \right\},
\end{equation}
where $R$ is the shared set of extreme rays for each of the $P^i$. Therefore, we will restrict our attention to the case where each of the $P^i$ are bounded, as formulating the unbounded case is a straightforward extension.

In the remainder of the paper, we will make the following assumptions on $\scrS$ that are without loss of generality.

\begin{assumption}\label{basicassumption}
    We assume the following about $\scrS$.
    \begin{itemize}
        \item $\scrS$ is irredundant: there do not exist distinct $S,T \in \scrS$ such that $S \subseteq T$.
        \item $\scrS$ covers the ground set: $\bigcup_{S \in \scrS} S = J$.
    \end{itemize}
\end{assumption}

We will say that a set $S \subseteq J$ is a \emph{feasible set} with respect to $\CDC(\scrS)$ if $Q(S) \subseteq \CDC(\scrS)$ (equivalently, if $S \subseteq T$ for some $T \in \scrS$) and that it is an \emph{infeasible set} otherwise.

\section{Motivating examples} \label{sec:motivating-examples}
As mentioned, we can represent any polyhedral disjunctive constraint \eqref{disconst} as the union of $\scrV$-polyhedra. However, there are many disjunctive constraints for which the $\scrV$-form of \eqref{disconst} is especially natural. We now present some as running examples that we will return to throughout.

One motif to appear repeatedly will be the graph of a continuous piecewise linear function. That is, given some bounded domain $\Omega \subset \bbR^n$ and some polyhedral partition $\bigcup_{i=1}^d P^i = \Omega$ (i.e. the relative interiors do not overlap\footnote{Formally, $\relint(P^i) \cap \relint(P^j) = \emptyset$ for each $i \neq j$.}), we are interested in modeling a continuous function $f : \Omega \rightarrow \bbR$ such that $x \in P^i \Longrightarrow f(x) = a^i \cdot x + b_i$ for some appropriate $a^i \in \bbR^n$ and $b \in \bbR$. In order to model the graph $\graph(f;\Omega) \defeq \{(x,f(x)) : x \in \Omega\}$, we can construct a formulation for $\CDC(\scrS)$ (with $\scrS = \{\ext(P^i)\}_{i=1}^d$) and express
\begin{equation} \label{eqn:pwl-graph}
    \graph(f) = \left\{ \sum_{v \in J} \lambda_v (v,f(v)) : \lambda \in \CDC(\scrS) \right\},
\end{equation}
where we will use the notation $\graph(f) \equiv \graph(f;\Omega)$ when $\Omega$ is clear from context.\footnote{We note that the results to follow can potentially be extended to certain discontinuous piecewise linear functions by working instead with the epigraph of $f$; we point the interested reader to \cite{Vielma:2010,Vielma:2008a} for further discussion.}

\subsection{Univariate piecewise linear functions and the SOS2 constraint}\label{sos2introsec}
Consider a univariate (nonconvex) piecewise linear function $f$ characterized by $N$ breakpoints $x^1 < x^2 < \ldots < x^N$. We may model the graph of this function via \eqref{eqn:pwl-graph}, where $\scrS = \{\{x^j,x^{j+1}\} : j \in \llbracket N-1 \rrbracket\}$.

As long as the ordering of the breakpoints is preserved, it is easy to see that condition \eqref{biyectioncond} will be satisfied for any transformations of the problem data. Furthermore, the only case in which knowledge of the specific data $\{(x^j,f(x^j))\}_{j=1}^n$ allows the simplification of the original disjunctive representation of $\graph(f)$ is when $f$ is affine in on two adjacent intervals, e.g. affine over $[x^j,x^{j+2}]$ for some $j \in \llbracket N-2 \rrbracket$. Therefore, the potential disadvantage of disregarding the specific data when formulating the constraint occurs only in rare pathological cases which are easy to detect.

For this reason, we will strip out the problem data and instead express the constraint with respect to the indices $j$ of the vertices $\{x^j\}_{j=1}^N$. That is, we take $J = \llbracket N \rrbracket $ and write $\scrS = \{\{\tau,\tau+1\} : \tau \in \llbracket N-1 \rrbracket\}$ to emphasize the data independence of the constraint and highlight the combinatorial structure. In this form, we can recognize  the special ordered set of type 2 (SOS2$(N)$) constraint of Beale and Tomlin \cite{Beale:1970}, which requires that at most two components of $\lambda$ may be nonzero, and  that these nonzero components must be consecutive in the ordering on $J$.

\subsection{SOS$k$}\label{soskintrosec}
A generalization of the special ordered sets considers the case where at most $k$ consecutive components of $\lambda$ may be nonzero at once. In particular, if $J = \llbracket N \rrbracket$, we have $\scrS = \{\{\tau,\tau+1,\ldots,\tau+k-1\} : \tau \in \llbracket N - k + 1 \rrbracket\}$. This constraint may arise, for example, in chemical process scheduling problems, where an activated machine may only be on for $k$ consecutive time units and must produce a fixed quantity during that period \cite{Floudas:2005,Kondili:1993}.

\subsection{Cardinality constraints}\label{cardintrosec}
An extremely common constraint in optimization is the cardinality constraint of degree $\ell$, where at most $\ell$ components of $\lambda$ may be nonzero. This corresponds to $\scrS = \{I \subseteq J : |I| = \ell\}$. A particularly compelling application of the cardinality constraint is in portfolio optimization \cite{Bertsimas:2007a,Bienstock:1996,Chang:2000,Vielma:2008}, where it is often advantageous to limit the number of investments to some fixed number $\ell$ to minimize transaction costs, or to allow differentiation from the performance of the market as a whole.

\subsection{Discretizations of multilinear terms}\label{multilinearsec}
Consider a multilinear function $f(x_1,\ldots,x_\eta) = \prod_{i=1}^\eta x_i$ defined over some box domain $\Omega \defeq [l,u] \subset \bbR^\eta$. This function appears often in optimization models \cite{Foulds:1992}, but is nonconvex, and often leads to problems which are difficult to solve to global optimality in practice \cite{Androulakis:1995,Quesada:1995,Wicaksono:2008}. As a result, computational techniques will often ``relax'' the graph of the function $\graph(f)$ with a convex outer approximation, which is easier to optimize over \cite{Sahinidis:1996}.\footnote{Note that these relaxations are useful in the context of global optimization, coupled with algorithmic techniques such as spatial branch-and-bound.}

For the bilinear case ($\eta=2$), the well-known McCormick envelope \cite{McCormick:1976} describes the convex hull of $\graph(f)$. Although traditionally stated in an inequality description, we may equivalently describe the convex hull via its four extreme points, which are readily available in closed form. For higher-dimensional multilinear terms, the convex hull has $2^\eta$ extreme points, and can be constructed in a similar manner (e.g. see equation (3) in~\cite{luedtke2012some} and the associated references).

Misener et al. \cite{Misener:2012,Misener:2011} propose a computational technique for optimizing problems with bilinear terms where, instead of modeling the graph over a single region $\Omega = [l,u] \subset \bbR^2$, they discretize the region in a regular fashion and apply the McCormick envelope to each subregion. They model this constraint as a union of polyhedra, where each subregion enjoys a tighter relaxation of the bilinear term. Additionally, they propose a logarithmically-sized formulation for the union. However, it is not ideal (see Appendix~\ref{app:misener}), it only applies for bilinear terms ($\eta=2$), and it is specialized for a particular type of discretization (namely, only discretizing along one component $x_1$, and with constant discretization widths $\{h^1_{j+1}-h^1_{j}\}_{j=1}^{d_1-1}$).

For a more general setting, we have that the extreme points of the convex hull of the graph $\Conv(\graph(f))$ are given by $\{(x,f(x)) : x \in \ext(\Omega)\}$ \cite[equation (3)]{luedtke2012some}, where it is easy to see that $\ext(\Omega) = \prod_{i=1}^\eta \{l_i,u_i\}$. Consider a grid imposed on $[l,u] \subset \bbR^\eta$; that is, along each component $i \in \llbracket \eta \rrbracket$, we partition $[l_i,u_i]$ along the points $l_i \equiv  h_1^i < h_2^i < \cdots < h_{d_i-1}^i < h_{d_i}^i \equiv u_i$. This yields $\prod_{i=1}^\eta \bra{d_i-1}$ subregions; denote them by $\scrR \defeq \left\{\prod_{i=1}^\eta [h^i_{k_i},h^i_{k_i+1}] : k \in \prod_{i=1}^\eta \llbracket d_i -1 \rrbracket\right\}$.

We can then take the polyhedral partition of $\Omega$ given by $P^R \defeq \Conv(\graph(f;R))$ for each subregion $R$, the sets as $\scrS = \left\{\ext(P^R)\right\}_{R \in \scrR}$, and the ground set as $J = \bigcup\{S \in \scrS\}$. In particular, we have that $J = \prod_{i=1}^\eta \{h^i_1,\ldots,h^i_{d_i}\}$. Analogously to the notational simplification we took with the SOS2 constraint, for the remainder we will take $J=\prod_{i=1}^\eta \llbracket d_i \rrbracket$ and  $\scrS=\set{\prod_{i=1}^\eta \set{k_i,k_i+1} : k \in \prod_{i=1}^\eta \llbracket d_i -1 \rrbracket}$. We also note that condition \eqref{biyectioncond} is satisfied as long as the ordering of the discretization is respected along each dimension.

\subsection{Piecewise linear functions in the plane and grid triangulations}\label{gridintrosec}

Consider a (potentially nonconvex) region $\Omega \subset \bbR^2$. We would like to model a (also potentially nonconvex) piecewise linear function $f$ with domain over $\Omega$. Take $\{P^i\}_{i=1}^d$ as the set of pieces of the domain, and the corresponding ground set as $J = \bigcup_{i=1}^d \ext(P^i)$ and sets as $\scrS = \{\ext(P^i)\}_{i=1}^d$. We may then model the piecewise linear function via the graph representation \eqref{eqn:pwl-graph}.

\begin{figure}[htpb]
    \centering
    \includegraphics[width=.25\linewidth]{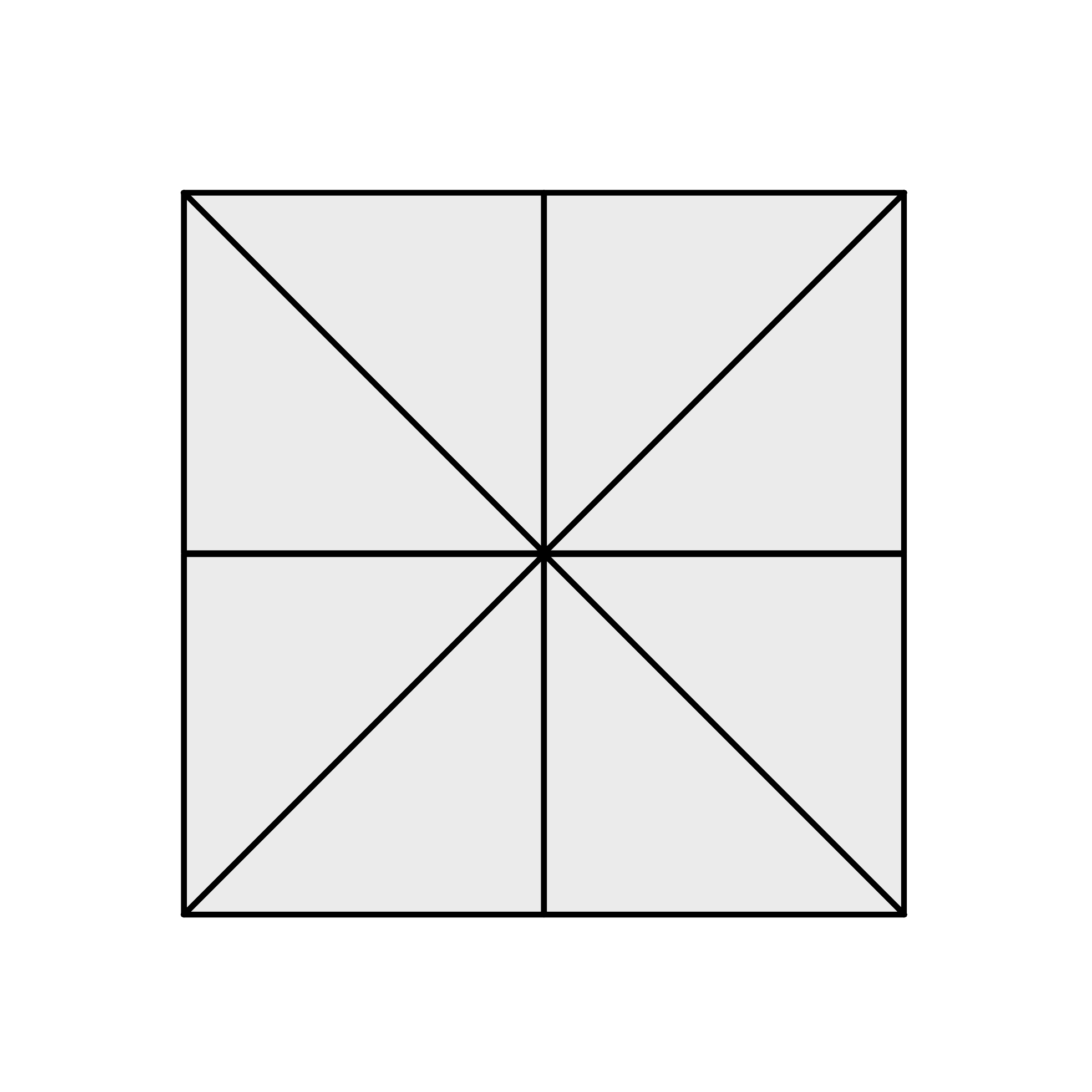}
    \includegraphics[width=.25\linewidth]{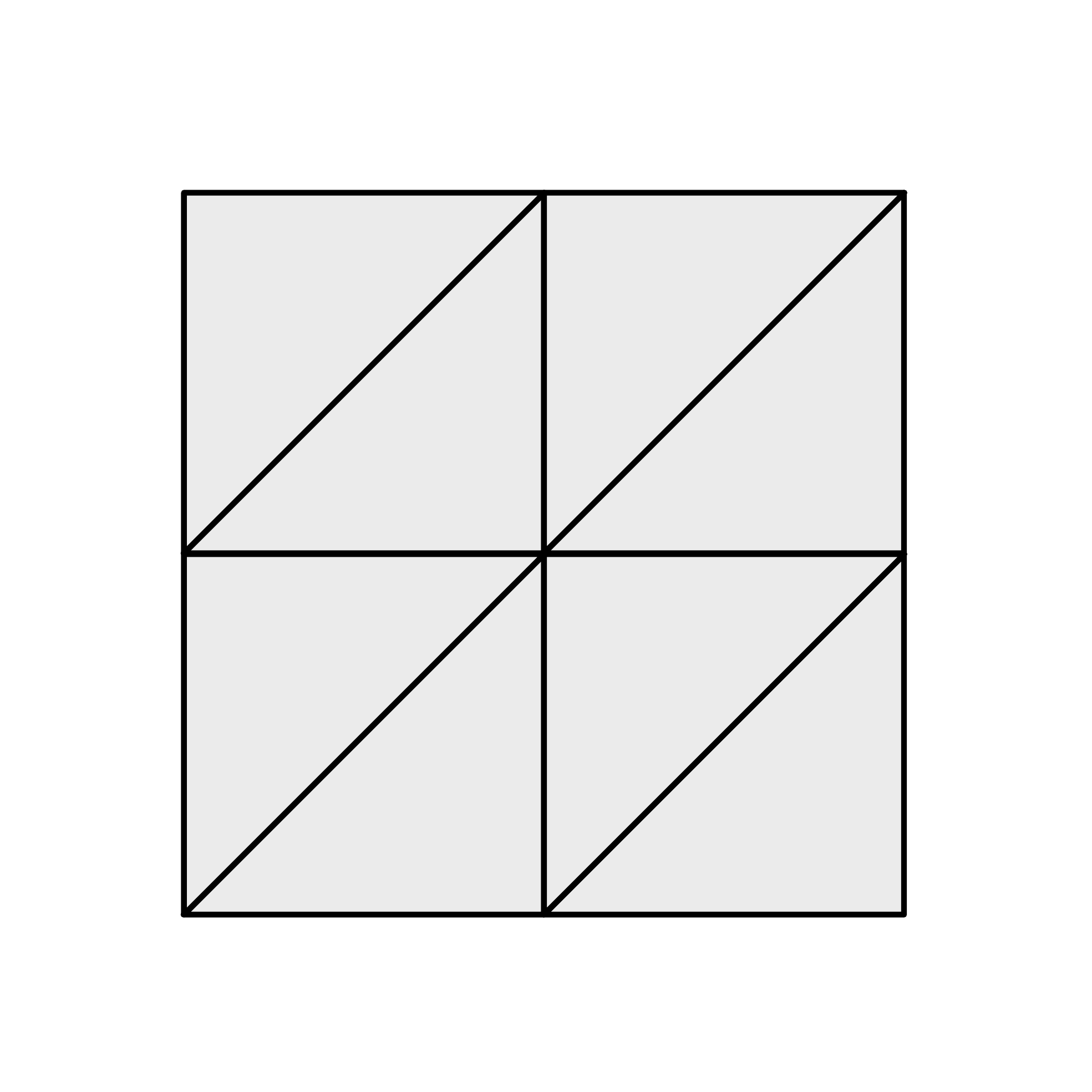}
    \includegraphics[width=.25\linewidth]{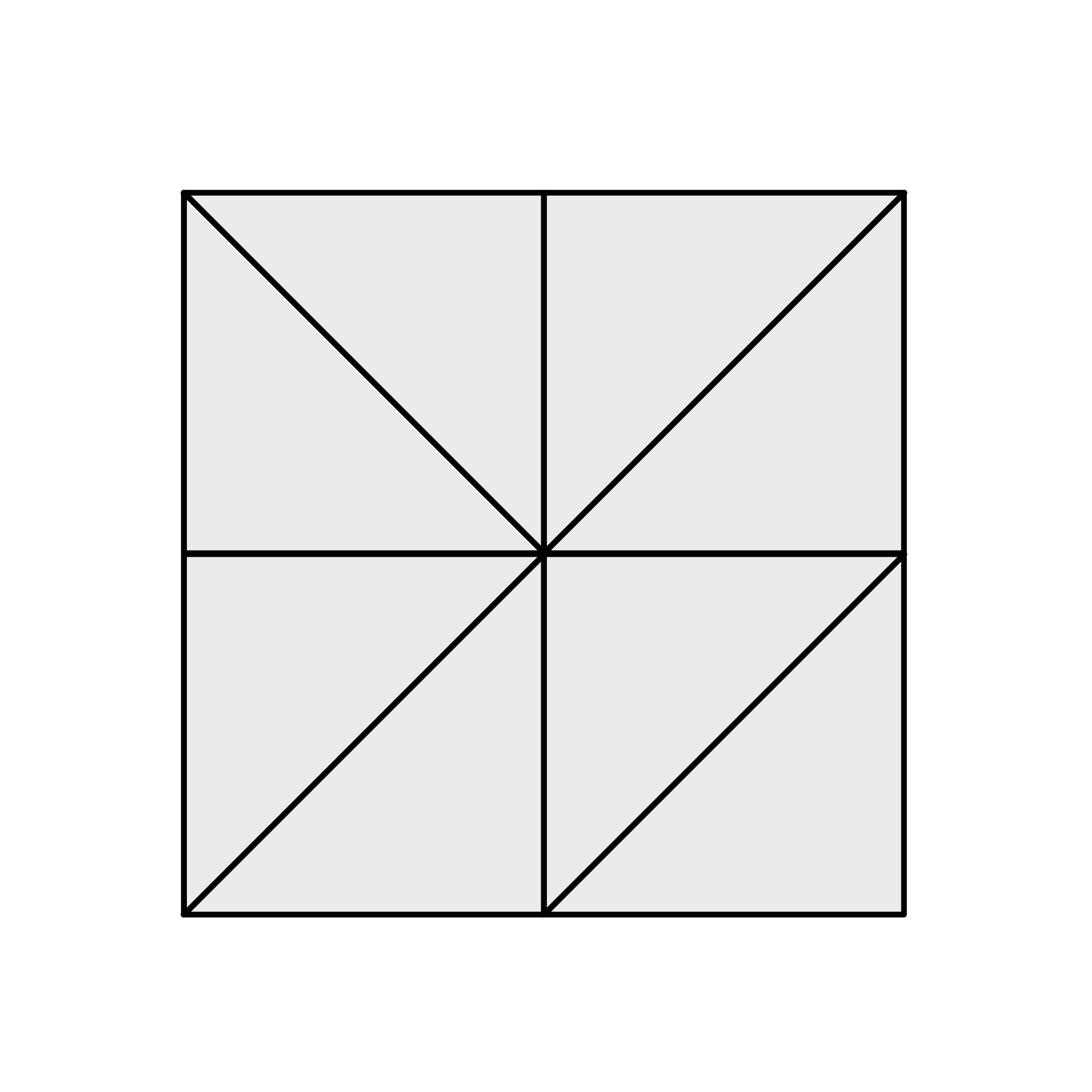}
    \caption{Three grid triangulations of $\Omega = [1,3]\times[1,3]$: the Union Jack (J1) \cite{Todd:1977} (Left), the K1 \cite{Kuhn:1960} (Center), and a more idiosyncratic construction (Right).}
    \label{fig:triangulations}
\end{figure}

An important special case occurs when the function $f$ is affine over a triangulation of a grid similar to the one used for multilinear terms. A description of the associated CDC that emphasizes its combinatorial structure can be obtained through the same simplification from Section~\ref{multilinearsec} as follows. Consider a rectangular region in the plane $\Omega = [1,M] \times [1,N]$, and the regular grid points $J = \{1,\ldots,M\} \times \{1,\ldots,N\}$. A \emph{grid triangulation} $\scrS$ of $\Omega$ is then a set $\scrS$ where:
\begin{itemize}
    \item Each $S \in \scrS$ is a triangle: $|S|=3$.
    \item $\scrS$ partitions $\Omega$: $\bigcup_{S \in \scrS} \Conv(S) = \Omega$ and $\relint(\Conv(S)) \cap \relint(\Conv(T)) = \emptyset$ for each distinct $S, T \in \scrS$.
    \item $\scrS$ is on a regular grid: $S \subset J$ for each $S \in \scrS$, and $||v-w||_\infty \leq 1$ for each $v,w \in S$.
\end{itemize}
As concrete examples, consider Figure~\ref{fig:triangulations}, where we depict three different triangulations with $M=N=3$.
Grid triangulations are often used to model bivariate, non-separable, piecewise linear functions using \eqref{eqn:pwl-graph}~\cite{Vielma:2010,Vielma:2009a}. Furthermore, any formulation constructed for a given grid triangulation can be readily applied to any other grid triangulation obtained by shifting the grid points in the plane, so long as the resulting triangulation is strongly isomorphic to, or compatible with, the original triangulation \cite{aichholzer2003towards}.

\subsection{Obstacle avoidance}\label{obstacleavoidancesec}
Consider an unmanned aerial vehicle (UAV) which you would like to navigate through an area with fixed obstacles. At any given time, you wish to impose the constraint that the location of the vehicle $x \in \bbR^2$ must lie in some (nonconvex) region $\Omega \subset \bbR^2$, which is the plane, less any obstacles in the area. MIP formulations of this constraint has received interest as a useful primitive for path planning \cite{Bellingham:2002,Deits:2015,Prodan:2016,Mellinger:2012}.

We may model $x \in \Omega$ by partitioning the region $\Omega$ with polyhedra such that $\Omega = \bigcup_{i=1}^d P^i$. Traditional approaches to modeling constraints \eqref{disconst} of this form use a linear inequality description for each of the polyhedra $P^i$ and construct a corresponding big-$M$ formulation \cite{Prodan:2012,Prodan:2016}, which will not be ideal in general. In the $\scrV$-polyhedra framework, we will instead take the $P^i$ as $\scrV$-polyhedra that partition $\Omega$, and be able to construct small, ideal formulations.

\section{MIP formulations for combinatorial disjunctive constraints}
Using standard MIP formulation techniques, we now present formulations for $\CDC(\scrS)$ for comparison with our approach. In particular, we will argue that the framework we will present later can lead to formulations that are smaller (in terms of the number of auxiliary variables), and enjoy other favorable properties we enumerate in Section \ref{ss:IB-schemes}.

\subsection{MIP formulations: Definitions, size, and strength}\label{ss:formulation-definitions}
Formally, we say that a (binary) MIP formulation $F$ for a constraint $x \in Q \subseteq \bbR^{n_1}$ is the composition of linear inequalities (the \emph{linear programming (LP) relaxation}, or just \emph{relaxation} in the context of this work)
\[
    R \defeq \left\{ (x,y,z) \in [l^x,u^x] \times [l^y,u^y] \times [0,1]^{n_3} : Ax + By + Cz \leq d \right\}
\]
with (binary) integrality conditions
\begin{equation}\label{eqn:generic-MIP-formulation}
    F \defeq R \cap \left(\bbR^{n_1} \times \bbR^{n_2} \times \{0,1\}^{n_3}\right)
\end{equation}
such that $\Proj_x(F) = Q$.\footnote{To handle the unbounded case, we allow the variable bounds $l^x \leq x \leq u^x$ and $l^y \leq y \leq u^y$ to take infinite values.} Finally, we assume $R$ is line-free (i.e. has at least one extreme point), which is satisfied by essentially all practical formulations.

Throughout, we will be interested in ways of understanding both the strength of a given formulation, as well as in quantifying the size or complexity of a formulation. We say that a formulation is \emph{ideal} if each extreme point of the relaxation naturally satisfies the integrality conditions, i.e. $\Proj_z(\ext(R)) \subseteq \{0,1\}^{n_3}$. The choice of name is apt, as this is the strongest possible MIP formulation we can expect.\footnote{An ideal formulation is also \emph{sharp}, i.e. its relaxation projects down to the convex hull of the set we are formulating ($\Proj_x(R) = \Conv(Q)$).}

As a measure of the complexity of the formulation, we count the number of auxiliary continuous variables $y$ and continuous binary variables $z$ used by the formulation, as well as the number of inequalities in our description of $R$. We ignore the size of $x$, since this is intrinsic to the constraint we wish to model. We will say that a formulation is \emph{extended} if there are auxiliary continuous variables $y$ in the representation $F$ (that is, $n_2 > 0$) and \emph{non-extended} otherwise ($n_2 = 0$). Furthermore, as suggested by the definition of $R$, we distinguish between variable bounds (e.g. $l^y \leq y \leq u^y$) and general inequalities ($Ax + By + Cz \leq d$), as modern MIP solvers are able to incorporate variable bounds with minimal extra computational cost. Finally, we note that from now on (with a single exception), any statements regarding formulation properties or size will be with respect to the combinatorial disjunctive constraint $\CDC(\scrS)$, rather than the $\scrV$-polyhedral disjunctive constraints itself. In other words, all statements are with respect a fixed combinatorial structure, with no possible simplifications from the specific data as discussed in Section~\ref{prelimsec}. The sole exception will be in Section~\ref{ss:partitions-of-the-plane}, where we will make use of some geometric properties inherent in a given realization of the data.

\subsection{Existing formulations for combinatorial disjunctive constraints}

A standard formulation for $\CDC(\scrS)$ adapted from Jeroslow and Lowe \cite{Jeroslow:1984} is
\begin{subequations} \label{eqn:jeroslow-formulation}
\begin{align}
    \lambda_v = \sum_{S \in \scrS : v \in S} \gamma^S_v &\quad\quad \forall v \in J \\
    z_S = \sum_{v \in S} \gamma^S_v &\quad\quad \forall S \in \scrS \\
    \sum_{S \in \scrS} z_S = 1 &\\
    \gamma^S \in \Delta^S &\quad\quad \forall S \in \scrS \\
    (\lambda,z) \in \Delta^J \times \{0,1\}^\scrS.
\end{align}
\end{subequations}
This formulation has $d=|\scrS|$ auxiliary binary variables, $\sum_{S \in \scrS} |S|$ auxiliary continuous variables, and no general inequalities. Additionally, it is ideal.

Using Proposition 9.3 from \cite{Vielma:2015}, we can construct an ideal MIP formulation with fewer auxiliary binary variables:
\begin{subequations} \label{eqn:prop-9.3}
    \begin{align}
        \lambda_v = \sum_{S \in \scrS : v \in S} \gamma^S_v &\quad\quad \forall v \in J \\
        \sum_{S \in \scrS} \sum_{v \in S} \gamma^S_v = 1 & \\
        \sum_{S \in \scrS} \sum_{v \in S} h^S \gamma^S_v = z & \\
        \gamma^S \geq 0 &\quad\quad \forall S \in \scrS \\
        z \in \{0,1\}^r,&
    \end{align}
\end{subequations}
where $\{h^S\}_{S \in \scrS} \subseteq \{0,1\}^r$ is some set of distinct binary vectors. This formulation is actually a generalization of \eqref{eqn:jeroslow-formulation}, which we recover if we take $h^S = \textbf{e}^S \in \bbR^\scrS$ as the canonical unit vectors. If instead we take $r$ to be as small as possible (while ensuring that the vectors $\{h^S\}_{S \in \scrS}$ are distinct), we recover $r = \lceil \log_2(d) \rceil$. Therefore, formulation \eqref{eqn:prop-9.3} yields an ideal extended formulation for \eqref{disconst} with $\lceil \log_2(d) \rceil$ auxiliary binary variables, $\sum_{S \in \scrS} |S|$ auxiliary continuous variables, and no general inequalities. The following corollary shows that this is the smallest number of auxiliary binary variables we may hope for.

\begin{proposition} \label{prop:log-bound}
    If the sets $\scrS$ are irredundant, then any binary MIP formulation for $\CDC(\scrS)$ must have at least $\lceil\log_2(d)\rceil$ auxiliary binary variables.
\end{proposition}
\proof{\textbf{Proof}}
    See Appendix~\ref{app:log-bound}.
\Halmos\endproof

The formulations thus far have been extended formulations, as they are constructed by formulating each polyhedra $Q\bra{S}$ separately and then aggregating them, rather than working with the combinatorial structure underlying the shared extreme points. Therefore, each of these formulations requires a copy of the multiplier $\gamma^S_v$ for each set $S \in \scrS$ for which $v \in S$, and so $\sum_{i=1}^d |\scrS|$ auxiliary continuous variables total.

In contrast, we can construct non-extended formulations for $\CDC(\scrS)$ that work directly on the $\lambda$ variables and the underlying combinatorial structure of $\scrS$. An example of a non-extended formulation for CDC is the widely used ad-hoc formulation (see \cite[Section 6]{Vielma:2015} and the references therein) given by
\begin{subequations} \label{eqn:disaggregated}
\begin{align}
    \lambda_v \leq \sum_{S \in \scrS: v \in S} z_S &\quad\quad \forall v \in J \label{eqn:disaggregated-1} \\
   \sum_{S \in \scrS} z_S = 1 \\
   (\lambda,z) \in \Delta^J \times \{0,1\}^{\scrS}.
\end{align}
\end{subequations}
This formulation  is not necessarily ideal, and it requires no auxiliary continuous variables, $d$ auxiliary binary variables, and $|J|$ general inequalities.

In summary, we have seen an ideal extended formulation \eqref{eqn:prop-9.3} for $\CDC(\scrS)$ with relatively few auxiliary binary variables, but relatively many auxiliary continuous variables. On the other end of the spectrum, we have a non-extended formulation \eqref{eqn:disaggregated} with no auxiliary continuous variables, but which requires relatively many auxiliary binary variables and which may fail to be ideal. However, we know that in special cases we can construct ideal, non-extended formulations with only $\scrO(\log(d))$ auxiliary variables and constraints (e.g. SOS1, SOS2, and particular 2-dimensional grid triangulations \cite{Vielma:2016,Vielma:2009a}). This work provides a framework for constructing such \emph{small, strong}, non-extended MIP formulations for $\CDC(\scrS)$, which are automatically ideal, and in the best case will have $\scrO(\log(d))$ auxiliary binary variables and general inequality constraints.

\section{Independent branching schemes} \label{ss:IB-schemes}
Vielma and Nemhauser \cite{Vielma:2009a} introduced the notion of an \emph{independent branching scheme} as a natural framework for constructing formulations for combinatorial disjunctive constraints. The independent branching scheme is a logically equivalent way of expressing a CDC in terms of a conjunction of dichotomies: that is, as a series of choices between two (simple) options. This approach is parsimonious: if you are given an independent branching scheme for a particular CDC, it is straightforward to construct an ideal formulation whose size is on the order of the number of dichotomies. For our purposes, we present a generalized notion, where we allow potentially more than two alternatives.

\begin{definition}
    A \emph{$k$-way independent branching scheme} for $\CDC(\scrS)$ is given by a family of sets $(L^j_1,\ldots,L^j_k)$ (where each $L^j_i \subseteq J$) for $j \in \llbracket t \rrbracket$, where
    \begin{equation} \label{eqn:multi-way-CDC}
        \CDC(\scrS) = \bigcap_{j=1}^t \left( \bigcup_{i=1}^k Q(L^j_i) \right).
    \end{equation}
\end{definition}
We say that such an IB scheme has \emph{depth $t$}, and that each $j \in \llbracket t \rrbracket$ yields a corresponding \emph{level} of the IB scheme $\bigcup_{i=1}^k Q(L^j_i)$, given by the $k$ \emph{alternatives} $Q(L^j_i)$.

An equivalent way of understanding these representations, which we will be using for the remainder of this work, is by eschewing the polyhedra $Q(L^j_i)$ and working directly on the underlying set $L^j_i$. That is, a valid $k$-way IB scheme satisfies the condition that
\[
    T \subseteq J \text{ is a feasible set} \Longleftrightarrow \forall j \in \llbracket t \rrbracket, \: \exists i \in \llbracket k \rrbracket \text{ s.t. } T \subseteq L^j_i.
\]
First, we observe that, due to our assumption that $\scrS$ covers the ground set, we have that for each element $v \in J$ and level $j \in \llbracket t \rrbracket$, there will be at least one alternative $i \in \llbracket k \rrbracket$ such that $v \in L^j_i$. We will use this extensively in the analysis to come, as it simplifies some otherwise tedious case analyses. Second, we see that this definition can capture potential schemes with a variable number of alternatives in each level by adding empty alternatives $L^j_i=\emptyset$, provided we take $k$ as the maximum number of alternatives for all levels. For notational simplicity, we say that a 2-way IB scheme is a \emph{pairwise IB scheme}, and in this case we write the sets as $\{(L^{j},R^{j})\}_{j=1}^t$ as in \cite{Vielma:2009a}. In contrast, we will call the case with $k > 2$ a \emph{multi-way IB scheme}.

In this form, we have replaced the monolithic constraint $\CDC(\scrS)$ by $t$ constraints, each of which require the selection between $k$ alternatives. We may then use standard techniques to construct a corresponding mixed-integer formulation.

\begin{proposition}
    Given an independent branching scheme $\{(L_1^j,\ldots,L_k^j)\}_{j=1}^t$ for $\CDC(\scrS)$, the following is a valid formulation for $\CDC(\scrS)$:
    \begin{subequations} \label{eqn:multiway-formulation}
    \begin{align}
        \sum_{v \not\in L^j_i} \lambda_v \leq 1-z^j_i &\quad \forall j \in \llbracket t \rrbracket, \forall i \in \llbracket k \rrbracket \label{eqn:multiway-formulation-1} \\
        \sum_{i=1}^k z^j_i = 1 &\quad \forall j \in \llbracket t \rrbracket \label{eqn:multiway-formulation-2} \\
        \lambda \in \Delta^J &\\
        z^j \in \{0,1\}^k &\quad \forall j \in \llbracket t \rrbracket.
    \end{align}
    \end{subequations}
\end{proposition}
The formulation is known to be ideal for $k=2$ \cite{Vielma:2010,Vielma:2009a}. It has no auxiliary continuous variables, $k t$ auxiliary binary variables, and $k t$ general inequalities. 

\subsection{Constraint branching via independent branching-based formulations}\label{constraintbranchingsec}
The canonical algorithmic technique for solving mixed-integer programming problems is some variation of branch-and-bound \cite{Land:1960}, which implicitly enumerates all possible values for the binary variables. In its simplest form, a sequence of problems are solved, starting with the relaxation of the MIP formulation, after which a binary variable $z_i$ is chosen for \emph{branching}. That is, the current problem is branched into two subproblems: one with the additional constraint $z_i \leq 0$, another with $z_i \geq 1$. Repeating this procedure, the subproblems form a (binary) tree whose leaves correspond to all $2^{n_3}$ possible values for the $n_3$ binary variables in formulation~\eqref{eqn:generic-MIP-formulation}. At any given subproblem, the augmented relaxation to be solved is described by the set of binary variables fixed to zero, and the set of those fixed to one.

The spirit of constraint branching is to allow richer branching decisions. For example, a branching decision might be between $k$ alternatives of the form $\{Q^i\}_{i=1}^k$, where each $Q^i$ is formed by adding a general inequality constraint to the existing relaxation at the current node. This concept has significant overlap with the broader field of constraint programming \cite{Apt:2003,Jaffar:1994}, which has been recognized and exploited in the mixed-integer programming literature \cite{Achterberg:2008,Appleget:2000,Hooker:2002,Ostrowski:2009,Ryan:1981}. More complex constraint branching can often lead to a more balanced branch-and-bound tree, which can significantly improve computational performance (see, for example, \cite[Section 8]{Vielma:2015} and \cite{Yildiz:2013} for more discussion). Combinatorial disjunctive constraints are a natural setting to apply constraint branching directly on the continuous $\lambda$ variables \cite{Beale:1970,DeFarias01,keha04,keha06,martin06}. Indeed, the classical examples of the SOS1 and SOS2 constraints~\cite{Beale:1970} show that we do not necessarily require a MIP formulation (or the auxiliary binary variables $z$) for modeling combinatorial disjunctive constraints, as the disjunction can be enforced directly through constraint branching on the $\lambda$ variables. These constraint branching approaches without auxiliary binary variables can be implemented in an ad-hoc branch-and-bound procedure, or through branching callbacks available in some MIP solvers, such as CPLEX. In theory, this approach should outperform a MIP formulation like \eqref{eqn:disaggregated} that introduces additional variables and constraints. However, realizing this performance advantage in practice can require significant effort and technical expertise. For instance, Vielma et al.~\cite{Vielma:2010,Vielma:2008a} observe that the basic formulation \eqref{eqn:disaggregated} clearly outperformed the SOS2 branching implementation in CPLEX v9.1. However, CPLEX v11 implemented an optimized version of SOS2 branching that used the advanced branch selection techniques available for variable branching, reversing this performance gap with respect the MIP formulation approach.

One way to avoid re-implementing the advanced branching selection techniques for a new constraint branching approach is by constructing a MIP model that automatically inherits the advanced constraint branch selection, but using the solvers traditional variable branching~\cite{Appleget:2000,Vielma:2009a}. For simplicity, assume that the constraint branching approach has $t$ branching options, each of which creates $k$ branches, and that each constraint added has support on the $\lambda$ variables, with variable coefficients in $\set{0,1}$ and a zero right-hand-side. That is, branch $i \in \llbracket k \rrbracket$ of branching option $j \in \llbracket t \rrbracket$ adds a constraint of the form $\sum_{v \not\in L^j_i} \lambda_v \leq 0$. This is equivalent to a \emph{multi-variable} branching approach that fixes groups of variables to zero, and it includes as special cases most constraint branching approaches, including SOS1/SOS2 branching. Then \eqref{eqn:multiway-formulation} is a MIP formulation for this multi-variable constraint branching scheme, as variable branching on $\{z^j_i\}_{i=1}^k$ enforces constraint branching option $j$ on the $\lambda$ variables.

This connection highlights the natural theoretical equivalence between a multi-variable branching and an independent branching formulation. A practical difference between the two is that direct multi-variable branching must implement an explicit branch selection and implementation routine, while an independent branching formulation inherits the variable branching selection and implementation routines of the MIP solver. The upshot of this is that the independent branching formulation must provide a complete catalog of all possible branching options up-front (i.e. through the formulation), while direct multi-variable branching can have a large catalog of branching options that are implicitly defined by the branching routines.

We finally note that, as discussed by Vielma and Nemhauser~\cite[Section 3]{Vielma:2009a}, variable branching on a non-IB formulation such as \eqref{eqn:prop-9.3} can fix components of $\lambda$ to zero in a fashion that is \emph{dependent} on previous branching decisions. This is not the case with independent branching formulations, hence the name. We review this independence property in detail in Appendix~\ref{IBappendix}.

\section{Independent branching scheme representability}
To start, we observe that the independent branching approach is not sufficiently general to capture every possible formulation for $\CDC(\scrS)$. In particular, there is the restriction that each alternative $Q(L^j_i)$ restricts the $\lambda$ variables to lie on a \emph{single} face of the standard simplex. A natural first question is then: given a family of sets $\scrS$, do any $k$-way IB schemes exist for $\CDC(\scrS)$? We provide an answer, based on a graphical characterization of the constraint.

\begin{definition}
    Let $H\defeq (J,\mathcal{E})$ be a hypergraph with hyperedge set $\mathcal{E}\subseteq 2^J$.
    \begin{itemize}
        \item The \emph{rank} of $H$ is $r\bra{H} \defeq \max\set{\abs{E}\,:\, E\in \mathcal{E}}$.
        \item A \emph{(weakly) independent set} of $H$ is a set $U\subseteq J$ that does not contain any element of $\mathcal{E}$ as a subset.
        \item The \emph{conflict hypergraph} of $\scrS$ is $H^c_\scrS \defeq (J,\mathcal{E}_\scrS)$, where $\mathcal{E}_\scrS \defeq \{E\subseteq J : \text{ $E$ is a minimal infeasible set}\}$.
    \end{itemize}
\end{definition}

\begin{lemma}
    The maximal independent sets $S$ in $H^c_\scrS$ are exactly the sets $S \in \scrS$.
\end{lemma}
\proof{\textbf{Proof}}

If $S \in \scrS$, it is obviously a feasible set, and so we have immediately that $S$ is an independent set in $H^c_\scrS$. If it is not maximal, then we could add some $v \in J \backslash S$ and maintain feasibility, which would violate our irredundancy assumption (i.e. $S \cup \{v\} \subseteq S' \in \scrS$, $S\in \scrS$, and  $ S \subsetneq S \cup \{v\} \subseteq S'$).

If $S$ is a maximal independent set in $H^c_\scrS$, then it must be a feasible set with respect to $\scrS$ as well. As it is maximal, there is no set $\hat{S} \in \scrS$ with $S \subsetneq \hat{S}$, and so we must have $S \in \scrS$ as well.
\Halmos\endproof

\begin{theorem}\label{thm:representability}
    A $k$-way IB scheme exists for $\CDC(\scrS)$ if and only if  $r\bra{H^c_\scrS} \leq k$. In particular, if $\mathcal{E}_\scrS=\set{E^j = \{e^j_1,\ldots,e^j_{|E^j|}\}}_{j=1}^t$ is the hyperedge set for the conflict hypergraph $H^c_\scrS$, then an $r\bra{H^c_\scrS}$-way IB scheme for $\CDC(\scrS)$ is given by
    \begin{equation}\label{basicIB}
        L^j_i = \begin{cases}
            J \backslash \{e^j_i\} & i \leq \abs{E^j} \\
            \emptyset & \text{o.w.}
        \end{cases} \quad \forall i \in \llbracket r\bra{H^c_\scrS} \rrbracket, \: j\in \llbracket t \rrbracket.
     \end{equation}
\end{theorem}
\proof{\textbf{Proof}}
To show the ``if'' direction, it suffices to show the validity of \eqref{basicIB}. First note that every minimally infeasible set $E^j \in \mathcal{E}_\scrS$ is rendered infeasible by level $j$, which implies that every infeasible set is rendered infeasible as well. Then note that for any  $S \in \scrS$ and for any $j\in \llbracket t \rrbracket$, we have $E^j \not\subseteq S$, so there exists $i \in \llbracket \abs{E^j} \rrbracket $ such that $e^j_i\in E^j \setminus S$. Hence, $S \in L^j_i$ and $S$ is feasible for level $j$.

To show the ``only if'' direction, assume for a contradiction that there exists a $k$-way IB scheme with $k\leq r\bra{H^c_\scrS}-1$. Take a minimal infeasible set $E =\{e_1,\ldots,e_r\} \in \mathcal{E}_\scrS$, where $r = r(H_\scrS^c)$. Then take $j \in \llbracket t \rrbracket$ as a level of the IB scheme that renders $E$ infeasible. By the minimality of $E$, we have that, for all $\ell \in \llbracket r \rrbracket$, there exists some $i(\ell) \in \llbracket k \rrbracket$ such that $E(\ell) \defeq E \setminus \set{e_\ell} \subseteq L^j_{i(\ell)}$. As $k < r$, we may apply the pigeonhole principle to see that there must exist some distinct $\ell^1,\ell^2 \in \llbracket r \rrbracket$ such that $i(\ell^1)=i(\ell^2)$, and such that $E(\ell^1) \subseteq L^j_{i(\ell^1)}$ and $E(\ell^2) \subseteq L^j_{i(\ell^2)}$. As $E = E(\ell^1) \cup E(\ell^2)$ and $L^j_{i(\ell^1)} = L^j_{i(\ell^2)}$, this implies that $E \subseteq L^j_{i(\ell^1)}$, which contradicts our supposition that level $j$ rendering $E$ infeasible.
\Halmos\endproof

Throughout, we will say that $\CDC(\scrS)$ is $k$-way IB-representable (or pairwise IB-representable for $k=2$) if it admits a $k$-way IB scheme.

\subsection{Cardinality constraints} \label{ss:cardinality}
Our first application of Theorem~\ref{thm:representability} is to derive a strong restriction on the existence of multi-way IB schemes for the cardinality constraint.

\begin{corollary} \label{prop:cardinality}
   A cardinality constraint of degree $\ell$ is $k$-way IB-representable if and only if $k > \ell$.
\end{corollary}
\proof{\textbf{Proof}}
    Direct from Theorem~\ref{thm:representability} by observing that the conflict hypergraph has rank $\ell+1$.
\Halmos\endproof
We observe that the IB scheme \eqref{basicIB}, when applied to the cardinality constraint, is a natural MIP formulation for the ``conjunctive normal form'' \cite{Balas:1985}, and is unlikely to be practical for even moderately large $\ell$. In addition, both specialized constraint branching schemes for cardinality constraints \cite{Ismael-R.-de-Farias:2003} and the binary variable branching induced by standard formulations for cardinality constraints are quite imbalanced. The existence of a pairwise independent branching scheme for cardinality constraints would likely have finally produced the sought-after balanced constraint branching. However, Corollary~\ref{prop:cardinality} implies that such a balanced constraint branching cannot be produced via IB schemes, or equivalently by  constraint branchings that do not use general inequalities (i.e. are only multi-variable branchings).

\subsection{Polygonal partitions of the plane} \label{ss:partitions-of-the-plane}
\begin{figure}[htpb]
    \centering
    \includegraphics[width=.32\linewidth]{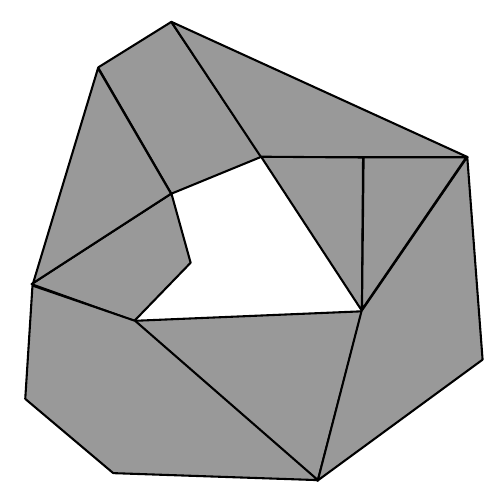}
    \includegraphics[width=.32\linewidth]{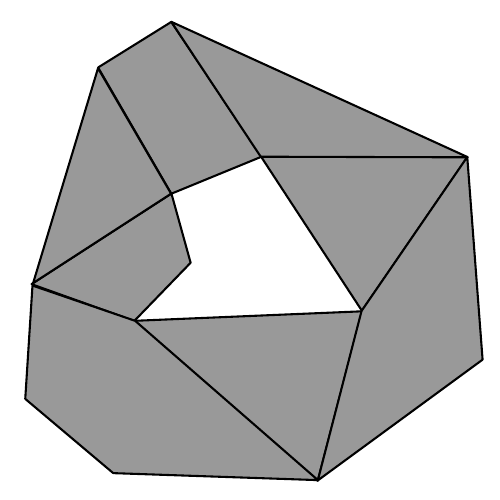}
    \includegraphics[width=.32\linewidth]{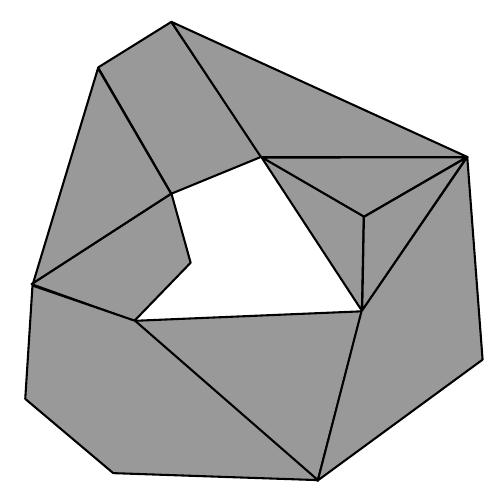}
    \caption{Partitions of a nonconvex region in the plane obtained by removing a central non-convex portion from a convex polyhedron. The first  partition does not satisfy the internal vertices condition \eqref{intvercon} (Left), the second partition admits a pairwise IB scheme (Center), and the third partition admits a 3-way IB scheme but not a pairwise one (Right).}
    \label{fig:partitioning}
\end{figure}

Consider a (nonconvex) bounded region in the plane $\Omega \subset \bbR^2$ that describes all possible locations for a UAV, as described in Section~\ref{obstacleavoidancesec}. Assume that
$\Omega$ can be partitioned into polyhedra $\{P^i\}_{i=1}^d$ such that $\bigcup_{i=1}^d P^i = \Omega$ and $\relint(P^i) \cap \relint(P^j) = \emptyset$ for each distinct $i,j \in \llbracket d \rrbracket$. We note that this partition will not, in general, be unique, and its selection can have a significant effect on questions of representability or formulation size. Figure~\ref{fig:partitioning} illustrates this for a convex region with a ``hole.'' The figure shows three ways to partition the resulting nonconvex region into convex polyhedra. Once this partition is fixed, we describe the associated polyhedra in $\scrV$-form, and so the corresponding combinatorial disjunctive constraint is given by $\scrS = \{\ext(P^i)\}_{i=1}^d$ and $J = \bigcup\{S \in \scrS\}$. We additionally forbid polyhedra with ``internal vertices'' by requiring that
\begin{equation}\label{intvercon}
v \in P^i \Longleftrightarrow v \in \ext(P^i) \quad \forall i \in \llbracket d \rrbracket, v \in J,
\end{equation}
so that $\scrS$ corresponds to the maximal elements of a polyhedral complex \cite[Section 5.1]{Ziegler:2007}. For example, the second and third partitions in Figure~\ref{fig:partitioning} satisfy this condition, while the first does not.

In this setting, minimal infeasible sets have a natural characterization.

\begin{theorem} \label{thm:partitioning}
    Take bounded $\Omega \subset \bbR^2$ and a polyhedral partition $\{P^i\}_{i=1}^d$ of $\Omega$ satisfying the internal vertex condition \eqref{intvercon}. If $\scrS = \{\ext(P^i)\}_{i=1}^d$, then $r\bra{H^c_\scrS}\leq 3$.
\end{theorem}

\proof{\textbf{Proof}}
    Take some minimal infeasible hyperedge $E \in \mathcal{E}_\scrS$ of $H^c_\scrS$, assuming for contradiction that $r = |E| > 3$, and label the points $E = \{v^i\}_{i=1}^r$. First, we show that the points may not be in general position, i.e. that w.l.o.g. $v^r \in \Conv(\{v^i\}_{i=1}^{r-1})$. Then, we argue that the points not being in general position implies that $\{v^i\}_{i=1}^{r-1}$ is also an infeasible set, violating the minimality condition.

    Assume for contradiction that the points are in general position; that is, that none can be written as a convex combination of the others. This implies that $\ext(\Conv(E)) = E$. Assume that the ordering $\{v^1,\ldots,v^r\}$ forms a path around the edges of $\Conv(E)$; that is, $v^i$ and $v^j$ both lie on an edge of $\Conv(E)$ if and only if $|i-j|=1$ or $\{i,j\}=\{1,r\}$.

    Choose some set $S^1 \in \scrS$ and some $2 < j < r$ such that $v^1,v^j \in S^1$ and $v^2 \not\in S^1$; the associated polyhedron is $P^1$. Such a set exists, else $E$ is not a minimal infeasible set (choose instead $E \backslash \{v^2\}$). Now choose $S^2 \in \scrS$ such that $v^2, v^r \in S^2$; the associated polyhedron is $P^2$. Such as set exists, as $\{v^2,v^r\} \subsetneq E$ and $E$ is minimal. As the nodes $v^1,v^2,v^j,v^r$ are interlaced along the boundary of $\Conv(E)$, we have that $\Conv(\{v^1,v^j\}) \cap \Conv(\{v^2,v^r\}) \subset \Conv(E)$ is nonempty. As each of the four points is on the boundary of $\Conv(E)$, and the points are in general position, it follows that $\Conv(\{v^1,v^j\}) \cap \Conv(\{v^2,v^r\}) = \relint(\Conv(\{v^1,v^j\})) \cap \relint(\Conv(\{v^2,v^r\}))$. Therefore, there must exist some point $y$ with $y \in \relint(\Conv(\{v^1,v^j\})) \subseteq \relint(P^1)$ and $y \in \relint(\Conv(\{v^2,v^r\})) \subseteq \relint(P^2)$. However, this implies that $\relint(P^1) \cap \relint(P^2) \neq \emptyset$, which contradicts the assumption that our sets partition the region $\Omega$.

    Finally, it just remains to show that $\{v^i\}_{i=1}^{r-1}$ is also an infeasible set, and therefore $\{v^i\}_{i=1}^r$ cannot be a minimal infeasible set. Assume for contradiction that it is not: i.e. that there exists some $j$ such that $\{v^i\}_{i=1}^{r-1} \subseteq \ext(P^j)$. But this implies that $v^r \in J$ and $v^r \in \Conv(E) \subseteq P^j$, yet $v^r \not\in \ext(P^j)$, a contradiction of the internal vertices assumption. \Halmos\endproof

In other words, every polyhedral partition of the plane is $3$-way independent branching-representable, and pairwise IB representability can be checked in polynomial time (for example, by enumerating the subsets of $J$ of cardinality 3). To illustrate, in Figure~\ref{fig:partitioning} we depict the three possible cases for a partition with respect to  Theorem~\ref{thm:partitioning}: 1) it does not satisfy the internal vertices condition, 2) it admits a pairwise IB scheme ($r(H^c_\scrS)=2$), or 3) it does not admit a pairwise IB scheme, but does admit a 3-way IB scheme ($r(H^c_\scrS)=3$).

Furthermore, we can argue that we can always represent a obstacle avoidance constraint in such a way that it admits a pairwise IB scheme. Inspecting Figure~\ref{fig:partitioning}, we see that the region $\Omega$ is the same in each, and it is only the partition of $\Omega$ that can potentially lead to constraints that are not pairwise IB-representable. Therefore, the obstacle avoidance constraint is invariant to the specification of the partition, and if any polyhedral partitioning exists, then it is always possible to construct one that satisfies the conditions of Theorem~\ref{thm:partitioning}.\footnote{Note that this result does not carry over to piecewise linear functions over $\Omega$, as the choice of the partition is intimately connected with the values the function may take.}

\subsection{SOS2} \label{ss:SOS2}

Our first example of a constraint that is always pairwise IB-representable is the SOS$2(N)$ constraint. Recall that $J = \llbracket N \rrbracket$ and $\scrS = \{\{\tau,\tau+1\}: \tau \in \llbracket N-1 \rrbracket\}$ for SOS2. Then $\mathcal{E}_\scrS=\set{\{\tau,\tau+t\}: \tau,\tau+t \in \llbracket N \rrbracket, \: t \geq 2}$, $r\bra{H^c_\scrS}=2$, and formulation \eqref{basicIB} has depth $t=\binom{N}{2}-N+1$.  However, Vielma and Nemhauser \cite{Vielma:2009a} construct a pairwise IB scheme for SOS$2(N)$ constraints of depth logarithmic in $N$. The construction is built around a Gray code \cite{Savage:1997}, or sequence of distinct binary vectors $\{h^i\}_{i=1}^{N-1} \subseteq \{0,1\}^{\lceil \log_2(N-1) \rceil}$ where each adjacent pair $(h^i,h^{i+1})$ differs in exactly one component. Notationally, here and throughout, take $h^0 \defeq h^1$ and $h^N \defeq h^{N-1}$. The pairwise IB scheme is then given by
\[
    L^j = \left\{\tau \in \llbracket N \rrbracket : h^{\tau-1}_j = 1 \text{ or } h^{\tau}_j = 1 \right\}, \quad R^j = \left\{\tau \in \llbracket N \rrbracket : h^{\tau-1}_j = 0 \text{ or } h^{\tau}_j = 0 \right\} \quad\quad \forall j \in \llbracket \lceil \log_2(N-1) \rceil \rrbracket.
\]

We observe that the resulting formulation matches the lower bound from Proposition~\ref{prop:log-bound} with respect to the number of auxiliary binary variables and is significantly smaller than formulation \eqref{basicIB}.

Indeed, formulation \eqref{basicIB} is likely to be unnecessarily large for the other constraints we consider as well, so we turn our attention to finding smaller IB schemes in Section~\ref{PISEC}. We end this section by applying Theorem~\ref{thm:representability} to succinctly prove the pairwise IB-representability of two other constraints.

\subsection{Other pairwise IB-representable constraints}

\paragraph{SOS$k$}

In this case $\mathcal{E}_\scrS=\set{\{\tau,\tau+t\}: \tau,\tau+t \in \llbracket N \rrbracket, \: t\geq k+1}$ and $r\bra{H^c_\scrS}=2$.

\paragraph{Grid triangulations}
We show that $r\bra{H^c_\scrS}=2$ by seeing that for any infeasible set $T \subseteq J$ there exist some distinct $v,w\in T$ such that $\set{v,w}$ is infeasible. Indeed, if there are some $v,w \in T$ such that $||v - w||_\infty > 1$, then there does not exist any triangle on the grid that contains both, so $\{v,w\}$ is also an infeasible set. Otherwise, we have that $T \subset \{r,r+1\}\times\{s,s+1\}$ for some $r,s$, and that $T$ contains elements in both of the triangles in this square. For each of the two triangles, we can select an element of $T$ that is not contained in the other triangle, which yields an infeasible pair contained in $T$. Therefore, any grid triangulation is pairwise IB-representable.

\paragraph{Discretization of multilinear terms}
Similarly as for the grid triangulation case, each infeasible set $T \subseteq J$ must necessarily contain two elements $v,w \in T$ with $||v - w||_\infty > 1$, and so we have that $T$ can be reduced to the infeasible pair $\{v,w\}$. Therefore, $r(H^c_\scrS)=2$, and so any discretization of this form is pairwise IB-representable.

\section{Pairwise independent branching schemes}\label{PISEC}
The pairwise independent branching scheme framework was initially introduced by Vielma and Nemhauser \cite{Vielma:2009a}, where it was used to model particularly structured piecewise linear functions. In the remainder of this work, we will focus on pairwise IB schemes and offer a complete picture of their expressive powers, along with an algorithmic framework for constructing them.

\subsection{Graphical representations of pairwise IB-representable CDCs} \label{ss:graph-representation}

From our covering assumption $J = \bigcup\set{S\in \scrS}$, we can see that $\abs{E}\geq 2$ for each $E\in\mathcal{E}_\scrS$. By applying Theorem~\ref{thm:representability}, we then immediately have that $\CDC(\scrS)$ is pairwise IB-representable if and only if $H^c_\scrS$ is (equivalent to) a graph. Along this line, for any constraint we may define a \emph{conflict graph} for any $\CDC(\scrS)$ as $G^c_\scrS \defeq (J,\bar{E})$, where $\bar{E} = \bar{E}_\scrS \defeq \{\set{u,v} \in [J]^2 : u \neq v, \: \{u,v\} \text{ is an infeasible set}\}$ is the set of all infeasible pairs of elements of $J$. Checking for pairwise IB-representability then reduces to verifying if $\mathcal{E}_\scrS=\bar{E}_\scrS $. The following corollary of Theorem~\ref{thm:representability} shows that this can also be verified by only working with $G^c_\scrS$.

\begin{corollary} \label{thm:cliques}
    $\CDC(\scrS)$ is pairwise IB-representable if and only if the sets $\scrS$ are exactly the maximal independent sets of $G_\scrS^c$.
\end{corollary}
\proof{\textbf{Proof}}
  If $\CDC(\scrS)$ is pairwise IB-representable, then $G_\scrS^c$ is equivalent to $H^c_\scrS$. By applying Theorem~\ref{thm:representability}, the maximal independent sets of $G_\scrS^c$ are exactly the elements of $\scrS$. For the converse, assume for a contradiction that $\scrS$ is exactly the maximal independent sets of $G_\scrS^c$, but that there exists some $E\in\mathcal{E}_\scrS$ with $\abs{E}\geq 3$. By the minimal infeasibility of $E$, we have that $\set{r,s}\notin \bar{E}_\scrS$ for any distinct $r,s\in E$, and therefore $E$ is an independent set in $G_\scrS^c$. This implies that $E$ is contained in a maximal independent $S$. By assumption, $S\in \scrS$, which contradicts the infeasibility of $E$.
\Halmos\endproof

Therefore, verifying general pairwise IB-representability reduces to enumerating the maximal independent sets of $G_\scrS^c$ and identifying them to exactly the sets $\scrS$. As an example, we can see that, for cardinality constraint of degree $\ell$ with $2 \leq \ell < |J|$, the only maximal independent set of $G_\scrS^c$ is the entire ground set $J$, which certainly cannot be identified with $\scrS = \{S \subset J : |S| = \ell \}$.

\subsection{Representation at a given depth}
Once a CDC has been shown to be pairwise IB-representable, a natural next question is: what is the smallest possible depth at which we may construct an IB scheme? That is, we ask if there exists a pairwise IB scheme for $\CDC(\scrS)$ of some given depth $t$. The answer to this question reduces to the existence of a graphical decomposition of the conflict graph $G_\scrS^c$.

\begin{definition}
    A \emph{biclique cover} of the graph $G = (J,E)$ is a collection of complete bipartite subgraphs $\left\{G^j\defeq (A^j\cup B^j,E^j)\right\}_{j=1}^t$ of $G$ that cover all the edges of $G$. Formally, this means that for each $j \in \llbracket t \rrbracket$, $\emptyset \subsetneq A^j,B^j \subsetneq J$, $A^j \cap B^j = \emptyset$, and $E^j =A^j * B^j\defeq \set{\set{a,b}\,:\, a\in A^j,\: b\in B^j}$, and that $\bigcup_{j=1}^t E^j = E$. 
    
    For notational simplicity, we will often refer to the sets $\{(A^j,B^j)\}_{j=1}^t$ as a biclique cover, as we can recover the graphs $G^j$ directly.
\end{definition}

The following theorem formalizes the equivalence between biclique covers and pairwise IB schemes.

\begin{theorem} \label{thm:biclique-IBSneq}
    If $\{(A^j,B^j)\}_{j=1}^t$ is  biclique cover of the conflict graph $G_\scrS^c$ for pairwise IB-representable $\CDC(\scrS)$, then a pairwise IB scheme for $\CDC(\scrS)$ is given by
    \begin{equation} \label{eqn:biclique-to-IBnew}
        L^j = J \backslash A^j, \quad R^j = J \backslash B^j \quad\quad \forall j \in \llbracket t \rrbracket.
    \end{equation}
    Conversely, if  $\{(L^j,R^j)\}_{j=1}^t$ is a pairwise IB scheme for $\CDC(\scrS)$, then a  biclique cover of the conflict graph $G_\scrS^c$ is given by
    \begin{equation} \label{eqn:biclique-to-IBnewconv}
        A^j = J \backslash L^j, \quad B^j = J \backslash R^j \quad\quad \forall j \in \llbracket t \rrbracket.
    \end{equation}
\end{theorem}
\proof{\textbf{Proof}}
For the first part, take $\bar{E}$ as the edge set of $G_\scrS^c$. To see that any $S \in \scrS$ is feasible for the IB scheme \eqref{eqn:biclique-to-IBnew}, note that if $S\not\subseteq L^j$ and $S\not\subseteq R^j$, then there exist some $u\in A^j \cap S$ and $v\in B^j \cap S$. However, this implies that $\set{u,v}\in A^j * B^j \subseteq \bar{E}$, which is a contradiction of feasibility as $\set{u,v}\subseteq S$ and $S \in \scrS$. Furthermore, as $\{(A^j,B^j)\}_{j=1}^t$ is a biclique cover of $G^c_\scrS$, for every $\set{u,v}\in \bar{E}$ we have that there exists some level $j\in \sidx{t}$ such that w.l.o.g. $u\in A^j$ and $v\in B^j$. This implies that $u\notin L^j$ and $v\notin R^j$ by their construction, and as $\CDC(\scrS)$ is pairwise IB-representable, then any infeasible set for $\CDC(\scrS)$ is also infeasible for the proposed IB scheme. Therefore, \eqref{eqn:biclique-to-IBnew} is a valid pairwise IB scheme.

For the second part, note that $A^j \cap B^j = \emptyset$ for all $j \in \llbracket t \rrbracket$, and that the covering portion of Assumption~\ref{basicassumption} implies that $L^j\cup R^j=J$. Therefore, it only remains to show that $\bar{E} = \bigcup_{j=1}^t \bar{E}^j$, where $\bar{E}^j = A^j * B^j$. For that, first note that as $L^j\cup R^j=J$, we have that $A^j=R^j\setminus L^j$ and $B^j=L^j\setminus R^j$. The containment $\bar{E} \subseteq \bigcup_{j=1}^t \bar{E}^j$ then follows by noting that, as $\{(L^j,R^j)\}_{j=1}^t$ is a valid pairwise IB scheme, each minimal infeasible set $\{u,v\}\in \bar{E}$ has some level $j \in \llbracket t \rrbracket $ such that $\{u,v\}\not\subseteq L^j$ and $\{u,v\}\not\subseteq R^j$. Then, as $L^j\cup R^j=J$, we have (w.l.o.g.) that  $u \in L^j \backslash R^j \equiv B^j$ and $b \in R^j \backslash L^j \equiv A^j$, and so $\{a,b\} \in \bar{E}^j$. For the reverse containment $\bigcup_{j=1}^t \bar{E}^j \subseteq \bar{E}$, take some arbitrary $j \in \llbracket t \rrbracket$ and some edge $\{a,b\} \in \bar{E}^j$. From the definition of our biclique cover, we have that w.l.o.g. $a \in A^j \equiv R^j \backslash L^j$ and $b \in B^j \equiv L^j \backslash R^j$. Therefore, $\{a,b\}$ is an infeasible set for the IB scheme, and thus for $\CDC(\scrS)$ as well, and so $\{a,b\} \in \bar{E}$.
\Halmos\endproof

We can now naturally frame the problem of finding a minimum depth pairwise IB scheme as the minimum biclique cover problem \cite{Fishburn:1996,Garey:1979}. Unfortunately, the decision version of this problem is known to be NP-complete \cite{Orlin:1977} and inapproximable within a factor of $|J|^{1/3-\epsilon}$ if $P \neq NP$~\cite{Gruber:2007}, even for bipartite graphs. However, we note that it is simple to construct a MIP feasibility problem for finding a pairwise IB scheme of a given depth $t$, which gives us a way to algorithmically find the smallest pairwise IB scheme for a specific (fixed) CDC. We present such a formulation in Proposition~\ref{prop:feasibility-IP} in Appendix~\ref{app:feasibility-IP}. Additionally, Cornaz and Fonlupt \citep{Cornaz:2006} present a MIP formulation (with an exponential number of constraints that can be efficiently separated) to find the minimum level biclique cover of a graph.

Furthermore, we can restate the MIP formulation from \cite{Vielma:2009a} (which is a special case of \eqref{eqn:multiway-formulation} with $k=2$) in terms of biclique covers of $G_\scrS^c$.

\begin{proposition}[Theorem 5, \cite{Vielma:2009a}; Theorem 1, \cite{Vielma:2010}] \label{prop:ideal-formulation-for-IBS}
    If $\CDC(\scrS)$ is pairwise independent branching-representable and $\{(A^j,B^j)\}_{j=1}^t$ is a biclique cover for $G_\scrS^c$, then the following is an ideal formulation for $\CDC(\scrS)$:
    \begin{subequations} \label{eqn:ideal-formulation-for-IBS}
        \begin{alignat}{2}
        \sum_{v \in A^j} \lambda_v &\leq z_j &\quad \forall j \in \llbracket t \rrbracket \label{eqn:ideal-formulation-for-IBS-1} \\
        \sum_{v \in B^j} \lambda_v &\leq 1 - z_j &\quad \forall j \in \llbracket t \rrbracket \label{eqn:ideal-formulation-for-IBS-2} \\
        (\lambda,z) &\in \Delta^J \times \{0,1\}^t. \label{eqn:ideal-formulation-for-IBS-3}
    \end{alignat}
    \end{subequations}
\end{proposition}

We end the section by noting that the relation between biclique covers and independent sets has also been exploited in the study of boolean functions, particularly in the equivalence between posiforms and maximum weighted stable sets  (e.g. \cite[Theorem 13.16]{crama2011boolean}). In fact, formulation \eqref{eqn:ideal-formulation-for-IBS} is reminiscent of formulation (13.45--13.50) in \cite[Theorem 13.13]{crama2011boolean}. The main difference between these formulations is that in the context of \cite{crama2011boolean} the $\lambda$ variables will be binary variables not constrained to lie in the unit simplex. For this reason inequalities (\ref{eqn:ideal-formulation-for-IBS-1}--\ref{eqn:ideal-formulation-for-IBS-2}) appear disaggregated in \cite[Theorem 13.13]{crama2011boolean} in the form $\lambda_v \leq z_j$ for all $v\in A^j$, $j \in \llbracket t \rrbracket$. However, the resulting formulation is not ideal (See \cite[Section 5]{Vielma:2009a} for more details). Still, the combinatorial aspects of this connection could prove useful for constructing small IB schemes.

In the next section, we will explore instances where we can, in closed form, construct small (asymptotically optimal) IB schemes for families of particularly structured CDCs.

\section{Illustrative examples} \label{s:examples}
With a framework to construct pairwise independent branching schemes for arbitrary pairwise IB-representable CDCs, we now return to some of our motivating examples. We will apply our methodology to these specific structures, and produce small, closed-form IB schemes. In particular, this allows us to construct novel, small MIP formulations for these constraints.

\subsection{A simple IB scheme and its limitations}
To start, we show that \emph{any} pairwise IB-representable CDC admits an IB scheme of depth $|J|$. If $|J|$ is smaller than $|\scrS|$, this already offers a drop in size from \eqref{eqn:disaggregated}. This IB scheme covers all edges incident to node with the simple biclique corresponding to the star centered at that node.

\begin{figure}[htpb]
    \centering
    \includegraphics[width=.37\linewidth]{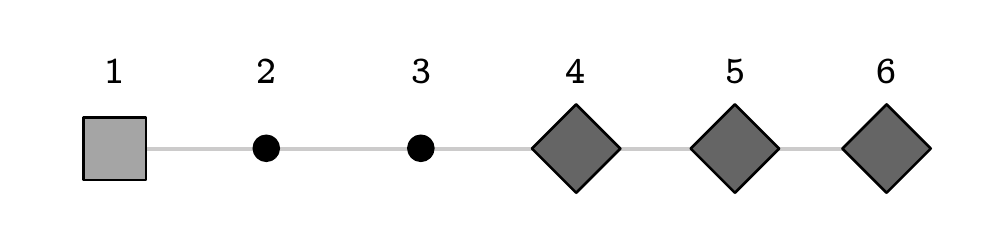} \hspace{2em}
    \includegraphics[width=.55\linewidth]{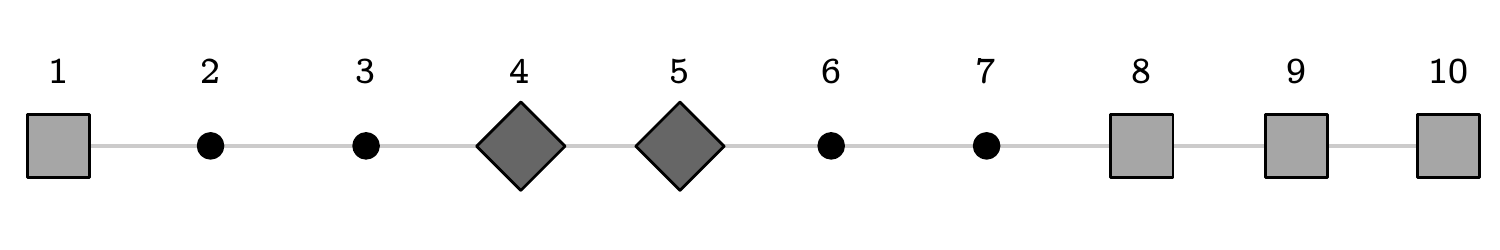} \\
    \includegraphics[width=.37\linewidth]{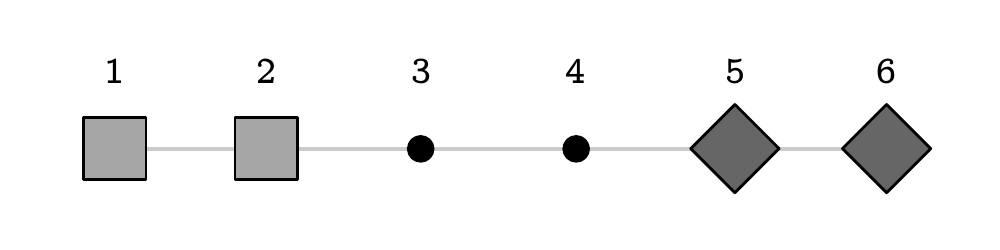} \hspace{2em}
    \includegraphics[width=.55\linewidth]{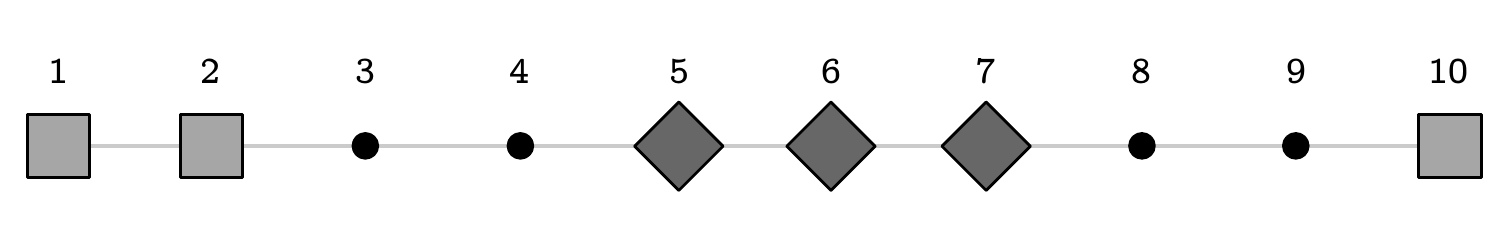} \\
    \includegraphics[width=.37\linewidth]{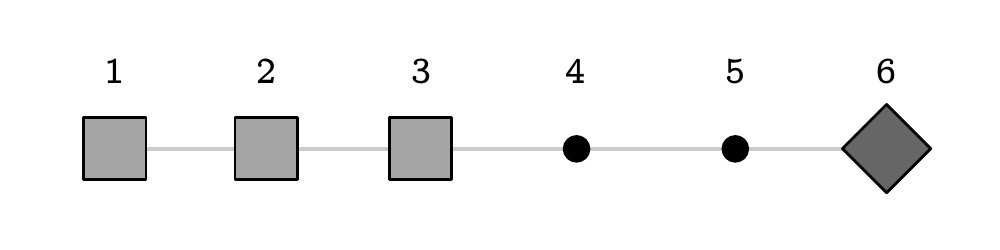} \hspace{2em}
    \includegraphics[width=.55\linewidth]{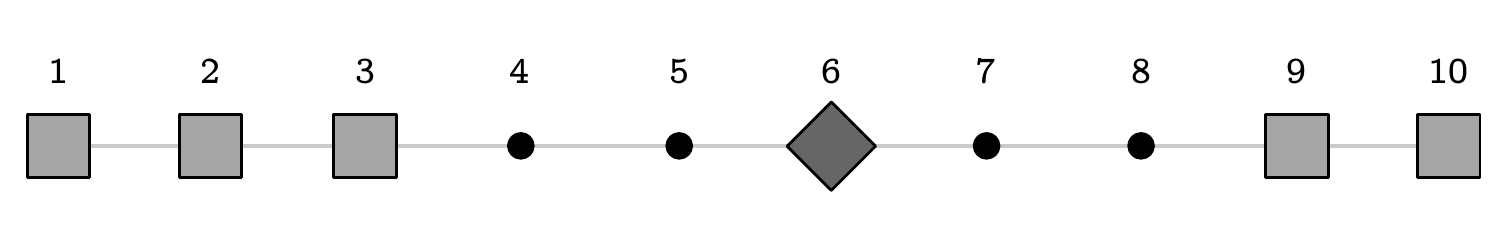} \\
    \includegraphics[width=.37\linewidth]{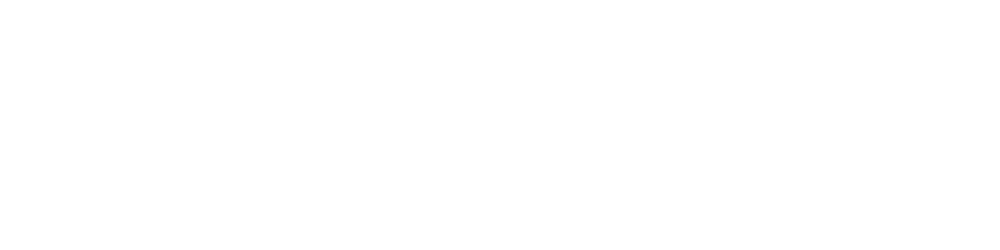} \hspace{2em}
    \includegraphics[width=.55\linewidth]{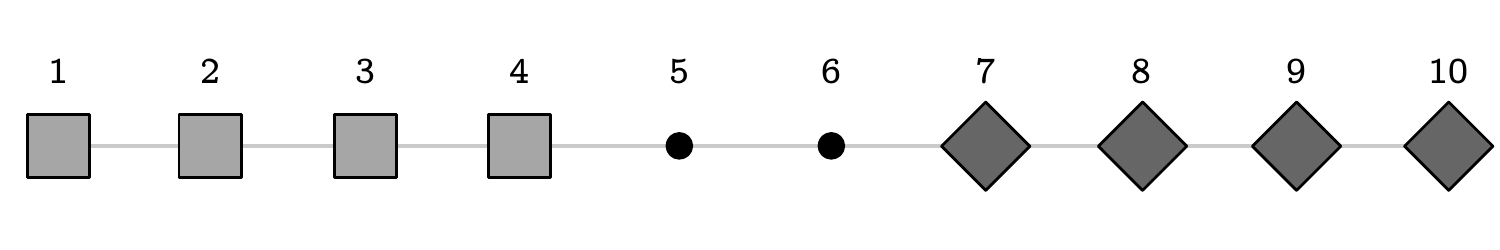}
    \caption{Visualizations of the biclique covers presented in the text for SOS$3(6)$ (Left) and SOS$3(10)$ (Right). Each row corresponds to some level $j$, and the elements of $A^j$ and $B^j$ are the squares and diamonds, respectively.}
    \label{fig:sos3}
\end{figure}

\begin{proposition}[Covering with Stars] \label{prop:simple-pIBS}
    For pairwise IB-representable $\CDC(\scrS)$, a biclique cover for $G_\scrS^c$ is given by:
    \begin{subequations} \label{eqn:simple-pIBS}
    \[
        A^v = \{v\},\quad\quad
        B^v = \left\{u \in J : \set{u,v} \in \bar{E}\right\} \quad\quad \forall v \in J.
    \]
    \end{subequations}
\end{proposition}
\proof{\textbf{Proof}}
    By construction of the sets, we see that each $\set{r,s} \in \bar{E}^v \equiv A^v * B^v$ corresponds to an infeasible edge: that is, $\bar{E}^v \subseteq \bar{E}$ for each $v$, and so $\bigcup_{v \in J} \bar{E}^v \subseteq \bar{E}$. Furthermore, each infeasible edge $\set{r,s} \in \bar{E}$ is infeasible for levels $r$ and $s$, and so $\bar{E} \subseteq \bigcup_{v \in J} \bar{E}^v$. Therefore, this construction forms a valid biclique cover of the conflict graph.
\Halmos\endproof
This gives us an upper bound of $|J|$ on the minimum depth for any pairwise IB-representable CDC. However, if we exploit the specific structure of a CDC, we can typically get much smaller formulations. For instance, consider the following two instances of the SOS$3(N)$ constraint for small values of $N$. First, consider the instance with $N=6$, where $\abs{J}=6$ and $\scrS = \{\{1,2,3\},\{2,3,4\},\{3,4,5\},\{4,5,6\}\}$. Therefore, $|\scrS|=4$, yielding a lower bound of depth $\log_2(4)=2$ from Proposition~\ref{prop:log-bound}. However, there does not exist a biclique cover of depth 2 (which can be verified via Proposition \ref{prop:feasibility-IP}), though one of depth 3 does exist:
\begin{alignat*}{3}
    A^1 &= \{1\}, \quad &B^1 &= \{4,5,6\} \\
    A^2 &= \{1,2\}, \quad &B^2 &= \{5,6\}   \\
    A^3 &= \{1,2,3\}, \quad &B^3 &= \{6\}.
\end{alignat*}
We can see the proposed IB scheme on the left side of Figure \ref{fig:sos3}. For clarity, the associated MIP formulation for the CDC from Proposition~\ref{prop:ideal-formulation-for-IBS} is
\begin{alignat*}{2}
    \lambda_1 &\leq z_1 \quad\quad \lambda_4 + \lambda_5 + \lambda_6 &\leq 1 - z_1 \\
    \lambda_1 + \lambda_2 &\leq z_2 \quad\quad \lambda_5 + \lambda_6 &\leq 1 - z_2 \\
    \lambda_1 + \lambda_2 + \lambda_3 &\leq z_3 \quad\quad \lambda_6 &\leq 1 - z_3 \\
    (\lambda,z) &\in \Delta^6 \times \{0,1\}^3.
\end{alignat*}

Next, we consider $N=10$, where we also cannot attain the $\log_2(8)=3$ lower bound. However, a biclique for this the conflict graph of this constraint is
\begin{subequations} \label{eqn:sos3-6-biclique}
\begin{alignat}{3}
    A^1 &= \{1,8,9,10\}, \quad &B^1 &= \{4,5\} \\
    A^2 &= \{1,2,10\}, \quad &B^2 &= \{5,6,7\}   \\
    A^3 &= \{1,2,3,9,10\}, \quad &B^3 &= \{6\} \\
    A^4 &= \{1,2,3,4\}, \quad &B^4 &= \{7,8,9,10\},
\end{alignat}
\end{subequations}
as seen on the right side of Figure \ref{fig:sos3}. The corresponding MIP formulation is
\begin{alignat*}{2}
    \lambda_1 + \lambda_8 + \lambda_9 + \lambda_{10} &\leq z_1 \quad\quad \lambda_4 + \lambda_5 &\leq 1 - z_1 \\
    \lambda_1 + \lambda_2 + \lambda_{10} &\leq z_2 \quad\quad \lambda_5 + \lambda_6 + \lambda_7 &\leq 1 - z_2 \\
    \lambda_1 + \lambda_2 + \lambda_3 + \lambda_9 + \lambda_{10} &\leq z_3 \quad\quad \lambda_6  &\leq 1 - z_3 \\
    \lambda_1 + \lambda_2 + \lambda_3 + \lambda_4 &\leq z_4 \quad\quad \lambda_7 + \lambda_8 + \lambda_9 + \lambda_{10} &\leq 1 - z_4 \\
    (\lambda,z) &\in \Delta^{10} \times \{0,1\}^4.
\end{alignat*}

The ad-hoc construction for SOS$3(6)$ suggests a more general construction for SOS$k(N)$ when $k \leq N/2$ (assume for convenience that $N$ is even). Consider the sets given by
\[
    A^j = \{1,\ldots,j\} \cup \{j+N/2+k,\ldots,N\}, \quad\quad\quad
    B^j = \{j+k,\ldots,j+N/2\}
\]
for each $j \in \llbracket N/2 \rrbracket$. It is straightforward to see that this yields a biclique cover of the conflict graph for SOS$k(N)$ of depth $N/2$. Therefore, with this simple operation, we have constructed an ideal formulation for SOS$k(N)$ with size strictly smaller than $N$, the size of the na\"ive non-extended formulation \eqref{eqn:disaggregated}.

Based on the second example~\eqref{eqn:sos3-6-biclique}, we know that this construction is, in general, not the smallest possible. In Section~\ref{ss:sosk}, we will see how we can systematically construct small biclique covers (and MIP formulations) for SOS$k$($N$) with arbitrary $k$ and $N$, using techniques we will now develop.

\subsection{Systematic construction of biclique covers}

As discussed in Section~\ref{ss:SOS2}, there exists an IB scheme for the SOS2 constraint of optimal depth that can be constructed using a Gray code. The following proposition shows how the validity of this scheme can easily be proven by reinterpreting it via a biclique cover.
\begin{proposition} \label{prop:sos2-IB-scheme}
    Take a Gray code $\{h^i\}_{i=1}^{N-1} \subseteq \{0,1\}^{\lceil \log_2(N-1) \rceil}$ and let $h^0 \defeq h^1$ and $h^N \defeq h^{N-1}$. If $G_\scrS^c$ is the conflict graph of SOS2($N$), then a biclique cover for $G_\scrS^c$ of depth $\lceil \log_2(N-1) \rceil$ is given  by
    \begin{equation}\label{SOS2biclique}
        A^{j} = \left\{\tau \in \llbracket N \rrbracket : h^{\tau-1}_j =  h^{\tau}_j  = 0 \right\},\quad  B^{j} = \left\{\tau \in \llbracket N \rrbracket :  h^{\tau-1}_j =  h^{\tau}_j  = 1 \right\} \quad \forall j \in \llbracket \lceil \log_2(N-1) \rceil \rrbracket.
      \end{equation}
\end{proposition}
\proof{\textbf{Proof}}
For the SOS2($N$) constraint we have that $\bar{E}_\scrS=\set{\{r,s\} \in \llbracket N \rrbracket^2 : r+2 \leq s}$. Take any infeasible pair $\{r,s\}\in \bar{E}_\scrS$. As $r + 2 \leq s$, we conclude that $r-1<r<s-1<s$, and so it must be that $h^{r-1}, h^r \neq h^{s-1},h^s$. The set of components which flip values between the two pairs of adjacent codes $(h^{r-1},h^r)$ and $(h^{s-1},h^s)$ is $I = \{ j \in \llbracket \lceil \log_2(N-1) \rceil \rrbracket : h^{r-1}_j \neq h^r_j \text{ or } h^{s-1}_j \neq h^s_j\}$, and $|I| \leq 2$ as we have selected a Gray code. Now it must be the case that there is some component $j \in \llbracket \lceil \log_2(N-1) \rceil \rrbracket \backslash I$ wherein $h^{r-1}_j=h^r_j \neq h^{s-1}_j=h^s_j$, else we conclude that two of the vectors $h^i=h^{\ell}$ coincide for some $i \in \{r-1,r\}$ and $\ell \in \{s-1,s\}$, a contradiction of their uniqueness. Then $\{r,s\} \in E^j$, i.e. it is covered by the $j$-th level of the biclique. Furthermore, we observe that no edges of the form $\{r,r+1\}$ will be contained in the biclique cover, as it is not possible that $h_j^{r-1} = h_j^r = 0$ (resp. $=1$) and $h_j^r = h_j^{r+1} = 1$ (resp. $=0$) simultaneously.
\Halmos\endproof

\begin{figure}[htpb]
    \centering
    \includegraphics[width=.57\linewidth]{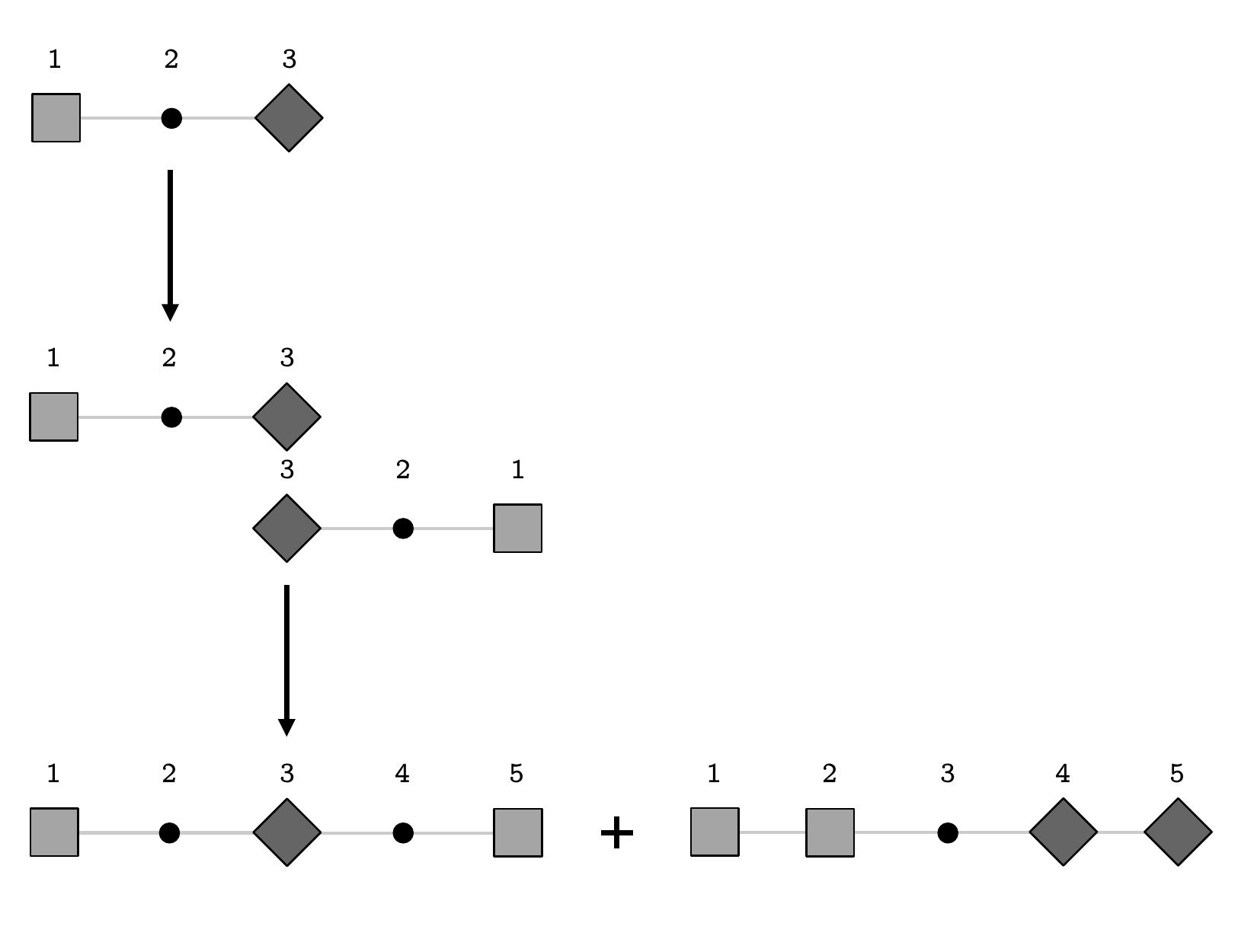}
    \caption{The recursive construction for biclique covers for SOS2. The first row is a single biclique that covers the conflict graph for SOS2(3) ($A^{1,1}$ are squares, $B^{1,1}$ are diamonds). The second row shows the construction which duplicates the ground set $\{1,2,3\}$ and inverts the ordering on the second copy. The third row shows the identification of the nodes that yields a valid biclique for SOS2(5). This biclique is then combined with a second that covers all edges between nodes identified with the first copy and those identified with the second, giving a biclique cover for SOS2(5) with two levels.}
    \label{fig:sos2-recursive}
\end{figure}

Interestingly, we can also view this construction recursively if we use a specific Gray code known as the binary reflected Gray code~\cite{Savage:1997}. For SOS2($2^k$), we will take $E^k$ as the edge set for the corresponding conflict graph. First, with $k=1$, $d=2^k=2$, $N=2^{k}+1=3$, then $E^1 = \{\{1,3\}\}$. A complete biclique cover is given by the single biclique $A^{1,1}=\{1\}$ and $B^{1,1} = \{3\}$. As we see in Figure~\ref{fig:sos2-recursive}, we can construct a biclique cover for SOS2 with $k=2$ by stitching together two copies of the biclique. We construct two copies of the node set for $k=1$, invert the second, and identify the last node from the first set with the first node with the second set. Then we can readily construct a mapping of the biclique $(A^{1,1},B^{1,1})$ for $k=1$ to a biclique for $k=2$, using the node identification, as $A^{2,1} = \{1,5\}$ and $B^{2,1} = \{3\}$. This will cover all edges in $E^2$ with both incident nodes in the first half of the nodes, or both in the second half of the nodes (along with some other edges in $E^2$, as well). To cover all edges with one adjacent node in the first half, and the other in the second half, we construct a second biclique of the form $A^{2,2} = \{1,2\}$ and $B^{2,2} = \{4,5\}$.

We can repeat this construction with $k=3$ to get the three level biclique cover
\begin{alignat*}{3}
    A^{3,1} &= \{1,5,9\}, \quad\quad &B^{3,1} &= \{3,7\} \\
    A^{3,2} &= \{1,2,8,9\}, \quad\quad &B^{3,2} &= \{4,5,6\} \\
    A^{3,3} &= \{1,2,3,4\}, \quad\quad &B^{3,2} &= \{6,7,8,9\}.
\end{alignat*}
Iterating this construction gives a biclique cover for $E^{k+1}$ as $\{(A^{k+1,i},B^{k+1,i})\}_{i=1}^{k+1}$, where
\begin{alignat*}{3}
    A^{k+1,i} &= \bigcup_{u \in A^{k,i}} \left\{u, 2^{k+1}+2-u\right\}, \quad\quad &B^{k+1,i} =& \bigcup_{v \in B^{k,i}} \left\{v,2^{k+1}+2-v\right\} \quad \forall i \in \llbracket k \rrbracket \\
    A^{k+1,k+1} &= \left\{1,\ldots,2^k\right\}, \quad\quad &B^{k+1,k+1} =& \left\{ 2^k + 2, \ldots, 2^{k+1}+1 \right\}.
\end{alignat*}
In fact, we can readily state this recursive construction in a more general form, where we adapt a biclique cover for one graph into a biclique cover for another graph that is created in some specific way.
\begin{lemma}\label{lem:recursive-cover}
    Take some graph $G=(\llbracket m+1 \rrbracket,E)$, and define $G^2 = (\llbracket 2m+1 \rrbracket, E^2)$ for
    \[E^2 =E \cup \bigl\{\{2m+2-u,2m+2-v\}\,:\, \set{u,v}\in E\bigr\} \cup\bigl(\sidx{m}*\sidx{m+2,2m+1}\bigr) \]
    where $\sidx{a,b}\defeq \set{a,\ldots,b}$.

    If $\{(A^j,B^j)\}_{j=1}^t$ is a biclique cover of $G$, then $\{(\tilde{A}^j,\tilde{B}^j)\}_{j=1}^{t+1}$ is a biclique cover of $G^2$, where
    \begin{alignat*}{3}
        \tilde{A}^j &= \bigcup_{u \in A^j} \{u, 2m+2-u\}, \quad\quad &\tilde{B}^j =& \bigcup_{v \in B^j} \{v, 2m+2-v\} \quad \forall j \in \llbracket t \rrbracket \\
        \tilde{A}^{t+1} &= \{1,\ldots,m\}, \quad\quad &\tilde{B}^{t+1} =& \{m+2,\ldots,2m+1\}.
    \end{alignat*}
\end{lemma}
In the remainder of this work, we will see how we may apply similar graphical results to systematically construct small biclique covers for the conflict graphs of constraints by exploiting their specific structure.

\subsection{Biclique covers for graph products and discretizations of multilinear terms}

Consider the discretization of multilinear terms described in Section~\ref{multilinearsec}, given by $J=\prod_{i=1}^\eta \llbracket d_i \rrbracket$ and $\scrS=\set{\prod_{i=1}^\eta \set{k_i,k_i+1} : k \in \prod_{i=1}^\eta \llbracket d_i -1 \rrbracket}$. We can interpret this constraint as a $\eta$-dimensional version of the SOS2 constraint, or as the Cartesian product of $\eta$ SOS2 constraints. This can be formalized through the following definition and straightforward lemma.

\begin{definition}
    The \emph{(disjunctive) graph product} of a family of graphs $\{G^i = (J^i,E^i)\}_{i=1}^\eta$ is $\bigvee_{i=1}^\eta G^i \defeq (J_P,E_P)$, where $J_P = \prod_{i=1}^\eta J^i$ and
    \[
        E_P = \left\{ \{u,v\} \in J_P * J_P : \exists i \in \llbracket \eta \rrbracket \text{ s.t. } \{u_i,v_i\} \in E^i \right\}.
    \]
\end{definition}

\begin{lemma} \label{lemma:graph-productmulti}
 Let $J=\prod_{i=1}^\eta \llbracket d_i \rrbracket$ and  $\scrS=\set{\prod_{i=1}^\eta \set{k_i,k_i+1} : k \in \prod_{i=1}^\eta \llbracket d_i -1 \rrbracket}$ be a $\eta$-dimensional discretization of multilinear terms, and $G_\scrS^c$ be the corresponding conflict graph. If $G^i$ is the conflict graph of SOS2($d_i$) for each $i\in \llbracket \eta \rrbracket$, then $G_\scrS^c= \bigvee_{i=1}^\eta G^i$.
\end{lemma}

Using this characterization, we can easily construct an IB scheme for discretizations of multilinear terms by taking the graph products of IB schemes for the SOS2 constraint.

\begin{lemma} \label{lemma:graph-product}
    Take a family of graphs $\{G^i = (J^i,E^i)\}_{i=1}^\eta$, and a biclique cover $\{(\tilde{A}^{i,j},\tilde{B}^{i,j})\}_{j=1}^{t_i}$ for each $G^i$. Then a biclique cover for $\bigvee_{i=1}^\eta G^i$ is given by $\bigcup_{i=1}^\eta \{(A^{i,j},B^{i,j})\}_{j=1}^{t_i}$, where
    \[
        A^{i,j} = \left(\prod_{\ell=1}^{i-1} J^\ell\right) \times \tilde{A}^{i,j} \times \left(\prod_{\ell=i+1}^\eta J^\ell \right), \quad\quad  B^{i,j} = \left(\prod_{\ell=1}^{i-1} J^\ell\right) \times \tilde{B}^{i,j} \times \left(\prod_{\ell=i+1}^\eta J^\ell \right) \quad\quad \forall i \in \llbracket \eta \rrbracket, \: j \in \llbracket t_i \rrbracket.
    \]
\end{lemma}

\begin{corollary}\label{prop:multilinear-IB-scheme}
Let $J=\prod_{i=1}^\eta \llbracket d_i \rrbracket$ and  $\scrS=\set{\prod_{i=1}^\eta \set{k_i,k_i+1} : k \in \prod_{i=1}^\eta \llbracket d_i -1 \rrbracket}$ describe a $\eta$-dimensional discretization of multilinear terms, and take $G_\scrS^c$ as its conflict graph. If for each $i \in \llbracket \eta \rrbracket$ we have a biclique cover $\{(\tilde{A}^{i,j},\tilde{B}^{i,j})\}_{j=1}^{t_i}$ for the conflict graph of SOS2$(d_i)$, then a biclique cover for $G_\scrS^c$ of depth $\sum_{i=1}^\eta t_i$ is given by $\bigcup_{i=1}^\eta \{(A^{i,j},B^{i,j})\}_{j=1}^{t_i}$, where
 \[
        A^{i,j} = \left\{x \in J : x_i \in \tilde{A}^{i,j} \right\}, \quad
        B^{i,j} = \left\{x \in J : x_i \in \tilde{B}^{i,j} \right\} \quad\quad \forall i \in \llbracket \eta \rrbracket, j \in \llbracket t_i \rrbracket.
 \]

    In particular, if we take $\{h^{i,j}\}_{j=1}^{d_i-1} \subseteq \{0,1\}^{\lceil \log_2(d_i-1) \rceil}$ as a Gray code for each $i \in \llbracket \eta \rrbracket$, where $h^{i,0} \defeq h^{i,1}$ and $h^{i,d_i} \defeq h^{i,d_i-1}$, then a biclique cover for $G_\scrS^c$ of depth $\sum_{i=1}^\eta \lceil \log_2(d_i-1) \rceil$ is given by:
    \begin{align*}
        \hat{A}^{i,j} = \left\{x \in J : \exists \gamma \text{ s.t. } x_i = \gamma, \: h^{i,\gamma-1}_j = h^{i,\gamma}_j = 0 \right\},\quad
        \hat{B}^{i,j} = \left\{x \in J : \exists \gamma \text{ s.t. } x_i = \gamma, \: h^{i,\gamma-1}_j = h^{i,\gamma}_j = 1 \right\}
    \end{align*}
    for each $i \in \llbracket \eta \rrbracket$ and $j \in \llbracket \lceil \log_2(d_i-1) \rceil \rrbracket$.
\end{corollary}

We note that, since $|\scrS| = \prod_{i=1}^\eta (d_i-1)$, by Proposition \ref{prop:log-bound} this construction yields a formulation that is asymptotically optimal (with respect to number of auxiliary binary variables) for any possible MIP formulation, up to an additive factor of at most $\eta$.

Furthermore, we can specialize this to the bilinear case studied by Misener et al. \cite{Misener:2011}.
\begin{corollary}
    There exists a biclique cover for a grid discretization of a bilinear function ($\eta=2$) with $d_1=m+1$ and $d_2=1$ of depth $\lceil \log_2(m) \rceil$.
\end{corollary}
This result yields an ideal MIP formulation for the outer-approximation of bilinear terms with $\lceil \log_2(m) \rceil$ auxiliary binary variables, $2(m+1)$ auxiliary continuous variables (the $\lambda$ variables, one for element in $J$), and $2\lceil \log_2(m) \rceil$ general inequality constraints. In contrast, the logarithmic formulation from Misener et al. \cite{Misener:2011} has $\lceil \log_2(m) \rceil$ auxiliary binary variables, $2\lceil \log_2(m) \rceil + 1$ auxiliary continuous variables, at least $2\lceil \log_2(m) \rceil + 6$ general inequality constraints, and is not ideal in general (see Appendix~\ref{app:misener}). Therefore, we gain an ideal formulation with a naturally induced constraint branching at the price of a modest number of additional auxiliary continuous variables. Furthermore, our formulation generalizes readily to a discretization along the second dimension ($d_2 > 1$), for non-uniform discretizations, and for higher dimensional multilinear functions ($\eta > 2$).

\subsection{Completing biclique covers via graph unions}
Another useful graphical technique for our heuristic constructions will be to combine together biclique covers, each of which is designed to cover a substructure of the constraint. For example, the conflict graph of a grid triangulation of the plane is equivalent to the conflict graph of a $2$-dimensional grid discretization of multilinear terms, with one extra edge added for each subrectangle in the grid. Therefore, a biclique cover of a grid triangulation can be obtained from a biclique cover of a $2$-dimensional discretization of multilinear terms (i.e. from Corollary~\ref{prop:multilinear-IB-scheme}) by completing it with some number additional bicliques that cover those extra edges. This construction can be formalized in the following way.
\begin{definition}
    The \emph{graph union} of a family of graphs $\set{G^i=(J^i,E^i)}_{i=1}^\eta$ is $\bigcup_{i=1}^\eta G^i\defeq (J_U,E_U)$, where $J_U=\bigcup_{i=1}^\eta J^i$ and $E_U=\bigcup_{i=1}^\eta E^i$.
\end{definition}

\begin{lemma} \label{lemma:graph-union}
    Take a family of graphs $\{G^i = (J^i,E^i)\}_{i=1}^\eta$ and a corresponding biclique cover $\{(A^j_i,B^j_i)\}_{j=1}^{t_i}$ of $G^i$ for each $i \in \llbracket \eta \rrbracket$. Then $\bigcup_{i=1}^\eta \{(A^j_i,B^j_i)\}_{j=1}^{t_i}$ is a biclique cover of $\bigcup_{i=1}^\eta G^i$.
\end{lemma}

\begin{figure}[htpb]
    \includegraphics[width=.185\linewidth]{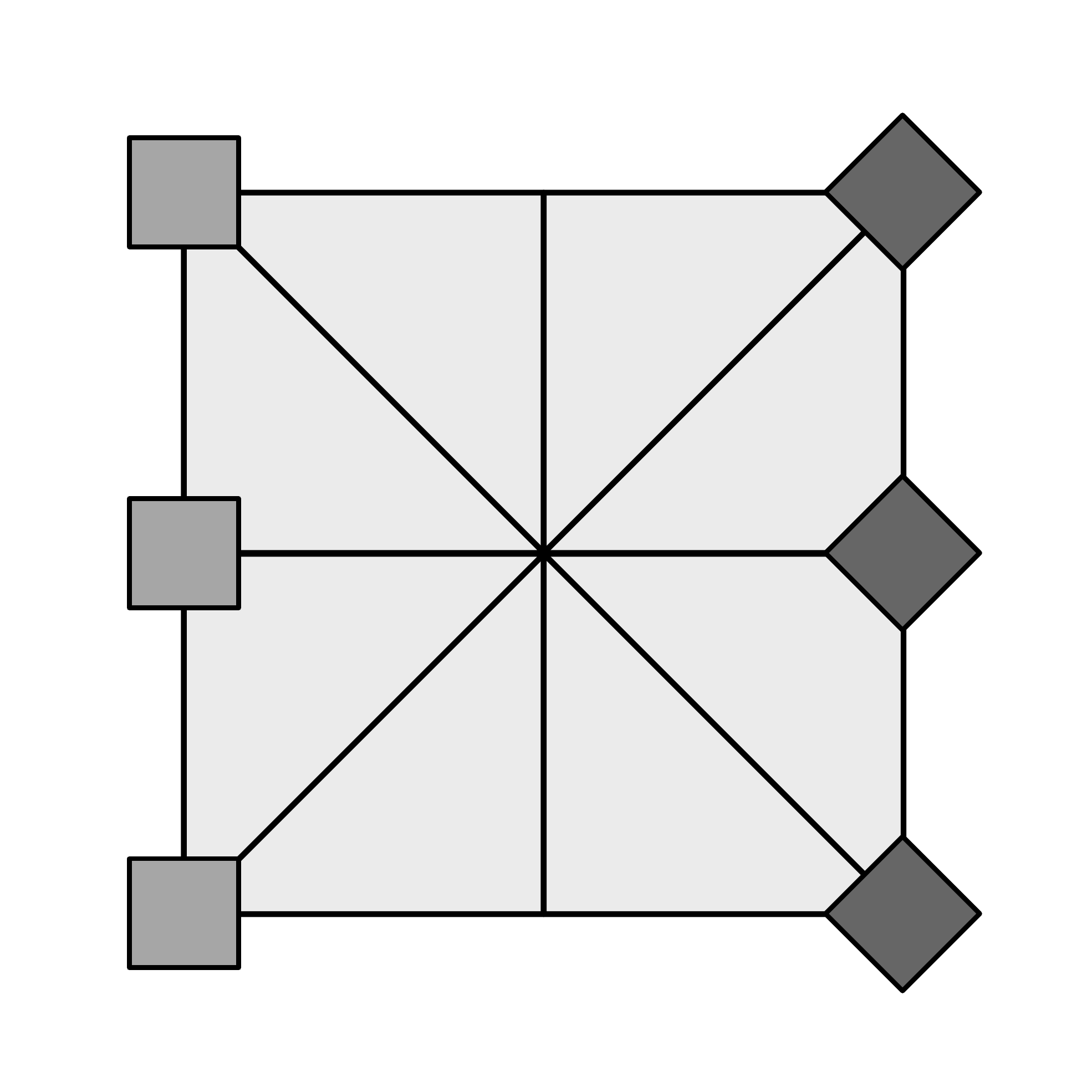}
    \includegraphics[width=.185\linewidth]{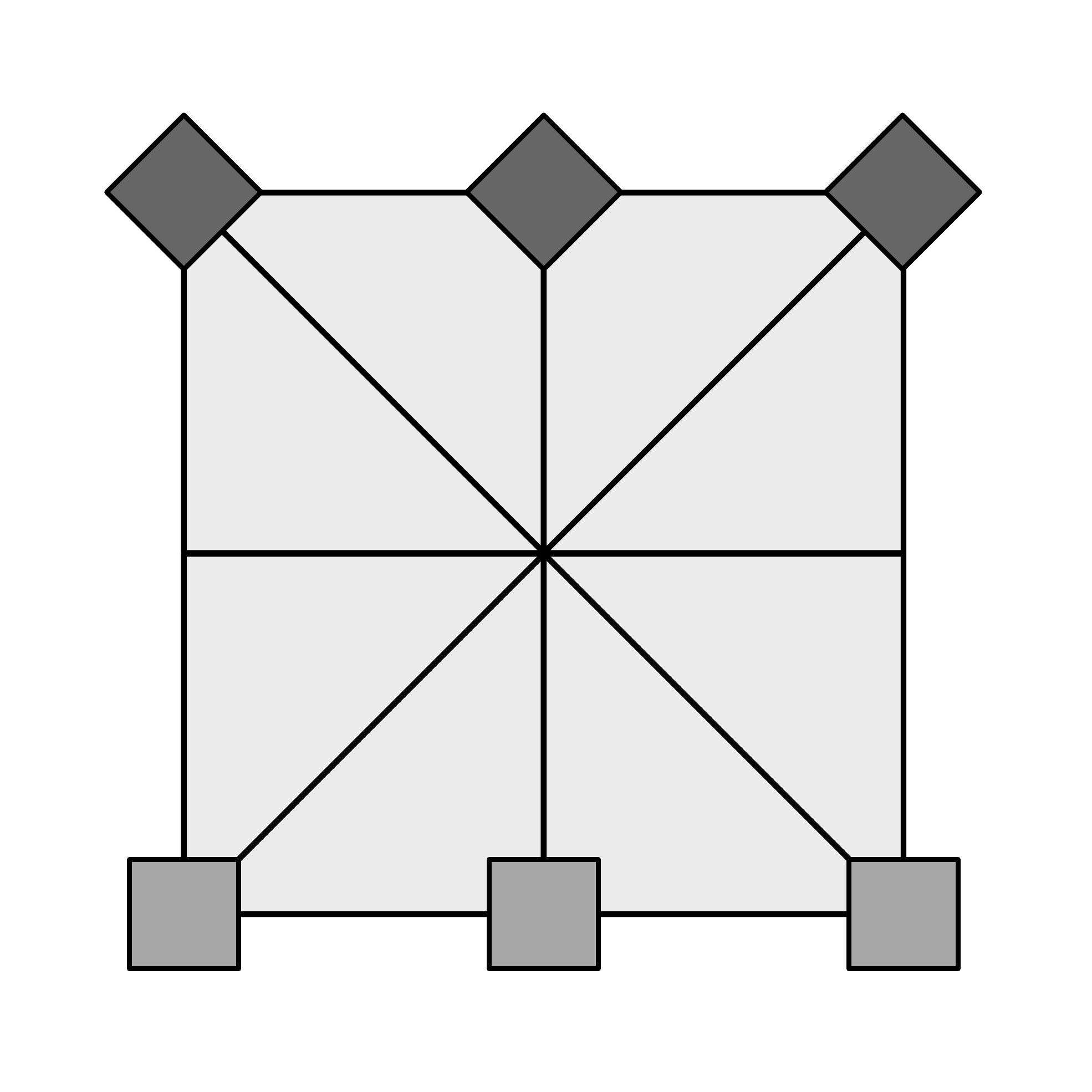}
    \includegraphics[width=.185\linewidth]{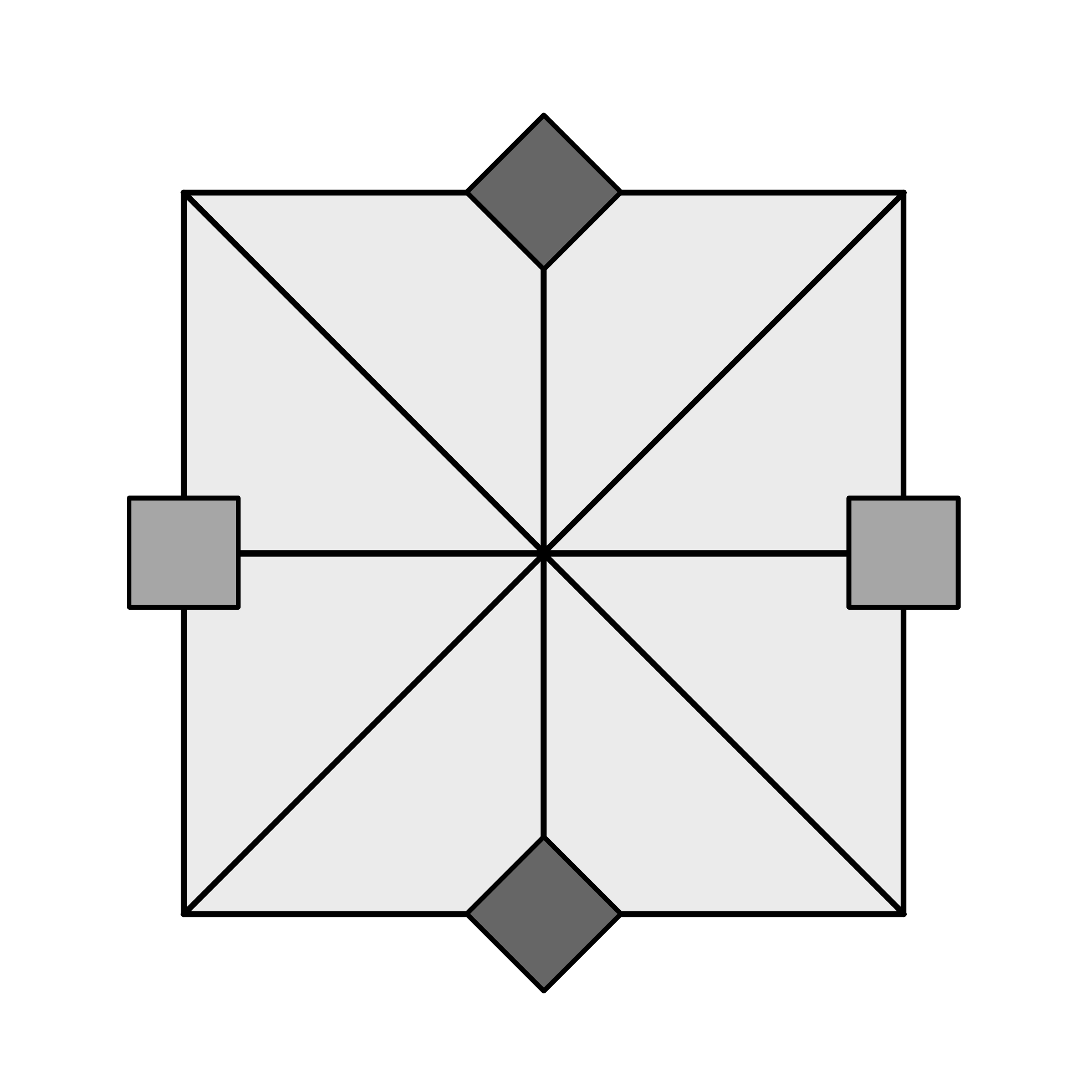} \\
    \includegraphics[width=.185\linewidth]{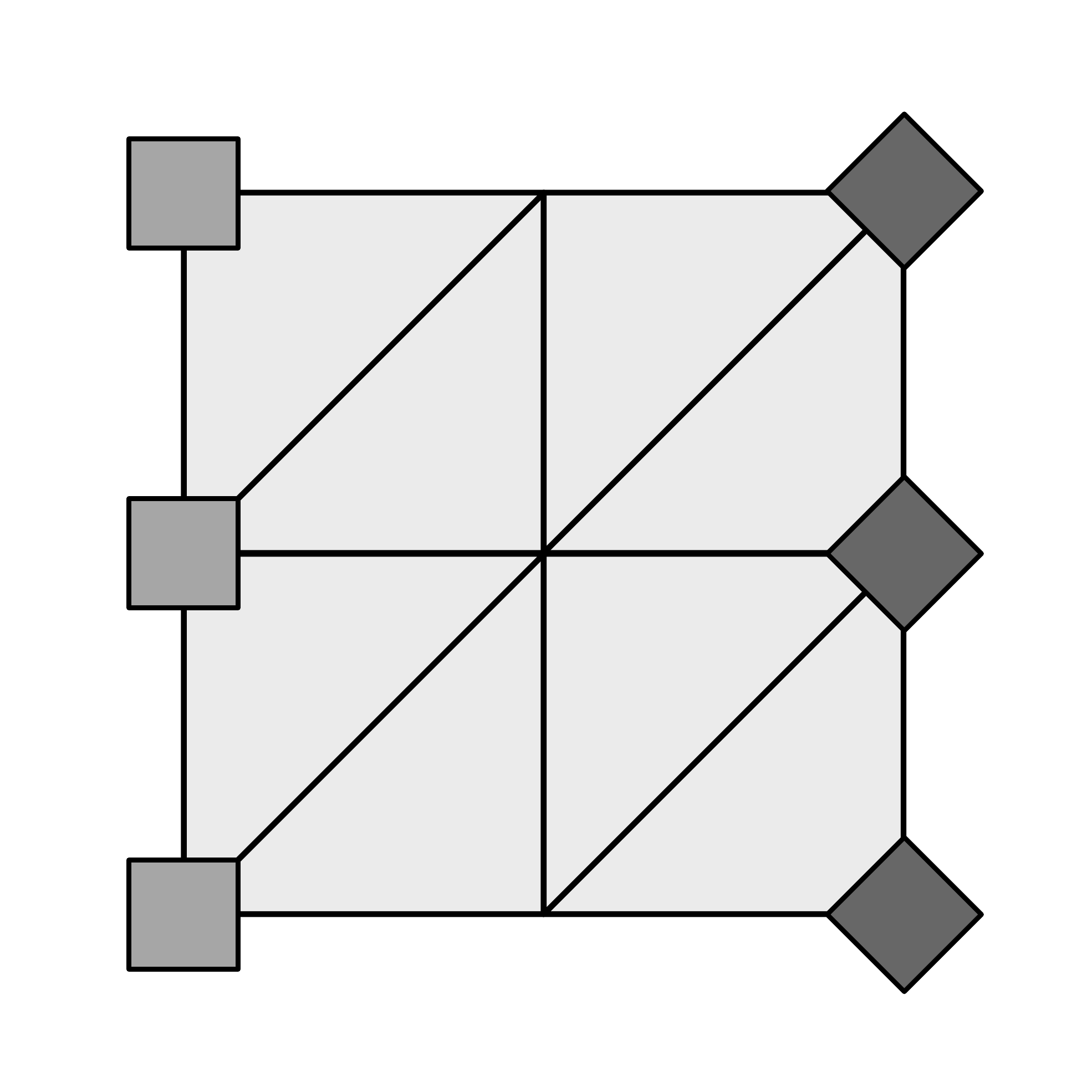}
    \includegraphics[width=.185\linewidth]{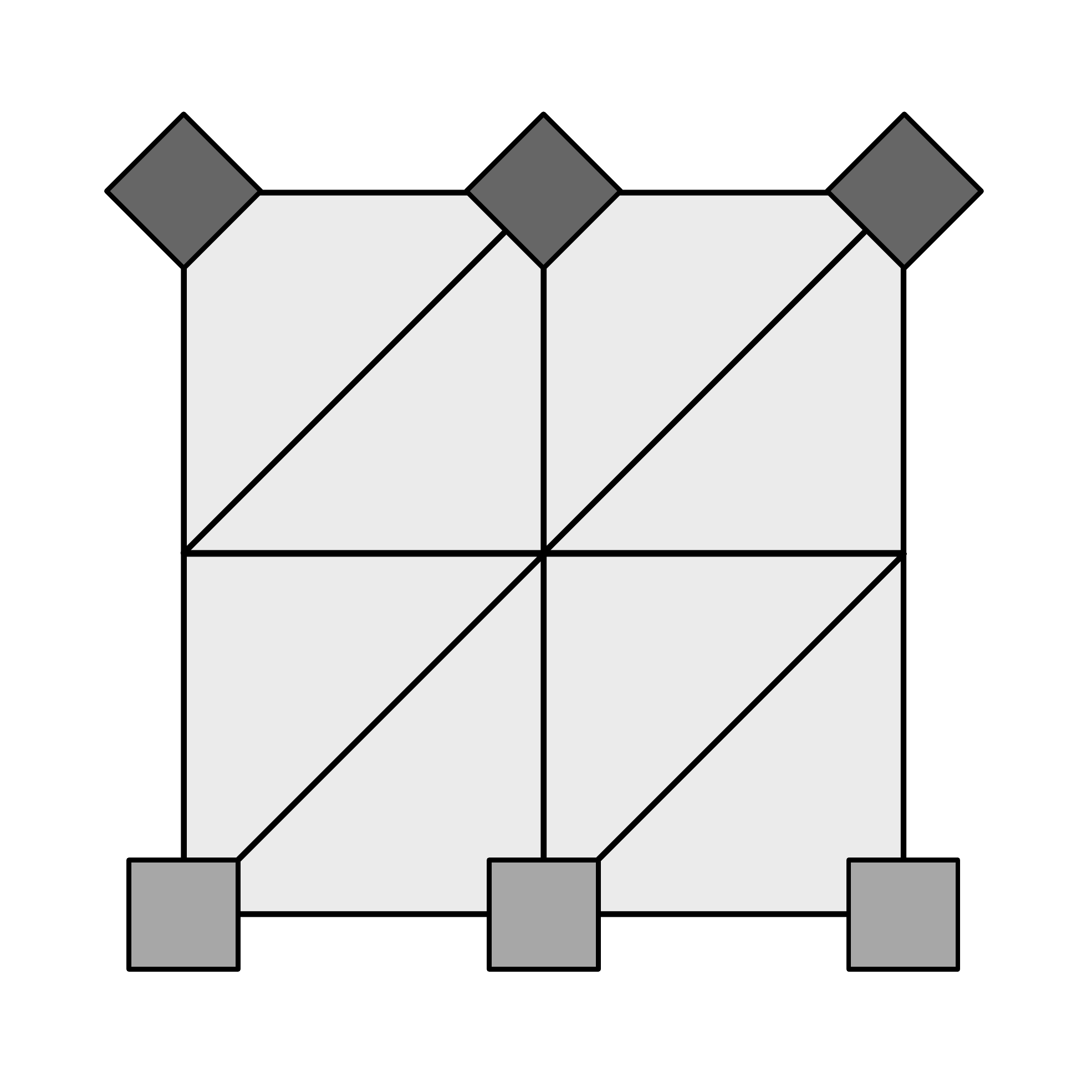}
    \includegraphics[width=.185\linewidth]{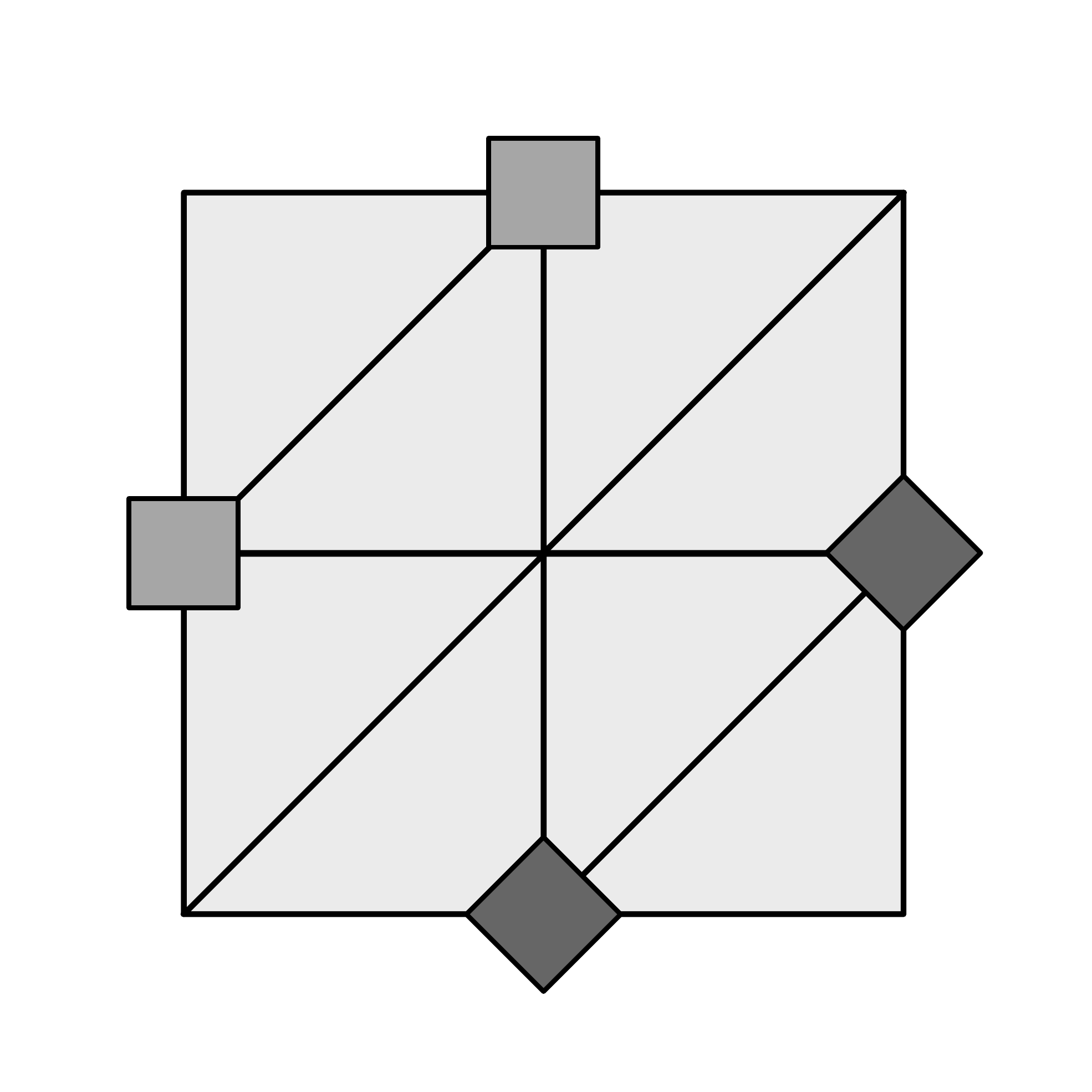}
    \includegraphics[width=.185\linewidth]{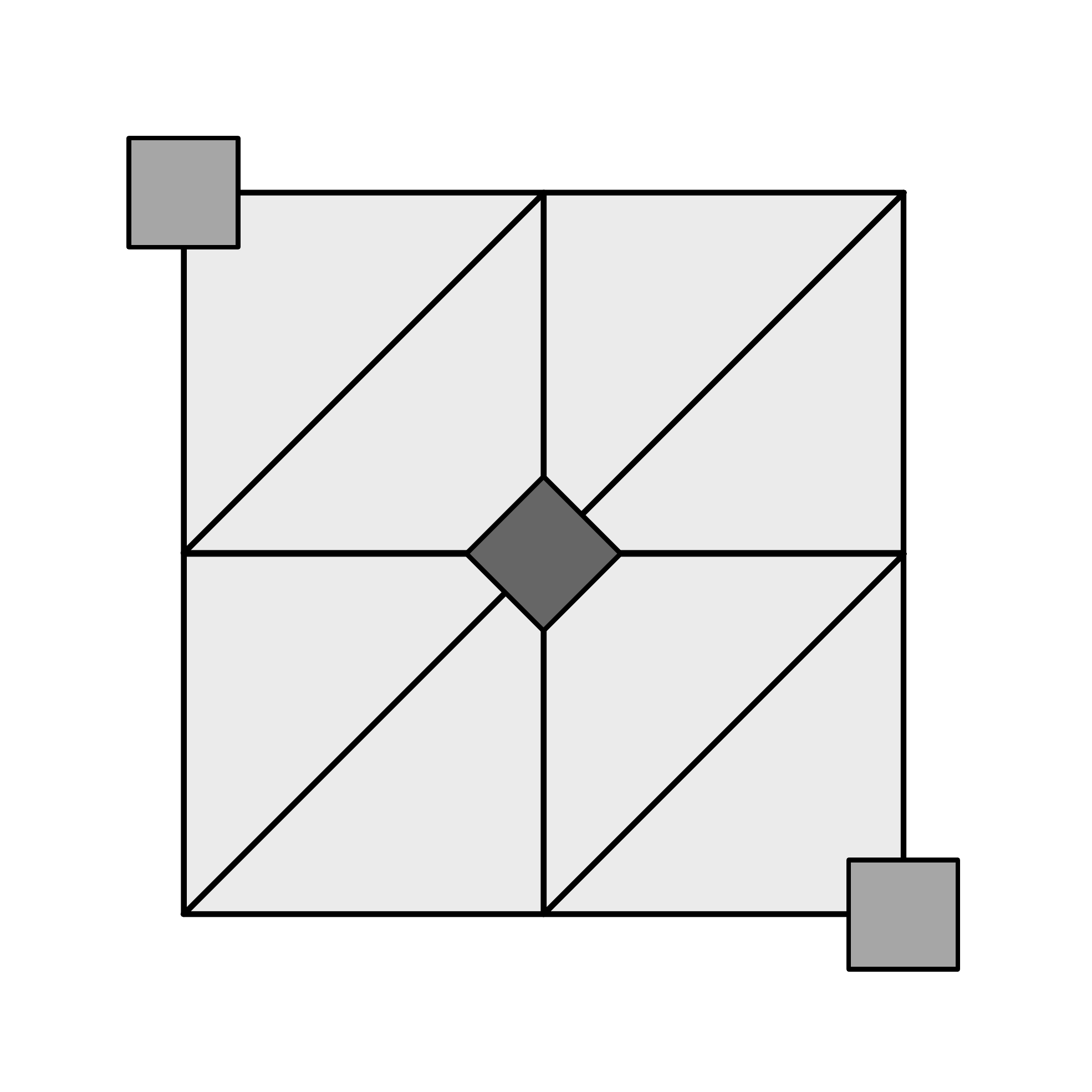} \\
    \includegraphics[width=.185\linewidth]{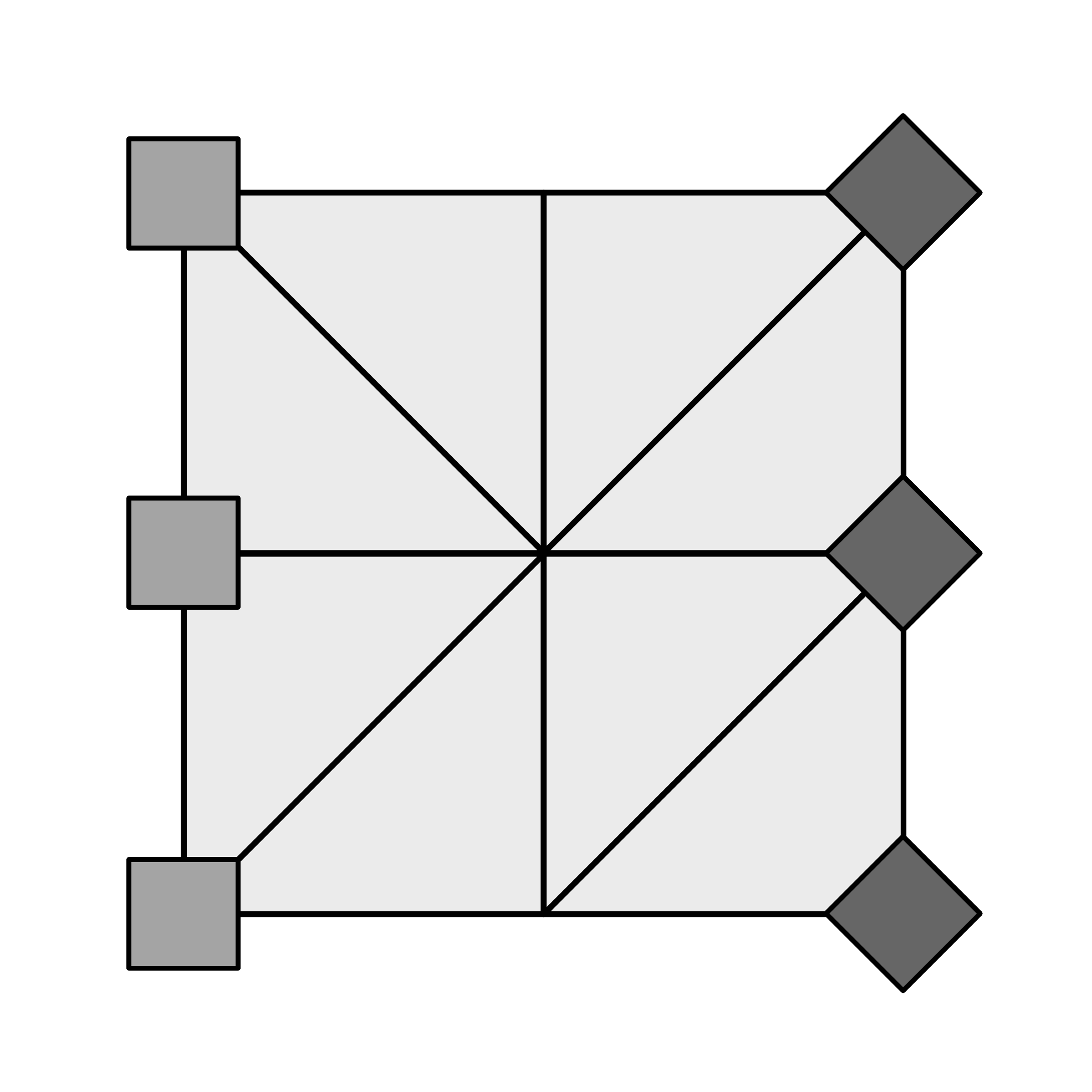}
    \includegraphics[width=.185\linewidth]{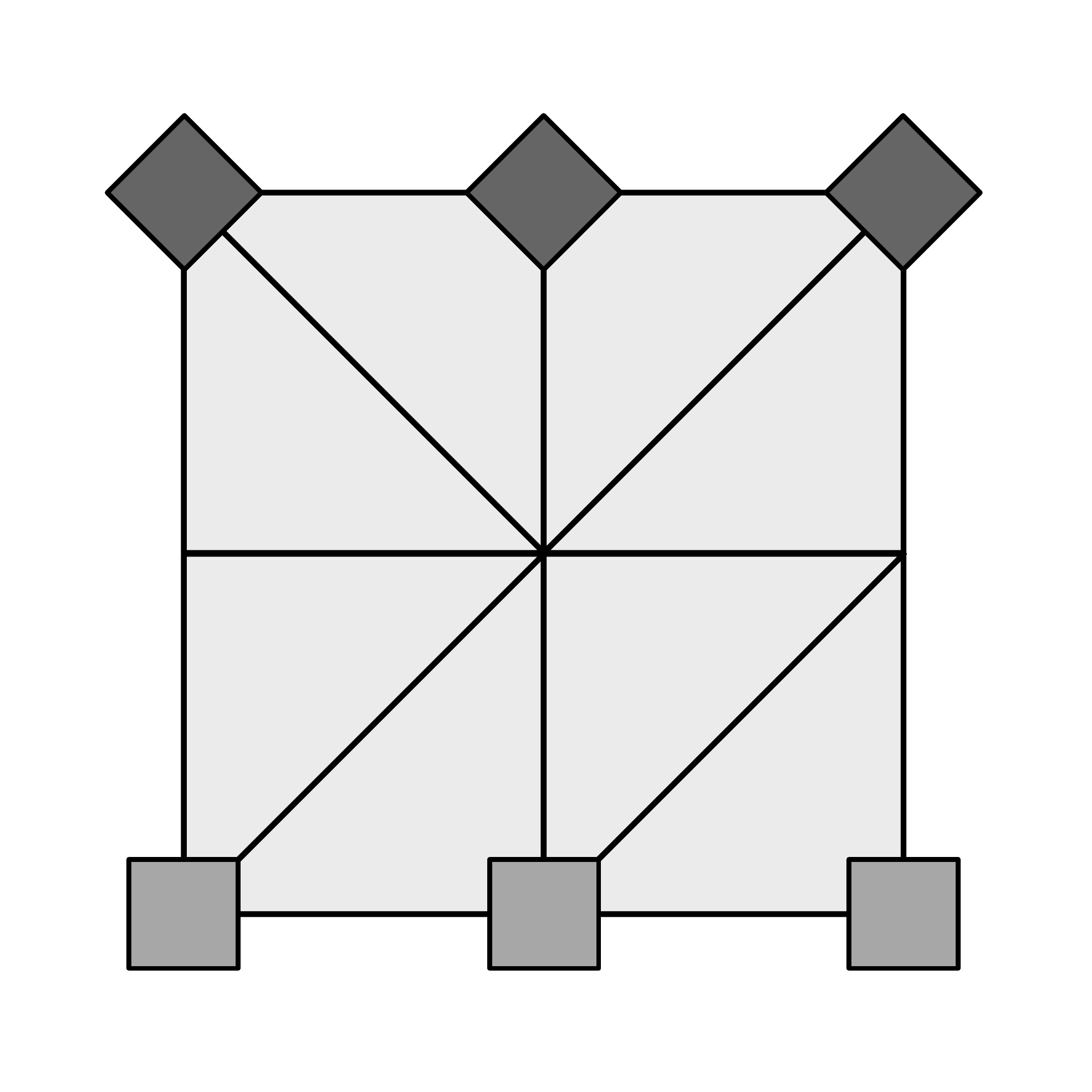}
    \includegraphics[width=.185\linewidth]{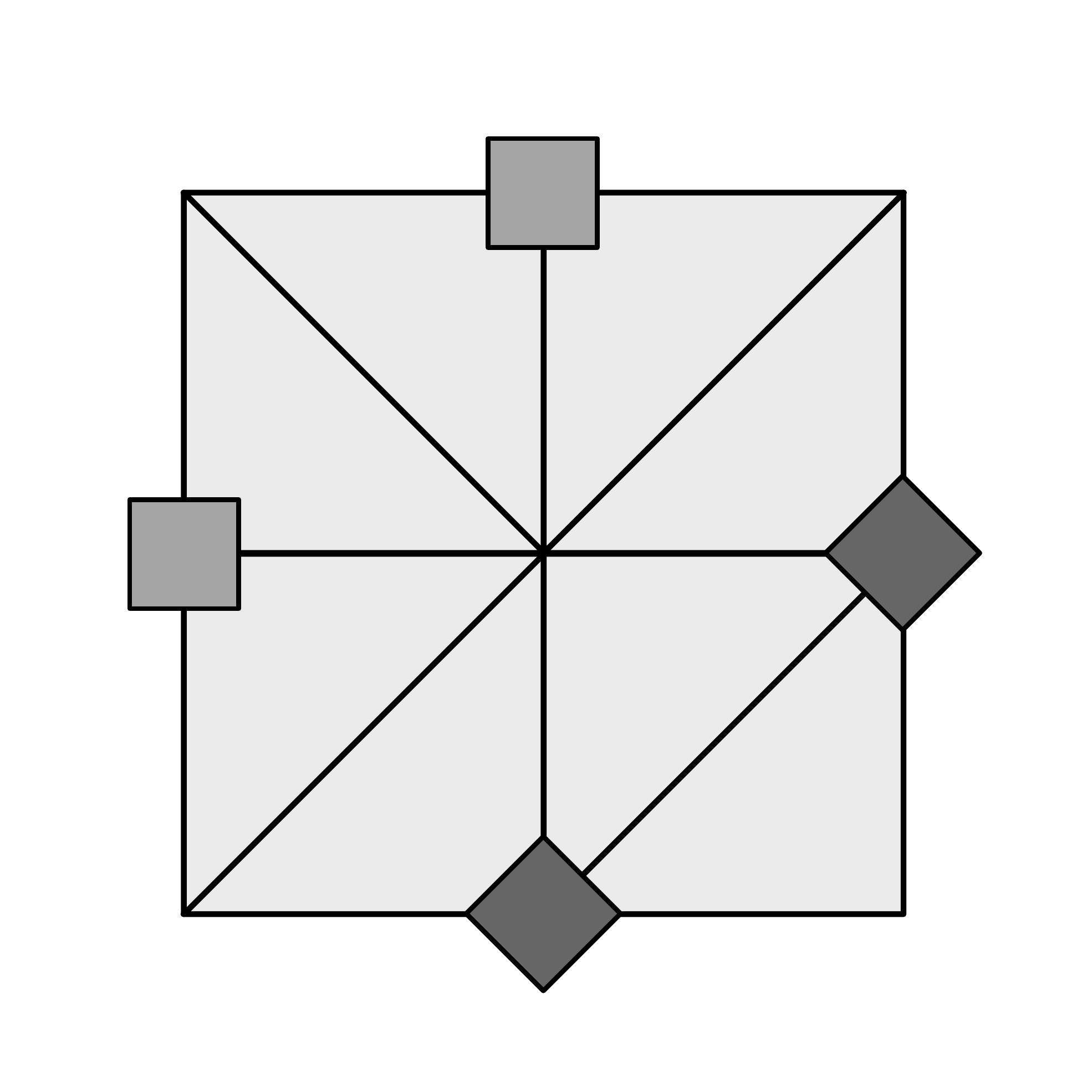}
    \includegraphics[width=.185\linewidth]{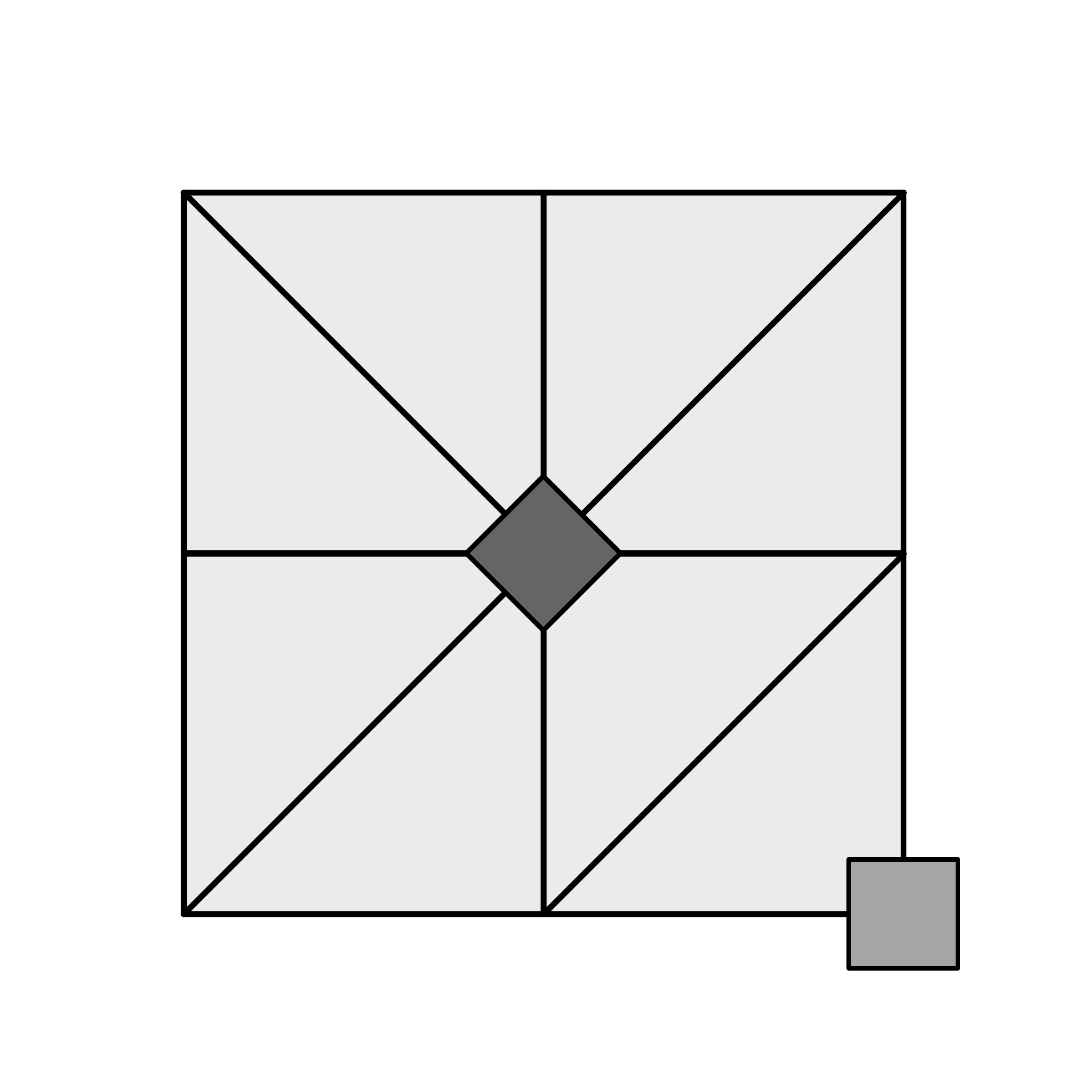}
    \includegraphics[width=.185\linewidth]{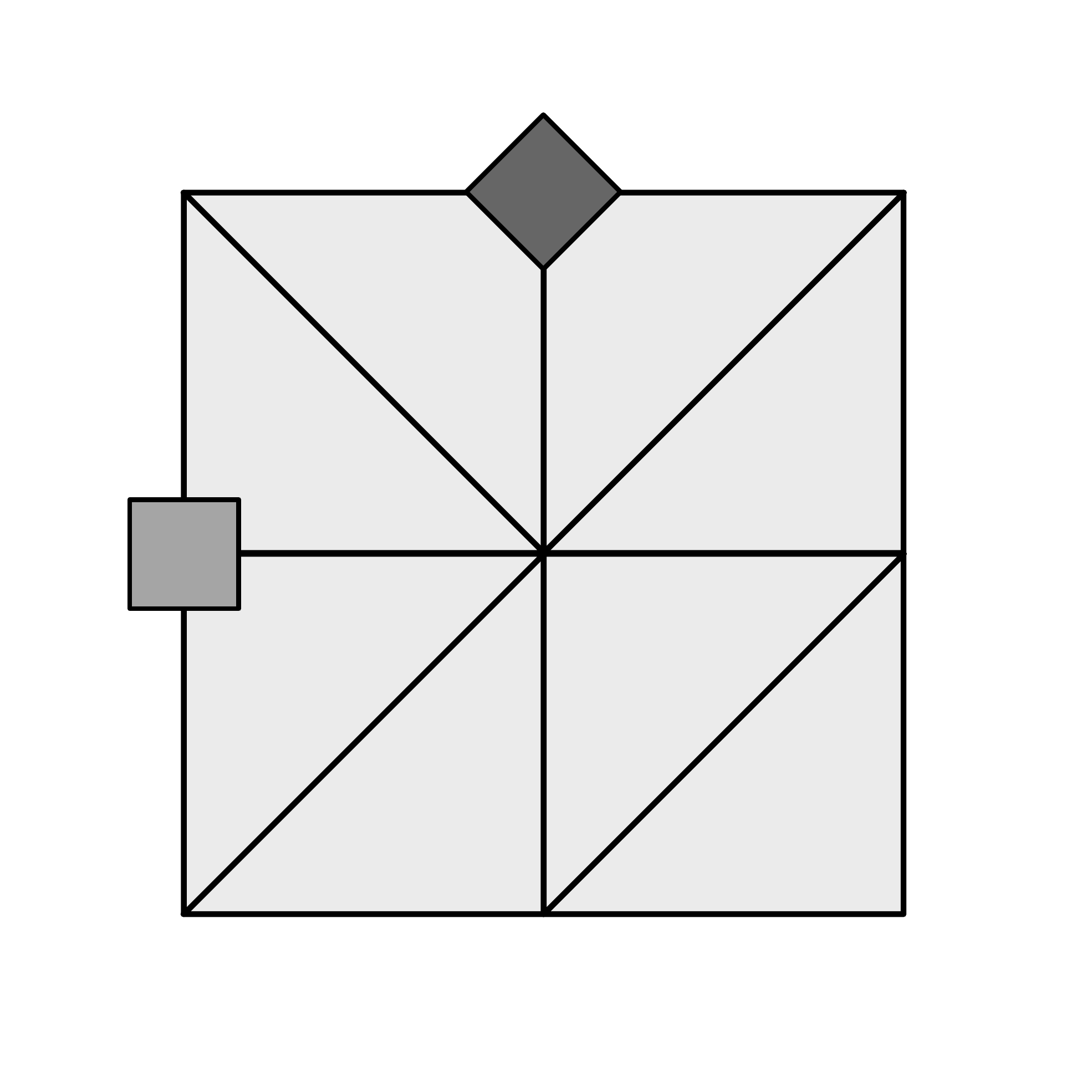}
       \caption{Independent branching schemes for the three triangulations presented in Figure \ref{fig:triangulations}, each given its own row.The sets $A^j$ and $B^j$ are given by the blue squares and green diamonds, respectively, in the $j$-th subfigure of the corresponding row.}
    \label{fig:triangulation-IBS}
\end{figure}

We can apply Lemma~\ref{lemma:graph-union} to construct biclique covers for the grid triangulations depicted in Figure~\ref{fig:triangulations}. First, we apply the biclique cover construction from Corollary~\ref{prop:multilinear-IB-scheme} to cover all edges not sharing a subrectangle. This is depicted in the first two subfigures of each row in Figure~\ref{fig:triangulation-IBS}. To cover the remaining 4 edges created by the triangulation, we see that the number of additional levels needed is dependent on the combinatorial structure. Additionally, in all three cases we can verify through Proposition~\ref{prop:feasibility-IP} that the resulting biclique cover is of the smallest possible depth.
The first example is  the ``Union Jack'' triangulation \cite{Todd:1977} for $N=M=3$ and the results in  \cite{Vielma:2009a} show that for this triangulation the biclique cover from Corollary~\ref{prop:multilinear-IB-scheme} can be completed with a single additional biclique cover for any $N$ and $M$. The second triangulation is a K1 triangulation \cite{Kuhn:1960} for $N=M=3$, and an early version of \cite{Vielma:2016} showed that for this triangulation the biclique can always be completed with two additional bicliques (See Section~\ref{chromaticsec} for a generalization of these results).

In contrast, for generic triangulations such as the third one, it was not previously known if the biclique cover can always be completed with fewer than the trivial $(M-1)(N-1)$ levels needed to cover each ``diagonal'' edge with its own additional biclique. First, we can adapt Proposition~\ref{prop:simple-pIBS} to cover the extra edges with stars, but in general this will result in $\Theta(M\cdot N)$ stars, and hence the same number of additional levels. To reduce this, we need a way to \emph{stick} the stars together into more complicated bicliques. In section~\ref{chromaticsec} we will see how we may use graph colorings for a broad class of triangulations (subsuming the Union Jack and K1 triangulations as special cases), to cover the extra edges with either one or two additional bicliques. In general, it turns out that we may cover the remaining edges for any grid triangulation with a constant number of additional levels by applying the simple following lemma.

\begin{lemma}\label{stickinglemma}
    Let $\{(A^j,B^j)\}_{j=1}^{t}$ be a family of bicliques of a graph $G$. If $(A^k,B^\ell)$ is also a biclique of $G$ for each $k,\ell\in\sidx{t}$, then  $\bra{\bigcup_{j=1}^t A^j,\bigcup_{j=1}^t B^j}$ is a biclique of $G$.
\end{lemma}

The strength of Lemma~\ref{stickinglemma} comes from the fact that many CDCs of practical interest have a local structure (i.e. sets in $\scrS$ have small cardinality, or, equivalently, the minimum degree of the conflict graph is close to the total number of nodes). In this case, the condition of Lemma~\ref{stickinglemma} will hold for families of stars centered at nodes that are located ``sufficiently far apart.''

\subsubsection{Grid triangulations of the plane} \label{ss:grid-triangulations}

\begin{figure}[htpb]
    \centering
    \includegraphics[width=.49\linewidth]{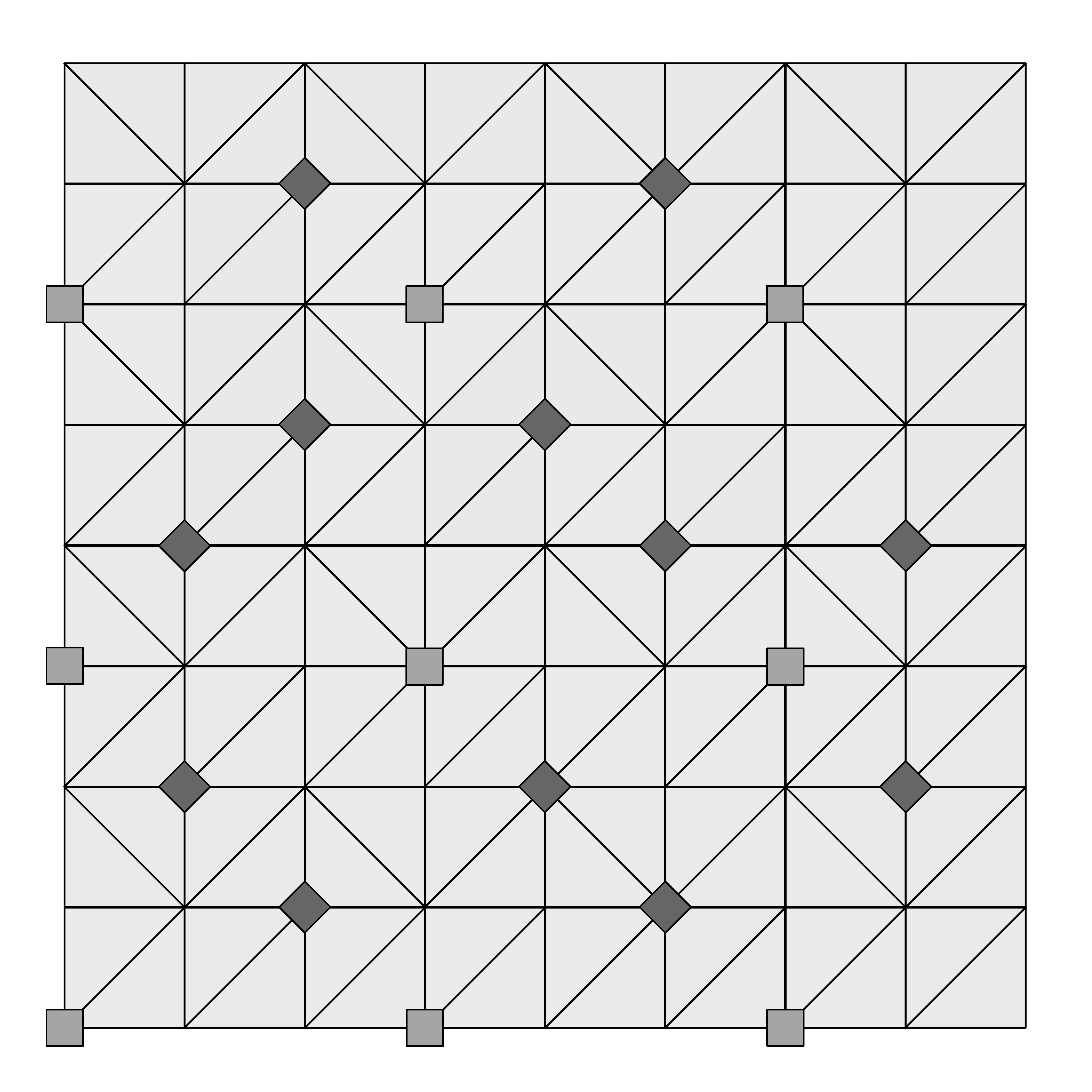}
    \includegraphics[width=.49\linewidth]{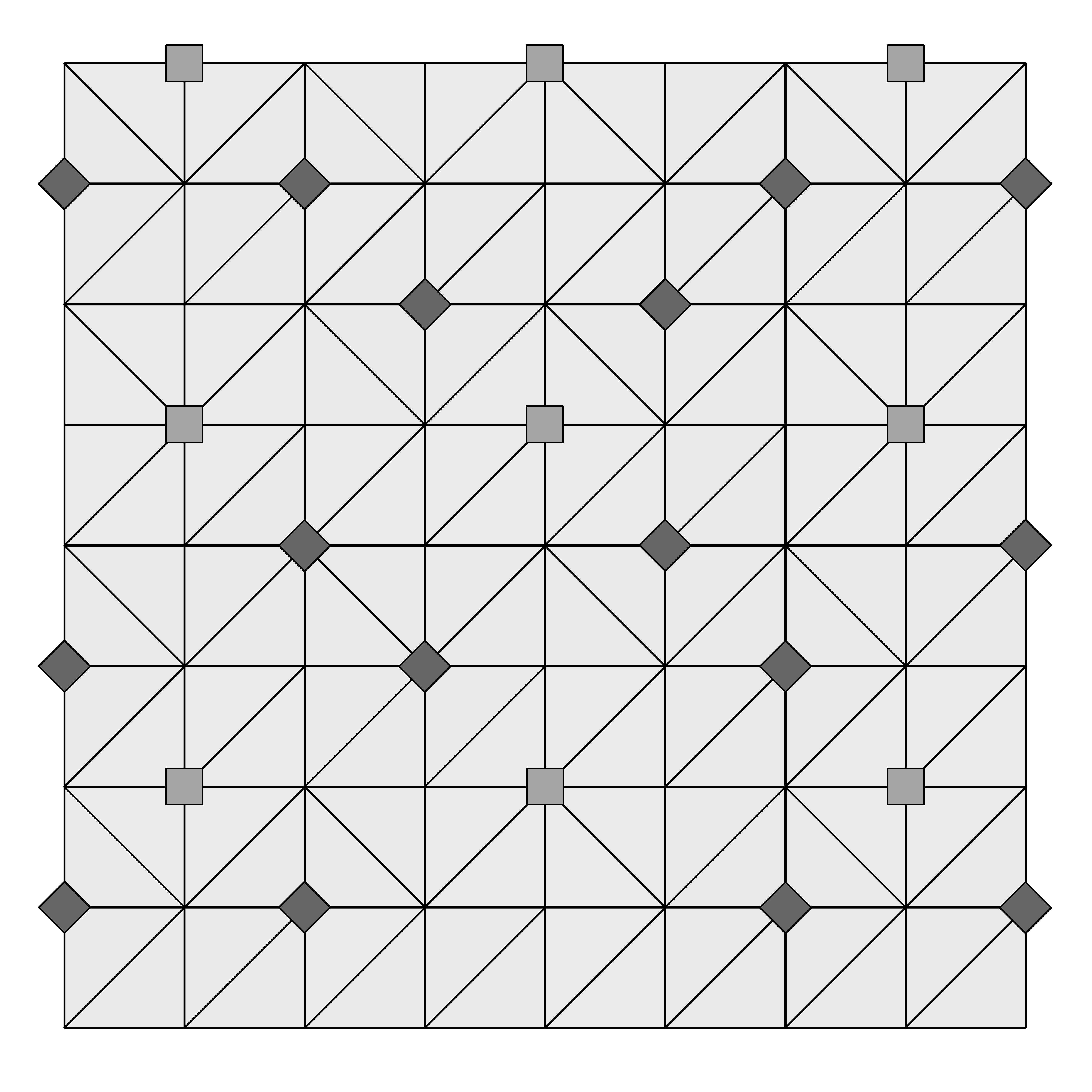}
    \caption{Two bicliques constructed via Corollary~\ref{stencil1coro} for a grid triangulation with $M=N=8$. On the left the construction follows by taking $u=(1,1)$; on the right, with $u=(2,3)$. For each level, the sets $A(u)$ and $B(u)$ are given by the squares and diamonds, respectively.}
    \label{fig:pwl}
\end{figure}

We may now present a biclique cover construction for generic grid triangulations, with no further assumptions on the structure of the triangles such as in \cite{Vielma:2010,Vielma:2009a}, whose depth scales like $\log_2(M) + \log_2(N) + \scrO(1)$. In the same way as depicted in Figure~\ref{fig:triangulation-IBS}, we construct the biclique cover by using Lemma~\ref{lemma:graph-union} to complete the construction of Corollary~\ref{prop:multilinear-IB-scheme}. For this, we will use the following corollary of Lemma~\ref{stickinglemma} that shows how to combine certain stars centered at sufficiently separated nodes.
\begin{corollary}\label{stencil1coro}
    Take a regular grid $J = \llbracket M \rrbracket \times \llbracket N \rrbracket$,  let $\scrS$ be a grid triangulation of $[1,M]\times[1,N]$, and take $G_\scrS^c=(J,\bar{E})$ as its conflict graph. For all $w\in J$, define $A(w)\defeq\set{w}$ and $B(w) \defeq \set{w+v\,:\, v\in \set{-1,1}^2, \set{w,w+v}\in\bar{E} }$. Then
    \[
        \left( \bigcup_{w \in J \cap (u+3\bbZ^2)} A(w), \bigcup_{w \in J \cap (u+3\bbZ^2)} B(w) \right)
    \]
    is a biclique of $G_\scrS^c$ for any $u\in J$.
\end{corollary}
\proof{\textbf{Proof}}
Direct from Lemma~\ref{stickinglemma} by taking $u \in J$ and the family of bicliques $\{(A(w),B(w)\}_{w \in J \cap (u + 3\bbZ^2)}$ and noting that, if $u,v\in J \cap (u+3\bbZ^2)$, then  $||u-v||_\infty \geq 3$, and so $(A(u),B(v))$ is also a biclique for $G^c_\scrS$.
\Halmos\endproof
Figure~\ref{fig:pwl} shows two possible bicliques that can be obtained from Corollary~\ref{stencil1coro}.

We can now use Lemma~\ref{lemma:graph-union} and Corollary~\ref{stencil1coro}, along with the biclique cover derived in Corollary~\ref{prop:multilinear-IB-scheme}, to obtain a biclique cover for any triangulation with an asymptotically optimal number of levels.

\begin{theorem} \label{thm:pwl}
Take  $J = \llbracket M \rrbracket \times \llbracket N \rrbracket$, and let $\scrS$ be a grid triangulation of $[1,M]\times[1,N]$. Take $G^c_\scrS=(J,\bar{E})$ as its conflict graph. Presume that $\{(\tilde{A}^{1,j},\tilde{B}^{1,j})\}_{j=1}^{t_1}$ and $\{(\tilde{A}^{2,j},\tilde{B}^{2,j})\}_{j=1}^{t_2}$ are biclique covers for the conflict graphs of the SOS2($M$) and SOS2($N$) constraints, respectively. Furthermore, define

    \begin{align*}
        A^{3,u} &= J \cap (u + 3\bbZ^2) \\
        B^{3,u} &= \bigcup_{w\in J\cap\bra{u+3\mathbb{Z}^2}}\set{w+v\,:\,v\in \set{-1,1}^2,\: \set{w,w+v}\in\bar{E} }
    \end{align*}
for each $u \in \{0,1,2\}^2$.
Then $\{(A^{1,j},B^{1,j})\}_{j=1}^{t_1} \cup \{(A^{2,j},B^{2,j})\}_{j=1}^{t_2} \cup \{(A^{3,u},B^{3,u})\}_{u \in \{0,1,2\}^2}$ is a biclique cover for $G^c_\scrS$, where
\begin{alignat*}{2}
    A^{1,j} &= \tilde{A}^{1,j} \times \llbracket N \rrbracket, \quad\quad B^{1,j} &= \tilde{B}^{1,j} \times \llbracket N \rrbracket, \\
    A^{2,j'} &= \llbracket M \rrbracket \times \tilde{A}^{2,j'}, \quad\quad B^{2,j'} &= \llbracket M \rrbracket \times \tilde{B}^{2,j'},
\end{alignat*}
for each $j \in \llbracket t_1 \rrbracket$ and $j' \in \llbracket t_2 \rrbracket$.

 In particular, if $\{h^{1,i}\}_{i=1}^{M-1} \subseteq \{0,1\}^{\lceil \log_2(M-1) \rceil}$ and $\{h^{2,i}\}_{i=1}^{N} \subseteq \{0,1\}^{\lceil\log_2(N-1)\rceil}$ are  Gray codes, where $h^{1,0} \defeq h^{1,1}$, $h^{1,M} \defeq h^{1,M-1}$, $h^{2,0} \defeq h^{2,1}$, and $h^{2,N} \defeq h^{2,N-1}$, then a biclique cover of $G_\scrS^c$  of depth $\lceil \log_2(M-1) \rceil + \lceil \log_2(N-1) \rceil + 9$ is given by:
    \begin{subequations}\label{generalPWL}
    \begin{align}
        A^{1,j} &= \big\{ (x,y) \in J : h^{1,x-1}_j = h^{1,x}_j = 1 \big\} \\
        B^{1,j} &= \big\{ (x,y) \in J : h^{1,x-1}_j = h^{1,x}_j = 0 \big\} \\
        A^{2,j'} &= \big\{ (x,y) \in J : h^{2,y-1}_{j'} = h^{2,y}_{j'} = 1 \big\} \\
        B^{2,j'} &= \big\{ (x,y) \in J : h^{2,y-1}_{j'} = h^{2,y}_{j'} = 0 \big\} \\
        A^{3,u} &= J \cap (u + 3\bbZ^2) \\
        B^{3,u} &= \bigcup_{w\in J\cap\bra{u+3\mathbb{Z}^2}}\set{w+v\,:\,v\in \set{-1,1}^2,\: \set{w,w+v}\in\bar{E} }
    \end{align}
    \end{subequations}
    for all $j \in \llbracket \lceil \log_2(M-1) \rceil \rrbracket$, $j' \in \llbracket \lceil \log_2(N-1) \rceil \rrbracket$, and $u \in \{0,1,2\}^2$.
\end{theorem}
\proof{\textbf{Proof}}
    Let $G^x \defeq (\llbracket M \rrbracket,E^x)$ and $G^y \defeq (\llbracket N \rrbracket,E^y)$ be the conflict graphs for SOS2($M$) and SOS2($N$), respectively. Furthermore, let
    \[
        G^3 \defeq \bigcup_{u \in \{0,1,2\}^2}(J,A^{3,u} * B^{3,u})=\left(J, \bigcup_{u \in \{0,1,2\}^2}(A^{3,u} * B^{3,u})\right).
    \]
    Then we see that $G_\scrS = (G^x \times G^y) \cup G^3$ by noting that all \emph{diagonal} edges of $E$ (i.e. those of the form $\set{w,w+v}\in E$ for $w\in J$ and $v\in \set{-1,1}^2$) are included in $G^3$, and observing that $G^3$ is a subgraph of $G^c_\scrS$. The result then follows from Lemma~\ref{lemma:graph-product}, Lemma~\ref{lemma:graph-union}, and Corollary~\ref{stencil1coro}.
\Halmos\endproof

By referring to Proposition \ref{prop:log-bound}, we recover a $\lceil \log_2(2(M-1)(N-1)) \rceil \geq \lceil\log_2(M-1)\rceil + \lceil\log_2(N-1)\rceil$ lower bound on the depth of any biclique cover for a grid triangulation, and see that our construction yields a MIP formulation that is within a constant additive factor of the smallest possible. Finally, we note that similarly to the results in \citep{Vielma:2010,Vielma:2009a} for the  Union Jack triangulation, formulation \eqref{eqn:ideal-formulation-for-IBS} for biclique cover \eqref{generalPWL} can provide a significant computational advantage for general grid triangulations \cite{Huchette:2017}.

\subsubsection{SOS$k$}\label{ss:sosk}
In this subsection, we will see how we may use the graph union construction to produce an IB scheme for SOS$k$($N$) of depth $\log_2(N/k)+\scrO(k)$, for any $k$ and $N$. Similar to the construction for grid triangulations, we first construct an initial family of bicliques based on the SOS2 constraint. Next, we expand this onto a larger node set by the graph product construction. Finally, we complete the biclique cover by combining a family of sufficiently separated stars. For grid triangulations, this approach meant applying SOS2 constraints horizontally and vertically, and taking a graph product of the two. One way to interpret this is as an SOS2 constraint applied to groups of aggregated nodes in the ground set (e.g. when SOS2 is applied horizontally, we group all elements with the same horizontal coordinate into a single group). For the SOS$k$($N$) constraint, we will apply the SOS2 constraint to the groups obtained by partitioning the $N$ original ground elements into $\lceil N/k \rceil$ subsets of $k$ consecutive elements. The following simple lemma shows how this grouping can also be represented through a graph product. For the remainder of the section, we assume that $N/k$ is integer; if this is not true, we artificially introduce $\lceil N/k \rceil k - N$ nodes such that this is the case, construct the formulation in Theorem~\ref{thm:sosk}, and remove the artificial nodes from the formulation afterwards.

\begin{lemma}\label{soskproduct}
    Let $J = \llbracket N \rrbracket$, $k\leq N $,  $\scrS$ correspond to the SOS$k(N)$ constraint, and  $G_\scrS^c=(J,\bar{E})$ be the corresponding conflict graph. Let $G^1 = (\llbracket  N/k \rrbracket,E^1)$ be the conflict graph for SOS2($N/k$) and $G^2 = (\{0,\ldots,k-1\},\emptyset)$ be the empty graph on $k$ nodes. Then $G^1 \times G^2$ is isomorphic to a subgraph $\hat{G}=(J,\hat{E})$ of $G^c_\scrS$ wherein $\hat{E} \subseteq \bar{E}$, and each edge $\{u,v\} \in \bar{E}$ with $|u-v| \geq 2k$ is contained $\{u,v\} \in \hat{E}$.
\end{lemma}
\proof{\textbf{Proof}}
Let $G^1 \times G^2 = (J',E')$. Consider the bijection $f:J\to J'$ given by $f(u)=(\div(u,k),\rem(u,k))$, where $\div(u,k)=\lfloor u/k \rfloor$ and $\rem(u,k)=u-k \div(u,k)$ are the quotient and remainder of the division of $u$ by $k$, so that  $f^{-1}\bra{m,r}=km+r$. We have that $\set{(m,r),(m',r')} \in E'$ if and only if $\set{m,m'} \in E^1$, which in turn is equivalent to $|m-m'| \geq 2$. Therefore, for any $\set{(m,r),(m',r')} \in E'$, we have
    \[
       \abs{f^{-1}(m,r) -f^{-1}(m',r')}= |(km+r)-(km'+r')| = |k(m-m')+(r-r')| \geq k|m-m'| + |r-r'| \geq 2k,
    \]
    and hence $\set{f^{-1}(m,r),f^{-1}(m',r')} \in \bar{E}$, i.e. $\hat{E} \subseteq \bar{E}$. For the second condition, see that if $u,v\in J$ are such that $\abs{u-v}\geq 2k$, then $|\div(u)-\div(v)| \geq 2$, and therefore $\set{f(u),f(v)} \in E'$.
\Halmos\endproof
We can then cover the remaining edges with the following bicliques obtained by stitching together families of sufficiently separated stars.
\begin{corollary}\label{stencil2coro}
    Let $J = \llbracket N \rrbracket$, $k\leq N $, $\scrS$ correspond to the SOS$k(N)$ constraint, and  $G_\scrS^c$ be the corresponding conflict graph. For all $w\in J$, define $A(w) \defeq \set{w}$ and $B(w) \defeq \set{u\in J\,:\, k\leq \abs{u-w}< 2k }$. Then
    \[
        \bra{\bigcup_{w\in J\cap\bra{u+3k\bbZ}} A(w),\bigcup_{w\in J\cap\bra{u+3k\bbZ}} B(w)}
    \]
    is a biclique of $G_\scrS^c$ for any $u\in J$.
\end{corollary}
\proof{\textbf{Proof}}
Direct from Lemma~\ref{stickinglemma} by considering the family of bicliques $\{(A(w),B(w)\}_{w \in J \cap (u+3k\bbZ)}$ and noting that, for distinct $u,v\in J \cap (u+3k\bbZ)$, $|u-v| \geq 3k$, and so $(A(u),B(v))$ is also a biclique for $G^c_\scrS$.
\Halmos\endproof

\begin{figure}[htpb]
    \centering
    \includegraphics[width=\linewidth]{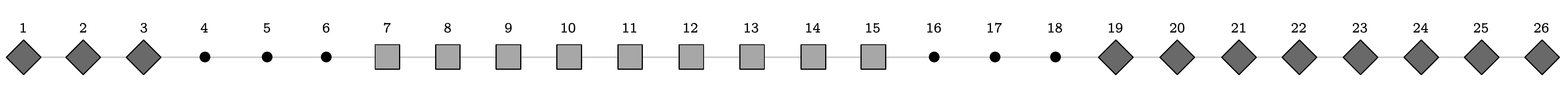}
    \includegraphics[width=\linewidth]{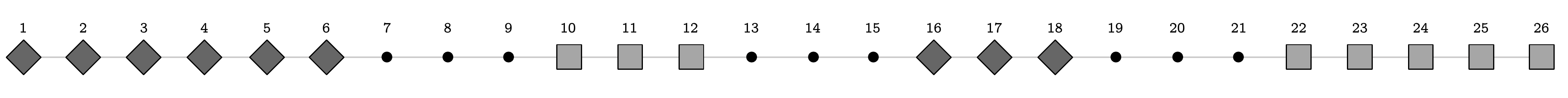}
    \includegraphics[width=\linewidth]{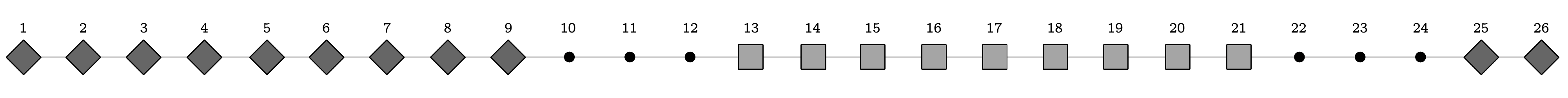}
    \includegraphics[width=\linewidth]{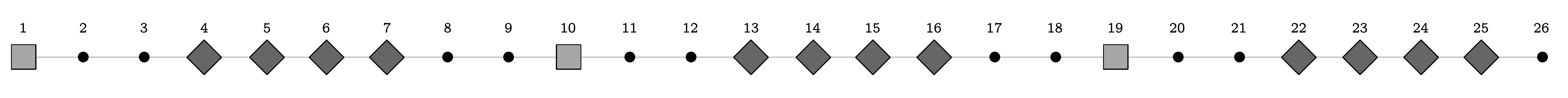}
    \includegraphics[width=\linewidth]{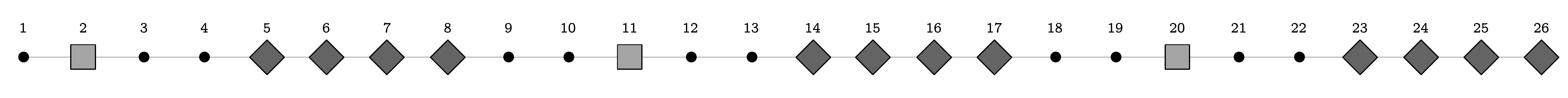}
    \includegraphics[width=\linewidth]{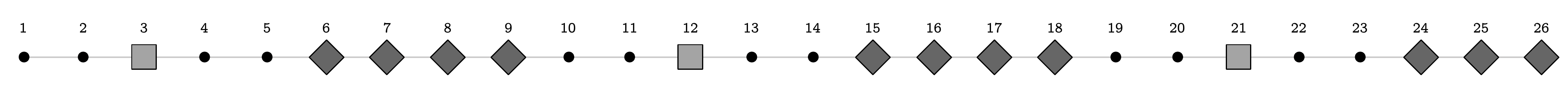}
    \includegraphics[width=\linewidth]{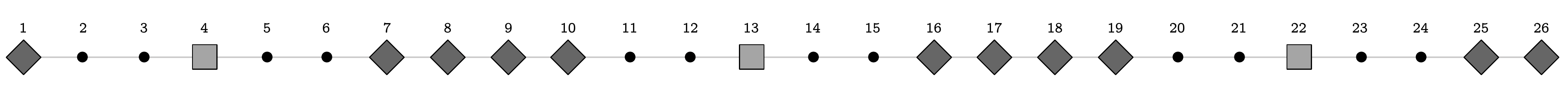}
    \includegraphics[width=\linewidth]{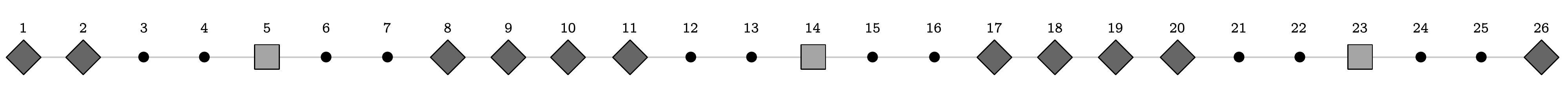}
    \includegraphics[width=\linewidth]{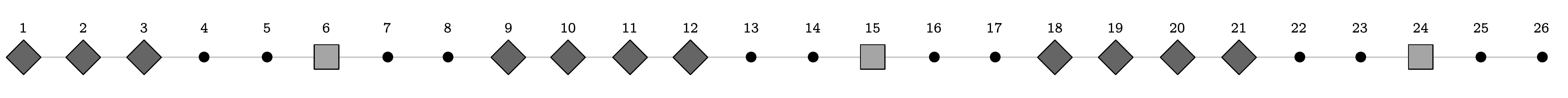}
    \includegraphics[width=\linewidth]{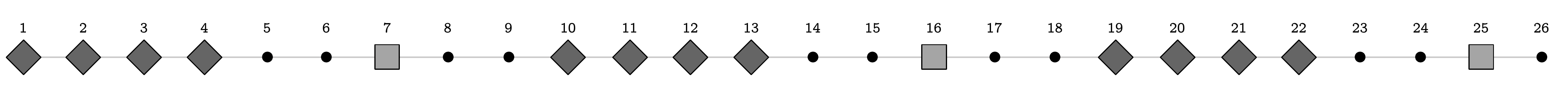}
    \includegraphics[width=\linewidth]{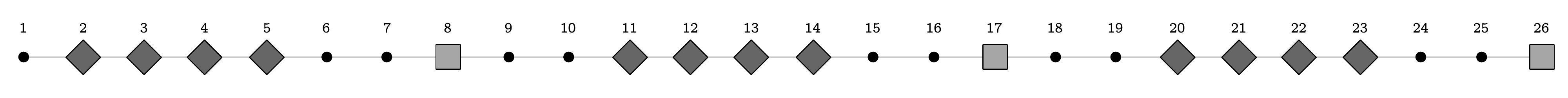}
    \includegraphics[width=\linewidth]{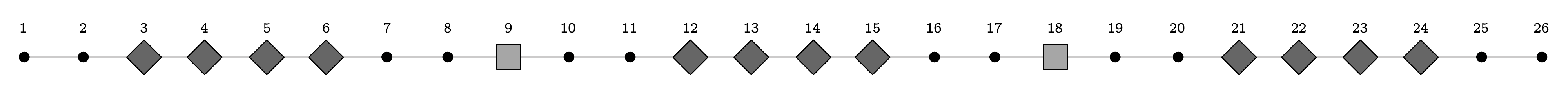}
    \caption{Visualizations of the biclique cover from the proof of Theorem~\ref{thm:sosk} for SOS$3(26)$. Each row corresponds to some level $j$, and the sets $A^j$ and $B^j$ are the squares and diamonds, respectively. The first three rows correspond the the ``first stage'' of the biclique cover $\{(A^{1,j},B^{1,j})\}_{j=1}^3$, and the second nine correspond to the ``second stage'' $\{(A^{2,j},B^{2,j})\}_{j=1}^9$.}
    \label{fig:sos3-26}
\end{figure}

Finally, we can combine both classes of bicliques with Lemma~\ref{lemma:graph-union} to construct a complete biclique cover for SOS$k(N)$.
See Figure~\ref{fig:sos3-26} for an example of the resulting construction.

\begin{theorem} \label{thm:sosk}
Let $J = \llbracket N \rrbracket$, $k\leq N $, $\scrS$ correspond to the SOS$k(N)$ constraint on $J$, and $G_\scrS^c$ be the corresponding conflict graph. Let $\{(\tilde{A}^{1,j},\tilde{B}^{1,j})\}_{j=1}^{t_1}$ be a biclique cover for the conflict graph of the SOS2($N/k$) constraint, and take
  \begin{align*}
        A^{2,j'} &= \bigcup_{i=0}^{\left\lceil \frac{N}{3k} \right\rceil} \left\{ \tau \in J : \tau = j' + (3i-3)k \right\} \\
        B^{2,j'} &= \bigcup_{i=0}^{\left\lceil \frac{N}{3k} \right\rceil} \left\{ \tau \in J : j'+(3i-2)k \leq \tau \leq j'+(3i-1)k \right\}
    \end{align*}
   for all $j' \in \llbracket 3k \rrbracket$. Then $\{(A^{1,j},B^{1,j})\}_{j=1}^{t} \cup \{(A^{2,j'},B^{2,j'})\}_{j'=1}^{\llbracket 3k \rrbracket} $ is a biclique cover for $G_\scrS^c$, where
\begin{alignat*}{2}
    A^{1,j} &= \set{ \tau \in J\,:\, \lceil \tau/k \rceil \in \tilde{A}^{1,j} }, \quad\quad B^{1,j} &=\set{ \tau \in J\,:\, \lceil  \tau/k \rceil \in\tilde{B}^{1,j}},
\end{alignat*}
for each $j \in \llbracket t \rrbracket$.

    In particular, if $\{h^i\}_{i=1}^{\lceil N/k \rceil-1} \subseteq \{0,1\}^{\lceil \log_2(\lceil N/k \rceil-1) \rceil}$ is a Gray code where $h^0 \defeq h^1$ and $h^{\lceil N/k \rceil } \defeq h^{\lceil N/k \rceil-1}$, then a  biclique cover of $G_\scrS^c$ of depth $\lceil \log_2(\lceil N/k \rceil-1)\rceil+3k$ is given by $\{(A^{1,j},B^{1,j})\}_{j=1}^{\lceil \log_2(\lceil N/k \rceil -1)\rceil} \cup \{(A^{2,j'},B^{2,j'})\}_{j'=1}^{\llbracket 3k \rrbracket} $, where
    \begin{align*}
        A^{1,j} = \left\{ \tau \in J : h^{\lceil \tau / k \rceil-1}_j = h^{\lceil \tau / k \rceil}_j = 0 \right\},\quad \quad
        B^{1,j} = \left\{ \tau \in J : h^{\lceil \tau / k \rceil-1}_j = h^{\lceil \tau / k \rceil}_j = 1 \right\}
    \end{align*}
    for all $j \in \llbracket \lceil \log_2(\lceil N/k \rceil-1) \rceil \rrbracket$.
\end{theorem}
\proof{\textbf{Proof}}
Take $G^1 \defeq (\llbracket N/k \rrbracket,E^1)$ as the conflict graph for SOS2($N/k$), $G^2 \defeq (\{0,\ldots,k-1\},\emptyset)$ as the empty graph on $k$ nodes, and $G^3 \defeq \bra{J, \bigcup_{j'=1}^{3k} A^{2,j'} * B^{2,j'}}$. Let $\hat{G}$ be the subgraph of $G_\scrS^c$ from Lemma~\ref{soskproduct}, which is isomorphic to $G^1 \times G^2$ through the bijection $g:\llbracket N/k \rrbracket\times \{0,\ldots,k-1\}\to J$ with  $g\bra{m,r}=km+r$. Then we have that  $G^c_\scrS= \hat{G}\cup G^3$, after applying Lemma~\ref{soskproduct} and using the fact that the edges of $G^3$ contain the edges of $G^c_\scrS$ not included in $\hat{G}$. The result then follows from Lemma~\ref{lemma:graph-product}, Lemma~\ref{lemma:graph-union}, Lemma~\ref{soskproduct}, and Corollary~\ref{stencil1coro}.
\Halmos\endproof

We note that, when $k = \scrO(\log(N))$, this biclique cover yields a MIP formulation that is asymptotically tight (with respect to the number of auxiliary binary variables) with our lower bound of $\lceil \log_2(N-k+1) \rceil$ from Proposition~\ref{prop:log-bound}. We can also show an absolute lower bound of depth $k$ for any biclique cover for SOS$k$. This implies that when $k = \omega(\log(N))$, although the formulation from Theorem~\ref{thm:sosk} is not tight with respect to the lower bound from Proposition~\ref{prop:log-bound}, it is asymptotically the smallest possible formulation in the pairwise IB framework.

\begin{proposition} \label{prop:sosk-lowerbound}
    Any biclique cover for the conflict graph of SOS$k(N)$ must have depth at least $\min\{k,N-k\}$.
\end{proposition}
\proof{\textbf{Proof}}
    Define $\gamma \defeq \min\{k, N-k\}$ and consider any possible biclique cover $\{(A^j,B^j)\}_{j=1}^t$. The biclique cover must separate the edges $\{(\tau,\tau+k)\}_{\tau=1}^{\gamma}$. Consider a level $j$ of the biclique cover that contains edge $\{\tau,\tau+k\}$ for some $\tau \in \llbracket \gamma \rrbracket$; w.l.o.g., $\tau \in A^j$ and $\tau+k \in B^j$. Consider the possibility that the same level $j$ separates another such edge in the set, e.g. $(\tau',\tau'+k)$ for $\tau' \in \llbracket \gamma \rrbracket$, where w.l.o.g. $\tau < \tau'$. That would imply that either $\tau' \in A^j$ or $\tau' \in B^j$. In the case that $\tau' \in A^j$, we have that $\bar{E}^j$ contains the edge $\{\tau',\tau+k\}$. However, since $|(\tau+k)-\tau'| = \tau+k - \tau' < \tau+k-\tau = k$, this implies that the biclique cover separates a feasible edge, a contradiction. In the case where $\tau' \in B^j$, we have that $\bar{E}^j$ contains the edge $\{\tau,\tau'\}$, and as $\tau'-\tau < k \leq \gamma$ from the definition of our set of edges, a similar argument holds. Therefore, each edge $\{\{\tau,\tau+k\}\}_{\tau=1}^{\gamma}$ must be uniquely contained in some level of the biclique cover, giving the result.
\Halmos\endproof

Furthermore, when $k = \lfloor N/2 \rfloor$, this proposition gives a lower bound on the depth of a biclique cover that is asymptotically tight with the upper bound of $N$ from Proposition~\ref{prop:simple-pIBS}. In other words, in this particular regime, we have that the SOS$k$ constraint admits a pairwise IB-based formulation, but only one that is relatively large ($\Omega(|\scrS|) = \Omega(|J|)$ auxiliary binary variables and constraints).

\subsection{Covering edges with a chromatic characterization of bicliques}\label{chromaticsec}
Our final contribution is to adapt a result of Cornaz and Fonlupt~\cite{Cornaz:2006} that gives a chromatic characterization of the set of edges that can be covered by a biclique.
\begin{theorem}[\cite{Cornaz:2006}]\label{cornaztheo}
Take the graph $G=(J,\bar{E})$, along with some edge subset $\bar{F} \subseteq \bar{E}$. Define $V(\bar{F}) = \bigcup\{\{u,v\} \in \bar{F}\}$ as all nodes incident to $\bar{F}$, and take $F = [(\bar{F} * \bar{F}) \backslash \bar{E}$ as all pairs incident to the edges $V(\bar{F})$ not contained in $\bar{E}$. Define both $E' = F \cup \bar{F}$ and $p:E'\to \set{0,1}$ such that $p(e)= \bbone[e \in \bar{E}]$\footnote{Where $\bbone[e \in \bar{E}]=1$ if $e \in \bar{E}$ and $\bbone[e \in \bar{E}]=0$ otherwise.}. Finally, take $\mathcal{C}\bra{E'}$ as the family of all cycles in $G'$. Then the following are equivalent:
\begin{enumerate}
    \item There exists a biclique $(A,B)$ of $G$ covering $(V(\bar{F}),\bar{F})$.
    \item For all $C\in \mathcal{C}\bra{E'}$, $\sum\nolimits_{u\in C} p(u)$ is even.
    \item There exists some $f:V(\bar{F})\to \set{0,1}$ such that
        \begin{itemize}
            \item $f(u)=f(v)$ for all $\set{u,v}\in F$,
            \item $f(u)\neq f(v)$ for all $\set{u,v}\in \bar{F}$, and
            \item $\bra{\set{u\in V(\bar{F})\,:\, f(u)=0},\set{u\in V(\bar{F})\,:\, f(u)=1}}$ is a biclique of $G$ covering $(V(\bar{F}),\bar{F})$.
        \end{itemize}
\end{enumerate}

\end{theorem}

Using this result, we can sometimes exploit the structure of a grid triangulations to obtain biclique covers that are smaller than those derived in Theorem~\ref{thm:pwl}.

\begin{theorem} \label{thm:pwl2}
Take a regular grid $J = \llbracket M \rrbracket \times \llbracket N \rrbracket$, the sets $\scrS$ given by a grid triangulation of $[1,M]\times[1,N]$ with $\bigcup\{S \in \scrS\} = J$, and take $G_\scrS^c=(J,\bar{E})$ to be the corresponding conflict graph. Furthermore, let $\{(\tilde{A}^{1,j},\tilde{B}^{1,j})\}_{j=1}^{t_1}$ and $\{(\tilde{A}^{2,j},\tilde{B}^{2,j})\}_{j=1}^{t_2}$ be biclique covers for the conflict graphs of the SOS2($M$) and SOS2($N$) constraints, respectively. Also define $J_{even} \defeq \set{u\in J\,:\, u_1 \equiv u_2 \mod 2}$ and $J_{odd}\defeq\set{u\in J\,:\, u_1 \not\equiv u_2 \mod 2}$ as the set of nodes whose sum of components is either even and odd, respectively. For each $s \in \{even,odd\}$, let $E_s\defeq \set{\set{u,v} \in [J_s]^2 : \left\|u-v\right\|_\infty=1}$, $\bar{F}_s\defeq E_s\cap\bar{E}$, and $F_s\defeq ([V\bra{\bar{F}_s}]^2 \cap E_s) \backslash \bar{E}$.

If, for each $s \in \{even,odd\}$, there exists $f_s:J_s\to \set{0,1}$ such that
\begin{equation}\label{pwlcolorcond}
f_s(u)=f_s(v) \quad \forall \set{u,v}\in F_s\quad \quad \text{ and } \quad\quad f_s(u)\neq f_s(v)\quad \forall\set{u,v}\in \bar{F}_s,
\end{equation}
then $\{(A^{1,j},B^{1,j})\}_{j=1}^{t_1} \cup \{(A^{2,j},B^{2,j})\}_{j=1}^{t_2} \cup \{(A^{3,i},B^{3,i})\}_{i\in\sidx{2}}$ is a biclique cover for $G_\scrS^c$, where
\begin{alignat*}{2}
    A^{1,j} &= \tilde{A}^{1,j} \times \llbracket N \rrbracket, \quad\quad& B^{1,j} &= \tilde{B}^{1,j} \times \llbracket N \rrbracket, \\
    A^{2,j'} &= \llbracket M \rrbracket \times \tilde{A}^{2,j'}, \quad\quad& B^{2,j'} &= \llbracket M \rrbracket \times \tilde{B}^{2,j'},\\
    A^{3,s} &= \set{u\in J_s\,:\, f_s(u)=0} , \quad\quad& B^{3,s} &=\set{u\in J_s\,:\, f_s(u)=1}
\end{alignat*}
for each $j \in \llbracket t_1 \rrbracket$, $j' \in \llbracket t_2 \rrbracket$, and $s \in \{even,odd\}$.

If we select $\{(\tilde{A}^{1,j},\tilde{B}^{1,j})\}_{j=1}^{t_1}$ and $\{(\tilde{A}^{2,j},\tilde{B}^{2,j})\}_{j=1}^{t_1}$ to correspond to the Gray code construction for SOS2 in Proposition~\ref{prop:sos2-IB-scheme}, then the resulting biclique cover has depth $\lceil \log_2(M-1) \rceil + \lceil \log_2(N-1) \rceil + 2$. Furthermore, if for some $s \in \{even,odd\}$ we have that $f_s(u)=\alpha$ for all $u\in J_s$ and some constant $\alpha$, we may reduce the depth to $\lceil \log_2(M-1) \rceil + \lceil \log_2(N-1) \rceil + 1$.

Finally, if we fix $s \in \{even,odd\}$ and $r \in \{even,odd\} \backslash \{s\}$, a sufficent condition for the existence of $f_s:J_s\to \set{0,1}$ satisfying \eqref{pwlcolorcond} is that
\begin{equation}\label{odddegreecond}
    d_{F_{r}}(u) \defeq \left|\left\{ e \in F_{r} : \exists v\in J_{r}\text{ s.t. } e=\set{u,v} \right\}\right| \text{ is even } \quad \forall u\in J_{r}\cap \llbracket 2,M-1 \rrbracket \times \llbracket 2,N-1 \rrbracket.
\end{equation}
\end{theorem}
\proof{\textbf{Proof}}
    Let $G^x \defeq (\llbracket M \rrbracket,E^x)$ and $G^y \defeq (\llbracket N \rrbracket,E^y)$ be the conflict graphs for the SOS2($M$) and SOS2($N$) constraints, respectively. Furthermore, let $G^{3} \defeq (\llbracket M \rrbracket \times \llbracket N \rrbracket,(A^{3,even} * B^{3,even}) \cup (A^{3,odd} * B^{3,odd}))$. Then we may see that $G^c_\scrS = (G^x \times G^y) \cup G^3$ by noting that all \emph{diagonal} edges of $\bar{E}$ (i.e. those of the form $\set{w,w+v}\in \bar{E}$ for $w\in J$ and $v\in \set{-1,1}^2$) are included in $G^3$, and that $G^3$ is a subgraph of $G^c_\scrS$. The first part of the theorem then follows from Lemma~\ref{lemma:graph-product}, Lemma~\ref{lemma:graph-union}, and Theorem~\ref{cornaztheo}.

   For the sufficient condition, w.l.o.g. consider the case where $s=even$ and $r=odd$. Define $p_{even}:E_{even}\to \set{0,1}$ as $p_{even}(e) = \bbone[e \in \bar{E}]$. The result will follow from Theorem~\ref{cornaztheo} by showing that \eqref{odddegreecond} satisfies condition 2 in the equivalence of the theorem. Let $E'_{even}=F_{even}\cup\bar{F}_{even}$, and assume for contradiction that there exists $C\in \mathcal{C}\bra{E'_{even}}$ such that $\sum\nolimits_{u\in C} p_{even}(u)$ is odd. If $\abs{C}=4$, we may assume without loss of generality that $V(C)=\set{u,u+(1,1),u+(1,-1),u+(2,0)}\subset J_{even}$ for some $u\in J_{even}$. Then $v=u+(0,1)\in J_{odd}$ is such that $d_{E_{odd}}(v)$ is odd, a contradiction of \eqref{odddegreecond}. If $\abs{C}>4$, note that $C\in \mathcal{C}\bra{E_{even}}$ and that there exists $e\in C$ such that $p_{even}(e)=1$. In addition, there exists $C_0\in \mathcal{C}\bra{E_{even}}$ such that $e\in C_0$, $\abs{C_0}=4$, $C_1=\bra{C_0\cup C}\setminus \bra{C_0\cap C}\in \mathcal{C}\bra{E_{even}}$ and $\Conv(C_1) \subsetneq \Conv(C)$.
   If $\sum\nolimits_{u\in C_0} p_{even}(u)$ is odd, we may make the same argument above as $|C_0|=4$ to derive a contradiction of \eqref{odddegreecond}. If not, then $C_1$ and $C$ have the same parity, and therefore $\sum\nolimits_{u\in C_1} p_{even}(u)$ is odd. We may then repeat this shrinking procedure recursively on $C_1$ until either the $4$-cycle $C_0$ has odd parity, or $C_1$ is itself a $4$-cycle. In either case, we have a $4$-cycle with odd parity, which implies some $u\in J_{odd}$ that violates \eqref{odddegreecond}, giving the result.
\Halmos\endproof

We note that we may use this coloring characterization to recover the biclique covers for  both the Union Jack and K1 triangulation example in Figure~\ref{fig:triangulation-IBS}. For the Union Jack example, we have that
\[
    \bar{F}_{even} = \emptyset, \quad\quad \bar{F}_{odd} = \left\{\{(1,2),(2,1)\},\:\{(2,1),(3,2)\},\:\{(1,2),(2,3)\},\:\{(2,3),(3,2)\}\right\},
\]
and so we may apply the simplification in Theorem~\ref{thm:pwl2} to construct a biclique cover of depth $\log_2(2)+\log_2(2)+1=3$. Indeed, the original formulation of Vielma and Nemhauser~\cite{Vielma:2009a} for the Union Jack triangulation can be reinterpreted analogously through the chromatic characterization of Theorem~\ref{thm:pwl2}.

For the K1 triangulation example, we have
\[
    \bar{F}_{even} = \left\{ \{(1,3),(2,2)\},\:\{(2,2),(3,1)\} \right\}, \quad\quad \bar{F}_{odd} = \left\{\{(1,2),(2,1)\},\:\{(2,3),(3,2)\} \right\},
\]
giving a biclique cover of depth $\log_2(2)+\log_2(2)+2=4$.

\begin{figure}[htpb]
    \centering
    \includegraphics[width=.32\linewidth]{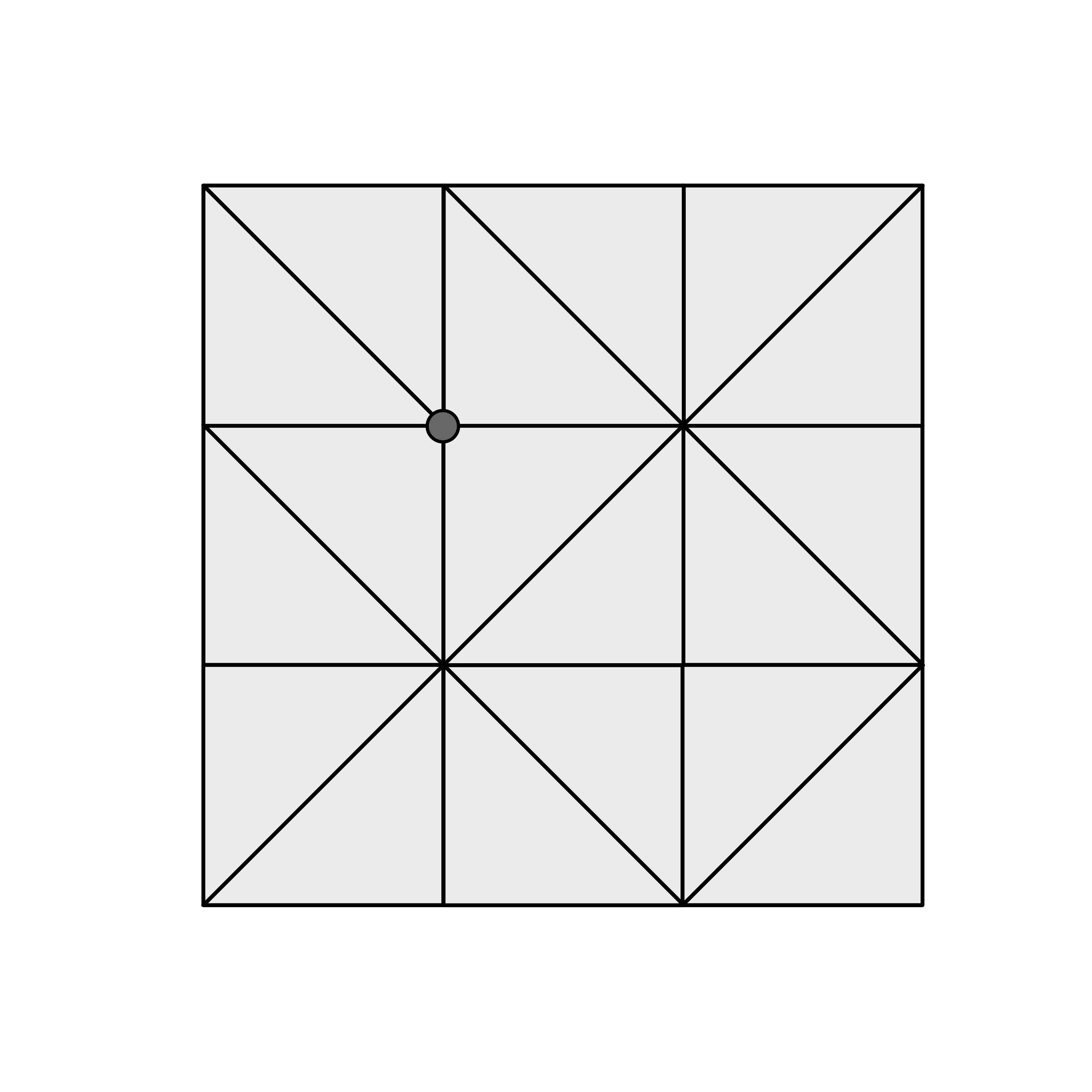}
    \includegraphics[width=.32\linewidth]{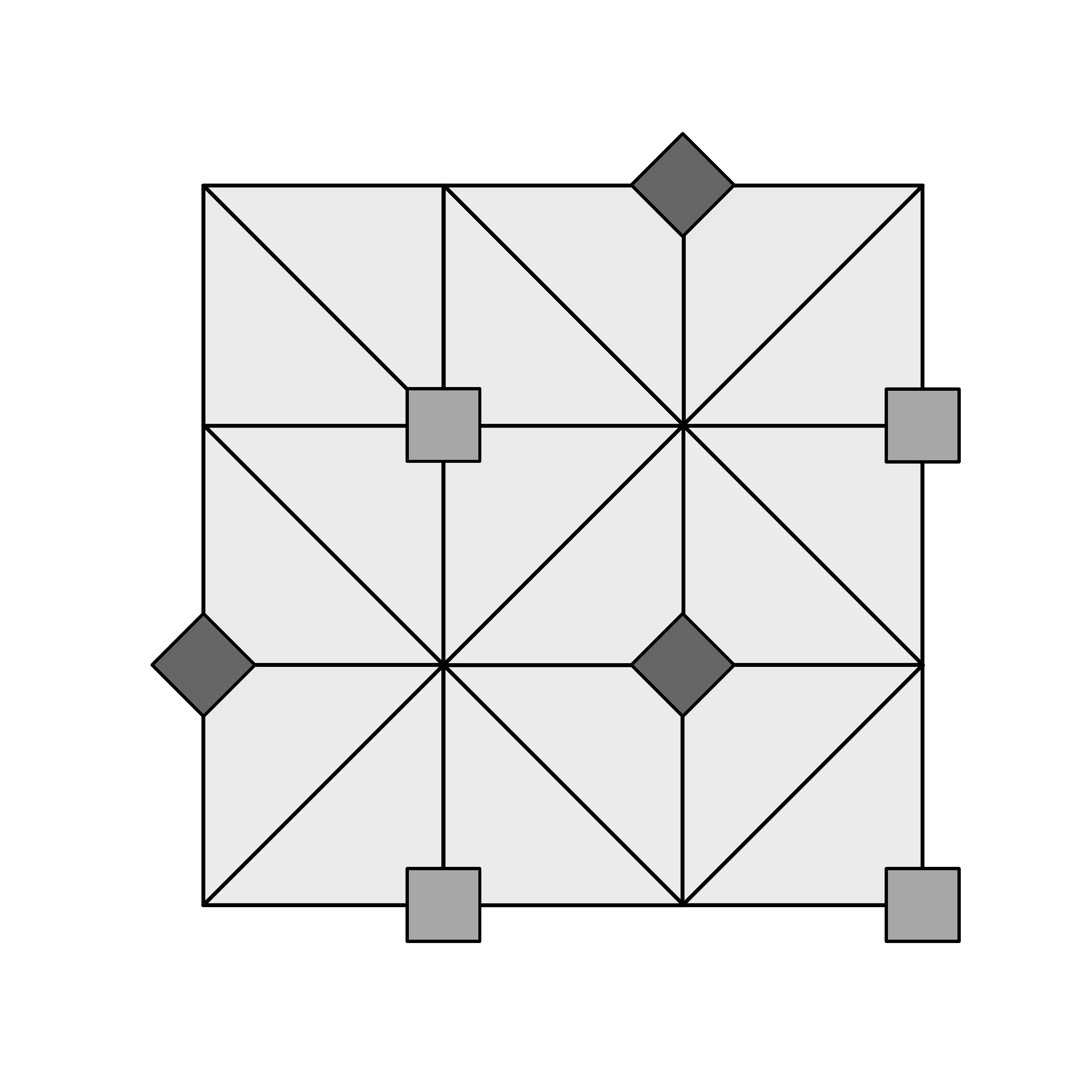}
    \includegraphics[width=.32\linewidth]{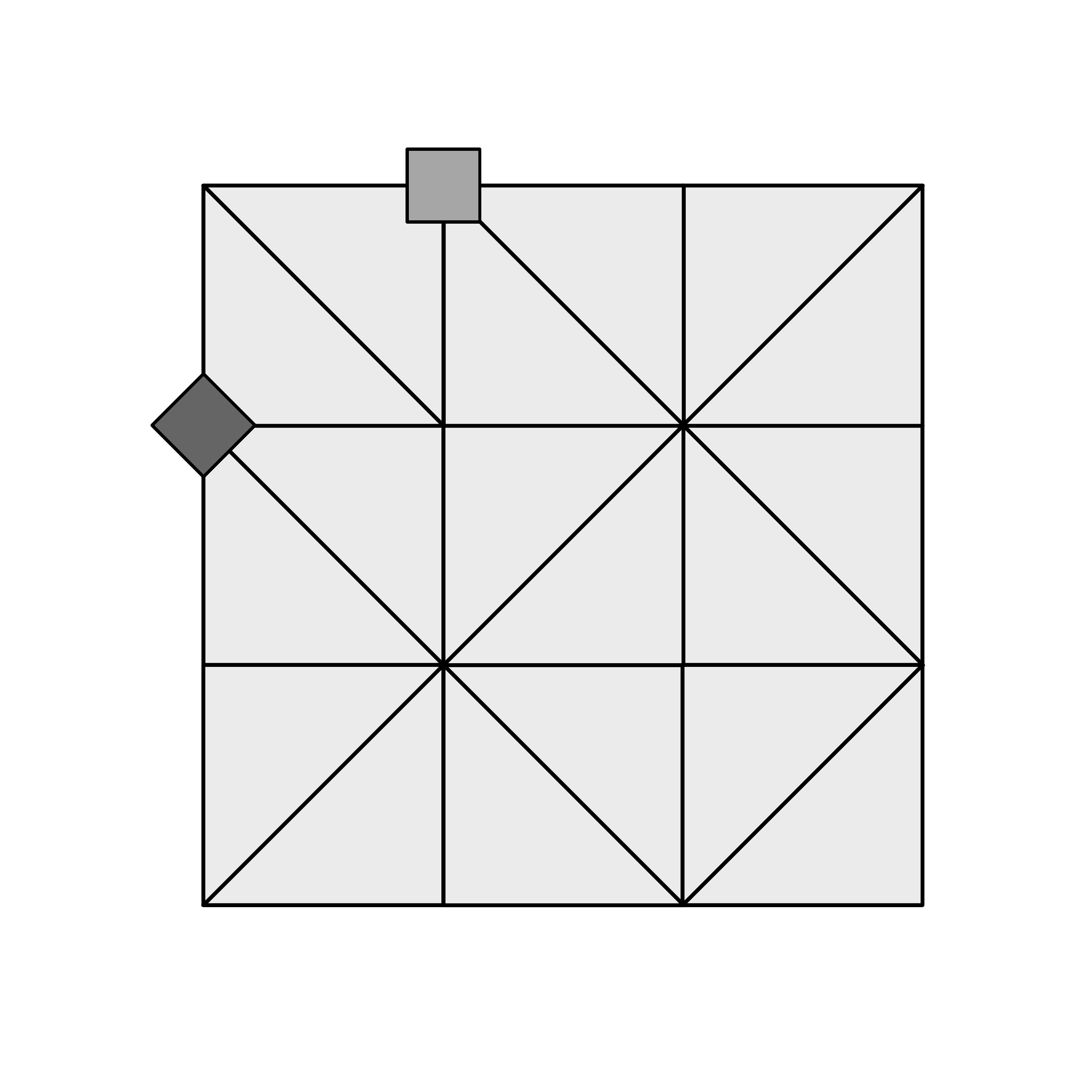}
    \caption{(Left) A grid triangulation that does not satisfy the sufficient condition of Theorem~\ref{thm:pwl2}, as $(2,3)$ (circle) has odd degree. However, there a coloring of the form described in Theorem~\ref{thm:pwl2}, leading to two bicliques (Center and Right) that cover all ``diagonal'' edges of the conflict graph.}
    \label{fig:coloring-example}
\end{figure}

Furthermore, we close by noting that the sufficient condition in Theorem~\ref{thm:pwl2} is, in general, not necessary. For example, in Figure~\ref{fig:coloring-example} we see a grid triangulation that does not satisfy \eqref{odddegreecond}, but for  which there exists a coloring given by $\{f_{even},f_{odd}\}$ that satisfies \eqref{pwlcolorcond}. That is, $(2,3)$ has odd degree (i.e. $d_{F_{odd}}((2,3)) = 1$), and so the sufficient condition is not satisfied. However, a coloring satisfying the conditions of Theorem~\ref{thm:pwl2} exists. This offers a generalization of the result that originally appeared in a preliminary version of \cite{Vielma:2016}, which only showed the sufficient condition \eqref{odddegreecond}.

\section*{Acknowledgements}
This material is based upon work supported by the National Science Foundation Graduate Research Fellowship under Grant No. 1122374, and by the National Science Foundation under Grant No. CMMI-1351619. The authors would like to thank Yves Crama for pointing out the relation between biclique covers, independent sets, and boolean functions.

\bibliographystyle{ormsv080}
\bibliography{master.bib}


\begin{APPENDICES}
\section{Logarithmic formulation from Misener et al. \cite{Misener:2011} is not ideal}\label{app:misener}
We show that the logarithmic formulation (16) from Misener et al.~\cite{Misener:2011} is not, in general, ideal. Using their notation, we take $N_P = 3$, $x^L = y^L = 0$, and $x^U = y^U = 3$ (and so $a=1$). Then formulation (16) is
\begin{subequations}
\begin{align}
    \lambda_1 + 2\lambda_2 &\leq x \\
    x &\leq 1 + \lambda_1 + 2\lambda_2 \label{freq1}\\
    1 + \lambda_1 + 2\lambda_2 &\leq 3 \\
    \Delta y_1 &\leq 3\lambda_1 \\
    \Delta y_2 &\leq 3\lambda_2 \label{freq2}\\
    \Delta y_1 &= y - s_1 \\
    \Delta y_2 &= y - s_2 \label{freq3}\\
    s_1 &\leq 3(1-\lambda_1) \\
    s_2 &\leq 3(1-\lambda_2) \\
    z &\geq \Delta y_1 + 2 \Delta y_2 \\
    z &\geq 3x + (y-3) + (\Delta y_1 - 3\lambda_1) + 2(\Delta y_2 - 3\lambda_2)\label{freq4} \\
    z &\leq y + \Delta y_1 + 2\Delta y_2 \\
    z &\leq 3x + (\Delta y_1 - 3\lambda_1) + 2(\Delta y_2 - 3\lambda_2) \\
    \lambda &\in \{0,1\}^2 \\
    \Delta y &\in [0,3]^2 \\
    s &\in [0,3]^2 \\
    (x,y) &\in [0,3] \times [0,3].
\end{align}
\end{subequations}
The feasible point for the relaxation $x = 3$, $y = 3$, $z = 9$, $\lambda = (1,0.5)$, $\Delta y = (3,1.5)$, and $s = (0,1.5)$ is a fractional extreme point, showing that the formulation is not ideal. Indeed, it satisfies at equality the set of linear independent constraints of the relaxation given by $x\leq 3$, $y\leq 3$, $\lambda_1\leq 1$, $\Delta y_1\leq 3$, $s_1\geq 0$, \eqref{freq1}, \eqref{freq2}, \eqref{freq3} and \eqref{freq4}.

\section{Proof of Proposition~\ref{prop:log-bound}} \label{app:log-bound}
First, we present a more general lemma.
\begin{lemma} \label{lemma:log-bound}
    If there do not exist polyhedra $\{Q^i\}_{i=1}^{d'}$ with $d' < d$ and $\bigcup_{i=1}^d P^i = \bigcup_{i=1}^{d'} Q^i$, then any binary MIP formulation for $\bigcup_{i=1}^d P^i$ must have at least $\lceil \log_2(d) \rceil$ binary variables.
\end{lemma}
\proof{\textbf{Proof}}
    Presume that formulation $F$ takes the form \eqref{eqn:generic-MIP-formulation}. For each $h \in \{0,1\}^{n_3}$, consider the preimage $\Pre(h) \defeq \{x \in \bbR^{n_1} : \exists y \in \bbR^{n_2} \text{ s.t. } (x,y,h) \in F\}$; it is clear from the definition of $F$ that this set is polyhedral. Furthermore, we need that $\bigcup_{h \in \{0,1\}^{n_3}} \Pre(h) = \bigcup_{i=1}^d P^i$. If the condition holds, then we have that the cardinality of the index set of the left side ($2^{n_3}$) must be at least $d$, which implies the result.
\Halmos\endproof
We may apply this lemma in the case where our irredundancy assumption holds.
\proof{\textbf{Proof of Proposition~\ref{prop:log-bound}}}
    For each $S \in \scrS$, take $\lambda^S = \frac{1}{|S|}\sum_{v \in S} {\bf e}^v$, where ${\bf e}^v \in \{0,1\}^J$ is the unit vector for component $v$. Assume that there is some $Q^i$ as in Lemma~\ref{lemma:log-bound} such that $\lambda^S, \lambda^{S'} \in Q^i$ for two $S,S' \in \scrS$. By convexity, $\frac{1}{2}(\lambda^S +\lambda^{S'}) \in Q^i$ as well. But this implies that there is a point in $Q^i$ with support over $S \cup S'$, which would violate the irredundancy of $S$ and $S'$, meaning that this cannot yield a formulation for $\bigcup_{i=1}^d P^i$. Therefore, each of the $d$ points $\lambda^S$ must be contained uniquely in some $Q^i$, and we may apply Lemma~\ref{lemma:log-bound} for the result.
\Halmos\endproof

\section{Independence in formulation-induced branching schemes}\label{IBappendix}

As discussed in \cite[Section 3]{Vielma:2009a}, the connection between  MIP formulations  and  branching schemes for CDCs can be used to explain in what sense an independent branching scheme is ``independent.'' As noted in Section~\ref{constraintbranchingsec}, the branching scheme on $\lambda$ induced by formulation \eqref{eqn:multiway-formulation} is precisely the multi-variable branching associated to the corresponding independent branching scheme. In contrast, formulations that are not based on IB schemes (e.g. \eqref{eqn:prop-9.3}) do not necessarily induce a multi-variable branching. However, we can  interpret  the induced effect on the $\lambda$ variables as a multi-way branching scheme that fixes the $\lambda$ variables in a \emph{non-independent} way.

For example, consider SOS$2(5)$ (i.e. $\scrS = \{\{1,2\},\{2,3\},\{3,4\},\{4,5\}\}$). For this particular instance and for  $\{h^S\}_{S \in \scrS}$ given by $h^{\{1,2\}}=(1,1)$, $h^{\{2,3\}}=(1,0)$, $h^{\{3,4\}}=(0,1)$, and $h^{\{4,5\}}=(0,0)$, formulation~\eqref{eqn:prop-9.3} is given by
\begin{subequations}\label{eqn:b-n-b-unbalanced}
\begin{gather}
    \lambda_1 = \gamma_1^{\{1,2\}}, \quad \lambda_2 = \gamma_2^{\{1,2\}} + \gamma_2^{\{2,3\}}, \lambda_3 = \gamma_3^{\{2,3\}} + \gamma_3^{\{3,4\}}, \\ \quad \lambda_4 = \gamma_4^{\{3,4\}} + \gamma_4^{\{4,5\}}, \quad \lambda_5 = \gamma_5^{\{4,5\}} \\
    \gamma_1^{\{1,2\}} + \gamma_2^{\{1,2\}} + \gamma_2^{\{2,3\}} + \gamma_3^{\{2,3\}} + \gamma_3^{\{3,4\}} + \gamma_4^{\{3,4\}} + \gamma_4^{\{4,5\}} + \gamma_5^{\{4,5\}} = 1 \\
    \gamma_1^{\{1,2\}} + \gamma_2^{\{1,2\}} + \gamma_2^{\{2,3\}} + \gamma_3^{\{2,3\}} = z_1 \\
    \gamma_1^{\{1,2\}} + \gamma_2^{\{1,2\}} + \gamma_3^{\{3,4\}} + \gamma_4^{\{3,4\}} = z_2 \\
    \gamma_v^S \geq 0 \quad \forall v, S \\
    (\lambda,z) \in \Delta^5 \times \{0,1\}^2
\end{gather}
\end{subequations}
and a pairwise IB-based formulation (simplified slightly from \eqref{eqn:multiway-formulation}) is
\begin{subequations}\label{eqn:b-n-b-balanced}
\begin{alignat}{2}
    \lambda_1 + \lambda_2 &\leq z_1, \quad\quad \lambda_4 + \lambda_5 &\leq 1-z_1 \\
    \lambda_3 &\leq z_2, \quad\quad \lambda_1 + \lambda_5 &\leq 1-z_2 \\
    &(\lambda,z) \in \Delta^5 \times \{0,1\}^2.
\end{alignat}
\end{subequations}
In Figure~\ref{fig:b-n-b-trees}, we see the first two levels of the branch-and-bound trees for both formulations for the cases where we choose either $z_1$ or $z_2$ for which to branch on first. We observe that, for formulation~\eqref{eqn:b-n-b-unbalanced}, the variables $\lambda_v$ that a given branching decision is able to prove are zero depends on the previous branching decisions in the branch-and-bound tree, while this is not the case for the independent branching formulation~\eqref{eqn:b-n-b-balanced}. For example, if we first branch down on $z_2$ ($z_2 \leq 0$), we are able to prove that $\lambda_1=0$. If we choose instead to branch down on $z_1$ ($z_1 \leq 0$), we are able to prove that $\lambda_1=\lambda_2=0$. However, if we branch down on $z_2$ and then branch down on $z_1$, we prove that $\lambda_1=\lambda_2=0$, but we are also able to prove that $\lambda_3=0$, which we could not prove without the combination of the two branching decisions. Indeed, we see that regardless of the branching decision we make, we will not be able to prove that $\lambda_3=0$ until the second level of the branching tree with formulation \eqref{eqn:b-n-b-unbalanced}.

\begin{figure}[htpb]
    \centering
    \includegraphics[width=0.48\linewidth]{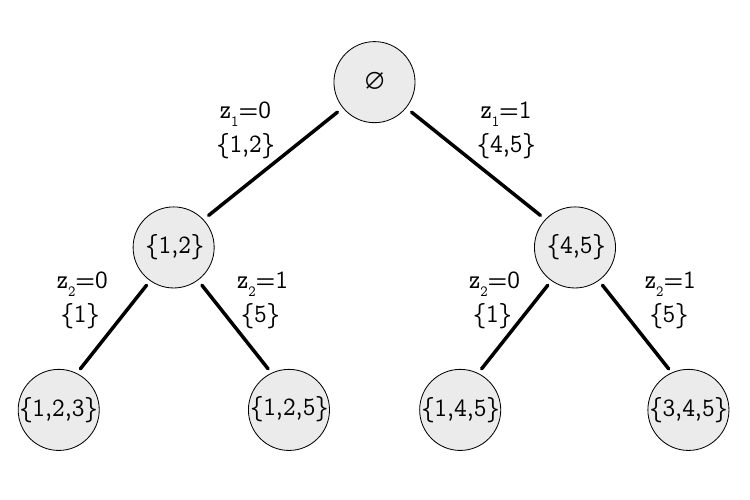} \hspace{1em}
    \includegraphics[width=0.48\linewidth]{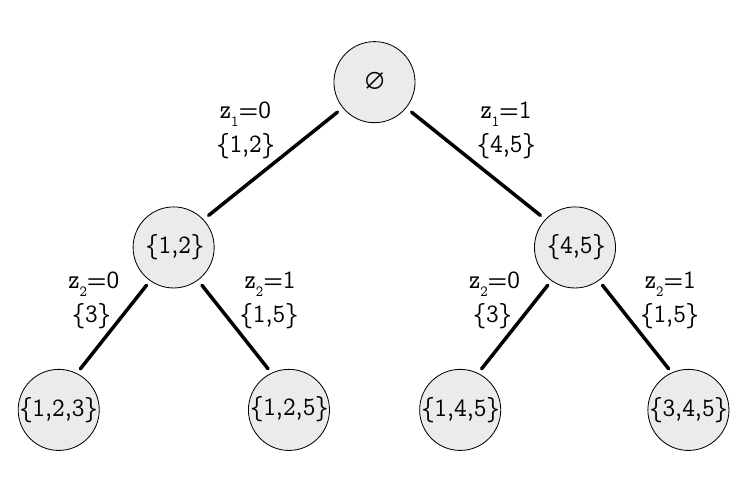} \\
    \includegraphics[width=0.48\linewidth]{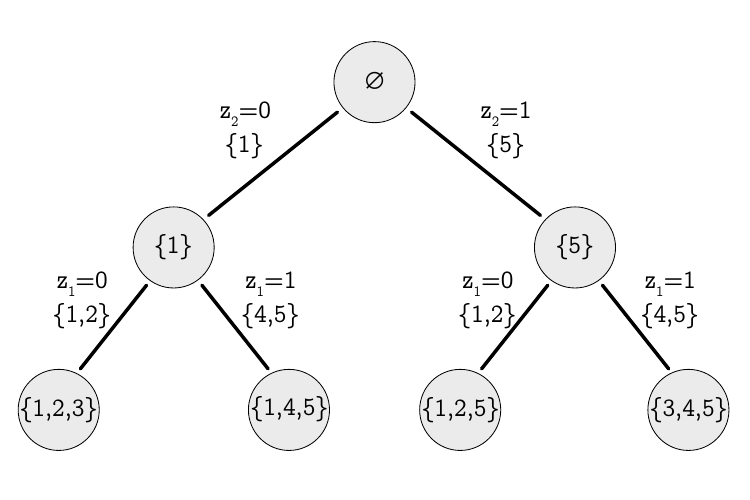} \hspace{1em}
    \includegraphics[width=0.48\linewidth]{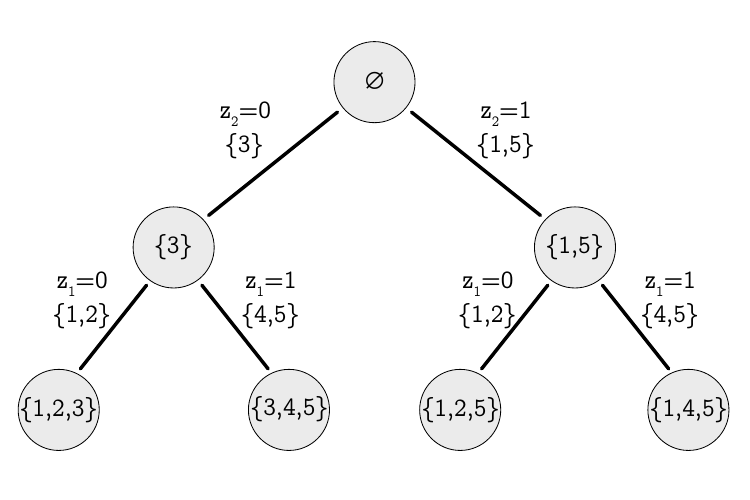}
    \caption{The branch-and-bound trees for \eqref{eqn:b-n-b-unbalanced} (Left) and \eqref{eqn:b-n-b-balanced} (Right), when $z_1$ is first to branch on, and then $z_2$ (Top row), and when $z_2$ is first to branch on, and then $z_1$ (Bottom row). Inside each node is the set $I \subset \llbracket 5 \rrbracket$ of all components $v$ for which the algorithm has been able prove that $\lambda_v = 0$ at this point in the algorithm via branching decisions. The text on the lines show the current branching decision (e.g. $z_2 \geq 1$), and the set of variables $\lambda_v$ for which the (a) subproblem is able to prove that $\lambda_v=0$ \emph{independently} of any other branching decisions (e.g. $z_2 \geq 1$ is the only additional branching constraint added to the original relaxation). This figure is adapted from \cite[Figure 2]{Vielma:2009a}.}
    \label{fig:b-n-b-trees}
\end{figure}

Contrastingly, each branching decision with the independent branching formulation \eqref{eqn:b-n-b-balanced} is able to fix components of $\lambda$ to zero, \emph{independent} of the location in the tree and of the previous branching decisions. For example, branching down or up on $z_2$ is always able to prove either $\lambda_3=0$ or $\lambda_1=\lambda_5=0$, independently. Consequentially, for every component $v \in J$, there exists a branching decision that is able to prove that $\lambda_v=0$ at the first level of the branch-and-bound tree, which is not the case with formulation \eqref{eqn:b-n-b-unbalanced} and $v=3$, as mentioned above. Having this independence  property is a restriction on the branching scheme, but has the potential to simplify branching rules (i.e. choosing which variable $z_i$ to branch on), a notoriously difficult and computationally important part of the algorithmic performance of a MIP solver (see, for example, \cite{Achterberg:2005}). Furthermore, we see that independent branching rules guarantee that the solver can prove any component of $\lambda$ is zero at the very beginning of the tree.

Finally, we note that MIP formulations that are not independent branching formulations can still exhibit the independent branching behavior. For example, if we had selected the encoding $\{h^S\}_{S \in \scrS}$ to be given by given by $h^{\{1,2\}}=(1,1)$, $h^{\{2,3\}}=(1,0)$, $h^{\{3,4\}}=(0,0)$, and $h^{\{4,5\}}=(0,1)$,  formulation~\eqref{eqn:prop-9.3} would satisfy the independent branching property. However, independent branching formulations provide an immediate proof that the property holds, which is not the case for general MIP formulations.

\section{Proposition~\ref{prop:feasibility-IP}} \label{app:feasibility-IP}
\begin{proposition} \label{prop:feasibility-IP}
    A biclique cover of depth $t$ exists for the conflict graph $G^c_\scrS = (J,\bar{E})$ of pairwise IB-representable $\CDC(\scrS)$ if and only if the following admits a feasible solution:
    \begin{subequations} \label{eqn:feasibility-IP}
    \begin{align}
        \left. \begin{array}{l} z^{r,s}_j \leq x^r_j + x^s_j \\
        z^{r,s}_j \leq x^r_j + y^r_j \\
        z^{r,s}_j \leq x^s_j + y^s_j \\
        z^{r,s}_j \leq y^r_j + y^s_j \\
        z^{r,s}_j \geq x^r_j + y^s_j - 1 \\
        z^{r,s}_j \geq x^s_j + y^r_j - 1 \end{array}\right\} &\quad \forall j \in \llbracket t \rrbracket, \forall \{r,s\} \in [J]^2 \label{eqn:feasibility-IP-1} \\
        x^r_j + y^r_j \leq 1 &\quad \forall j \in \llbracket t \rrbracket, \forall r \in J \label{eqn:feasibility-IP-4} \\
        \sum_{j=1}^t z^{r,s}_j \geq 1 &\quad \forall \{r,s\} \in \bar{E} \label{eqn:feasibility-IP-2} \\
        \sum_{j=1}^t z^{r,s}_j = 0 &\quad \forall \{r,s\} \in [J]^2 \backslash \bar{E} \label{eqn:feasibility-IP-3} \\
        x^r \in \{0,1\}^t &\quad \forall r \in J \\
        y^r \in \{0,1\}^t &\quad \forall r \in J \\
        z^{r,s} \in \{0,1\}^t &\quad \forall \{r,s\} \in [J]^2.
    \end{align}
    \end{subequations}
    Moreover, for any feasible solution $(x,y,z)$, a biclique cover for $G_\scrS^c$ is given by $A^j = \{r \in J : x^r_j = 1 \}$ and $B^j = \{r \in J : y^r_j = 1 \}$ for each $j \in \llbracket t \rrbracket$.
\end{proposition}

\proof{\textbf{Proof}}
    The interpretation of the decision variables is:
    \begin{subequations}\label{eqn:feasibility-IP-solution}
    \begin{align}
        x^r_j &= \bbone\left[ r \in A^j \right] \\
        y^r_j &= \bbone\left[ r \in B^j \right] \\
        z^{r,s}_j &= \bbone\left[x^r_j = y^s_j = 1 \text{ or } x^s_j = y^r_j = 1 \right].
    \end{align}
    \end{subequations}
    That is, $z^{r,s}_j = 1$ iff level $i$ separates infeasible edge $\{r,s\} \in \bar{E}$, which is enforced via (\ref{eqn:feasibility-IP-1}-\ref{eqn:feasibility-IP-4}). To show that the existence of a biclique cover implies that \eqref{eqn:feasibility-IP} is feasible, you may consider the proposed solution \eqref{eqn:feasibility-IP-solution} and see that it is feasible for \eqref{eqn:feasibility-IP}.

    To show that a feasible solution maps to a biclique cover, consider some $(x,y,z)$ feasible for \eqref{eqn:feasibility-IP}, and the corresponding sets $A^j = \{r \in J : x^r_j = 1 \}$ and $B^j = \{r \in J : y^r_j = 1 \}$ for each $j \in \llbracket t \rrbracket$. Inequalities \eqref{eqn:feasibility-IP-4} ensure that $A^j \cap B^j = \emptyset$ for each $j \in \llbracket t \rrbracket$. Constraints \eqref{eqn:feasibility-IP-3} ensure that $A^j * B^j \subseteq \bar{E}$ for each $j \in \llbracket t \rrbracket$. Therefore, each $(A^j,B^j)$ is a biclique of $G^c_\scrS$. Furthermore, \eqref{eqn:feasibility-IP-2} ensures that that there is at least one level $j$ that separates each infeasible edge $\{r,s\} \in \bar{E}$. Therefore, $\{(A^j,B^j)\}_{j=1}^t$ is a biclique cover of $G^c_\scrS$.
\Halmos\endproof

\end{APPENDICES}
\end{document}